\documentclass[11pt, letterpaper]{amsart}
\usepackage[colorlinks]{hyperref}
\AtBeginDocument{
  \hypersetup{
    linkcolor=blue,
    citecolor=green,
  }
}
\usepackage{mathscinet}
\usepackage{mathtools}
\usepackage{bm}
\usepackage{cleveref}
\crefname{equation}{}{}

\makeatletter
\let\c@author\relax
\makeatother
\usepackage[backend=bibtex,style=numeric,sorting=nyt,
giveninits=true,isbn=false,url=false]{biblatex}
\renewbibmacro{in:}{\ifentrytype{article}{}{\printtext{\bibstring{in}\intitlepunct}}}
\DeclareFieldFormat*{title}{#1}
\usepackage[mathlines,pagewise]{lineno}
\newcommand*\patchAmsMathEnvironmentForLineno[1]{%
  \expandafter\let\csname old#1\expandafter\endcsname\csname #1\endcsname
  \expandafter\let\csname oldend#1\expandafter\endcsname\csname end#1\endcsname
  \renewenvironment{#1}%
     {\linenomath\csname old#1\endcsname}%
     {\csname oldend#1\endcsname\endlinenomath}}%
\newcommand*\patchBothAmsMathEnvironmentsForLineno[1]{%
  \patchAmsMathEnvironmentForLineno{#1}%
  \patchAmsMathEnvironmentForLineno{#1*}}%
\AtBeginDocument{%
\patchBothAmsMathEnvironmentsForLineno{equation}%
\patchBothAmsMathEnvironmentsForLineno{align}%
\patchBothAmsMathEnvironmentsForLineno{flalign}%
\patchBothAmsMathEnvironmentsForLineno{alignat}%
\patchBothAmsMathEnvironmentsForLineno{gather}%
\patchBothAmsMathEnvironmentsForLineno{multline}%
}

\newtheorem{theorem}{Theorem}[section]
\newtheorem{lemma}[theorem]{Lemma}
\newtheorem{proposition}[theorem]{Proposition}
\newtheorem{corollary}[theorem]{Corollary}

\theoremstyle{definition}
\newtheorem{definition}[theorem]{Definition}

\newtheorem{notation}[theorem]{Notation}

\theoremstyle{remark}
\newtheorem{remark}[theorem]{Remark}

\numberwithin{equation}{section}

\begin{document}

\title[Regularity for locally uniformly elliptic equations]{Pointwise regularity for locally uniformly elliptic equations and applications}

\author{Yuanyuan Lian}
\address{Departamento de An\'{a}lisis Matem\'{a}tico,
Instituto de Matem\'{a}ticas IMAG, Universidad de Granada}
\email{lianyuanyuan.hthk@gmail.com; yuanyuanlian@correo.ugr.es;\newline
 MR ID:\href{https://mathscinet.ams.org/mathscinet/2006/mathscinet/search/author.html?mrauthid=1378049}
{1378049};
ORCID:\href{https://orcid.org/0000-0002-2276-3063}
{0000-0002-2276-3063}}

\author{Kai Zhang}
\address{Departamento de Geometr\'{i}a y Topolog\'{i}a,
Instituto de Matem\'{a}ticas IMAG, Universidad de Granada}
\email{zhangkaizfz@gmail.com; zhangkai@ugr.es;\newline
 MR ID:\href{https://mathscinet.ams.org/mathscinet/2006/mathscinet/search/author.html?mrauthid=1098004}
{1098004};
ORCID:\href{https://orcid.org/0000-0002-1896-3206}
{0000-0002-1896-3206}
}
\thanks{This research has been financially supported by the Project PID2020-118137GB-I00 and PID2020-117868GB-I00 funded by MCIN/AEI /10.13039/501100011033.}

\subjclass[2020]{Primary 35B65, 35D40, 35J60, 35J96, 35J93}

\date{}

\dedicatory{Dedicated to the book Fully Nonlinear Elliptic Equations (by Luis A. Caffarelli and Xavier Cabr\'{e}).}

\keywords{Fully nonlinear equation, regularity theory, viscosity solution, Schauder estimate, Monge-Amp\`{e}re equation}

\begin{abstract}
In this paper, we study the regularity for viscosity solutions of locally uniformly elliptic
equations and obtain a series of interior pointwise $C^{k,\alpha}$ ($k\geq 1$, $0<\alpha<1$) regularity with smallness assumptions on the solution and the right-hand term. As applications, we obtain various interior pointwise regularity for several classical elliptic equations, i.e., the prescribed mean curvature equation, the Monge-Amp\`{e}re equation, the $k$-Hessian equations, the $k$-Hessian quotient equations and the Lagrangian mean curvature equation. Moreover, the smallness assumptions are necessary in most cases (see \Cref{re1.0}, \Cref{re2.1}, \Cref{re4.5}, \Cref{re5.2} and \Cref{re6.2}).
\end{abstract}

\maketitle

\section{Introduction}
In this paper, we study the interior Schauder regularity for viscosity solutions of
\begin{equation}\label{e1.1}
  F(D^2u,Du,u,x)=f~~\mbox{ in}~B_1,
\end{equation}
where $B_1\subset \mathbb{R}^n$ is the unit open ball and $F$ is a \emph{locally uniformly} elliptic operator. Precisely, we use the following notion.
\begin{definition}\label{de10.1}
Let $\rho>0$. The $F:\mathcal{S}^n\times \mathbb{R}^n\times \mathbb{R}\times B_1\to
\mathbb{R}$ is called locally uniformly elliptic with $\rho$ (or $\rho$-uniformly elliptic) if there exist
constants $0<\lambda\leq \Lambda$ such that for any $|M|, |N|,|p|,|s|\leq \rho$ and $x\in B_1$,
\begin{equation}\label{e1.6}
  \mathcal{M}^-(N,\lambda,\Lambda)\leq F(M+N,p,s,x)-F(M,p,s,x)\leq \mathcal{M}^+(N,\lambda,\Lambda),
\end{equation}
where $\mathcal{S}^n$ denotes the set of $n\times n$ symmetric matrices (see \Cref{no1.1}). The $\mathcal{M}^-,
\mathcal{M}^+$ are the usual Pucci's extremal operators:
\begin{equation*}
\mathcal{M}^{+}(M,\lambda,\Lambda)=\Lambda\sum_{\lambda_i>0}\lambda_i
+\lambda\sum_{\lambda_i<0}\lambda_i, \quad
\mathcal{M}^{-}(M,\lambda,\Lambda)=\lambda\sum_{\lambda_i>0}\lambda_i
+\Lambda\sum_{\lambda_i<0}\lambda_i,
\end{equation*}
where $\lambda_i$ are the eigenvalues of $M$.
\end{definition}

\begin{remark}\label{re1.1}
If $\rho=\infty$, we arrive at the definition of uniformly elliptic operators (see \cite[Definition 2.1]{MR1351007}). If $F$ is a smooth operator and
\begin{equation*}
\lambda I\leq D_MF(M,p,s,x)\leq \Lambda I,~\forall ~|M|\leq 2\rho, |p|,|s|\leq \rho, x\in B_1,
\end{equation*}
then \cref{e1.6} holds by the Lagrange mean value theorem.

If $F$ is $\rho$-uniformly elliptic and $u$ is a solution of \cref{e1.1}, then $\tilde u=Ku$ is a solution of
\begin{equation*}
\tilde{F}(D^2\tilde{u},D\tilde{u},\tilde{u},x):=K F(\frac{1}{K}D^2u,\frac{1}{K}Du,\frac{1}{K}u,x)=f
\end{equation*}
and $\tilde{F}$ is $K\rho$-uniformly elliptic. That is, the scaling operation changes the uniform ellipticity. Hence, the scaling operation is restricted in deducing the regularity. This is the main obstacle for developing the regularity theory.
\end{remark}
~\\

The notion of viscosity solution is defined as follows. The main difference from the
classical definition is that the solution must be small and the selection of a test function $\varphi$ is more restricted.
\begin{definition}\label{de10.2}
Let $F$ be $\rho$-uniformly elliptic and $f\in C(B_1)$. We say that $u\in C(B_1)$ is a viscosity
supersolution of \cref{e1.1} if $\|u\|_{L^{\infty}(B_1)}\leq \rho/2$ and for any $x_0\in B_1$ and $\varphi\in C^2(B_1)$ with
\begin{equation}\label{e1.3}
\begin{aligned}
&\|\varphi\|_{C^{1,1}(\bar B_1)}\leq \rho, \quad \varphi(x_0)=u(x_0),\quad \varphi\leq u~~\mbox{ in}~B_1,\\
\end{aligned}
\end{equation}
we have
\begin{equation*}
F(D^2\varphi(x_0),D\varphi(x_0),\varphi(x_0),x_0)\leq f(x_0).
\end{equation*}
Similarly, we can define viscosity subsolution and viscosity solution as usual.
\end{definition}

Since we consider the pointwise regularity in this paper, let us recall the definition of pointwise $C^{k,\alpha}$ (see \cite[Definition 2.2]{lian2020pointwise} for the definition of $C^{k,\omega}$ with a general modulus of continuity).
\begin{definition}\label{d-f}
Let $k\geq 0, 0<\alpha\leq 1$. We say that $f:B_1\to \mathbb{R}$ is $C^{k,\alpha}$ at $0$ or $f\in C^{k,\alpha}(0)$ if there exist a polynomial $P\in \mathcal{P}_k$ (see \Cref{no1.1}) and a constant $K$ such that
\begin{equation}\label{m-holder}
  |f(x)-P(x)|\leq K|x|^{k+\alpha},~\forall~x\in B_{1}.
\end{equation}
We call $P$ the Taylor polynomial of $f$ at $0$ and define
\begin{equation*}
\begin{aligned}
&D^mf(0)= D^mP(0),\quad\|f\|_{C^{k}(0)}= \|P\|=\sum_{m=0}^{k}|D^m P(0)|,\\
&[f]_{C^{k,\alpha}(0)}= \min \left\{K\big | \cref{m-holder} ~\mbox{holds with}~K \right\},~~
\|f\|_{C^{k,\alpha}(0)}= \|f\|_{C^{k}(0)}+[f]_{C^{k,\alpha}(0)}.
\end{aligned}
\end{equation*}

In addition, we say that $f\in C^k(0)$ if there exist a polynomial $P\in \mathcal{P}_k$ such that
\begin{equation}\label{m-holder}
  |f(x)-P(x)|\leq |x|^{k}\omega(|x|),~\forall~x\in B_{1},
\end{equation}
where $\omega$ is a modulus of continuity, i.e., $\omega:\mathbb{R}_+\to \mathbb{R}$ and $\omega(r)\to 0$ as $r\to 0$.
\end{definition}

The regularity theory for uniformly elliptic equations has been well developed. Taking the $C^{k,\alpha}$ regularity for example, we have interior $C^{\alpha}$ regularity (\cite[Chapter 4]{MR1351007}, \cite{MR525227,MR563790}), interior $C^{1,\alpha}$ regularity (\cite[Chapter 8]{MR1351007}), interior $C^{2,\alpha}$ regularity (\cite{MR1391943,MR649348,MR661144,MR688919}, \cite[Chapter 8]{MR1351007}), interior $C^{k,\alpha}$ ($k\geq 3$) regularity (\cite{lian2020pointwise}), boundary $C^{\alpha}$ regularity (\cite{MR4163131}),  boundary $C^{1,\alpha}$ and $C^{2,\alpha}$ regularity (\cite{MR688919,MR4088470,MR3246039}), boundary $C^{k,\alpha}$ regularity (\cite{lian2020pointwise}).

For locally uniformly elliptic equations, the pointwise $C^{2,\alpha}$ regularity for viscosity of \cref{e1.1} was first obtained by Savin \cite{MR2334822}, which was found to have many applications to other problems, e.g.,
\begin{itemize}
  \item Partial regularity for fully nonlinear equations \cite{MR2928094}.
  \item Unique continuation for fully nonlinear equations \cite{MR2875864,MR3311746}.
  \item Regularity for singular nonlinear equations \cite{MR4537427}.
  \item Regularity for the $\sigma_k$-Loewner–Nirenberg problem \cite{MR4611402}.
  \item $W^{2,p}$ estimate for the complex Monge-Amp\`{e}re equation \cite{MR4645029}.
  \item Problems from geometry \cite{MR4718434,MR4460593,MR4688155,MR4153912}.
\end{itemize}
Later, it was extended to parabolic equations by Wang \cite{MR3158522} and nonlocal elliptic equations by Yu \cite{MR3605294}. They all considered homogenous equations (i.e., $f\equiv 0$). In this paper, we study the nonhomogeneous equations and derive a series of pointwise $C^{k,\alpha}$ ($k\geq 1$ and $0<\alpha<1$) regularity.

The motivation of studying locally uniformly elliptic equations is its wide applications. Besides the applications mentioned above, we can obtain various new pointwise regularity for non-uniformly elliptic equations (see \Crefrange{S.MS}{S.SL} for details). Moreover, it indicates that maybe regarding (transforming) non-uniformly elliptic equations as (to) locally uniformly elliptic equations to study the regularity is essential. In another word, based on the regularity theory for locally uniformly elliptic equations, we can obtain regularity for non-uniformly elliptic equations with \emph{necessary} assumptions (see \Cref{re1.0}, \Cref{re2.1}, \Cref{re4.5}, \Cref{re5.2} and \Cref{re6.2}). Note that these applications are just a few of examples. This technique/viewpoint has potential applications to other equations, boundary value problems, parabolic equations, complex equations and equations on Riemannian manifolds etc.

Next, we state our main results on the locally uniformly elliptic equations. For their applications, we state the theorems and give the proofs in \Crefrange{S.MS}{S.SL}. First, we consider the pointwise $C^{1,\alpha}$ regularity for the following equation in a special form:
\begin{equation}\label{e1.1-2}
  F(D^2u,Du,u,x)=A^{ij}(Du,u,x)u_{ij}+B(Du,u,x)=0\quad\mbox{in}~B_1,
\end{equation}
where $F$ is $\rho$-uniformly elliptic, i.e.,
\begin{equation*}
\lambda I\leq A(p,s,x) \leq \Lambda I,~\forall ~|p|,|s|<\rho,~x\in B_1.
\end{equation*}
If $A$ is continuous, we define the modulus of continuity
\begin{equation*}
  \begin{aligned}
\omega_A(r)=&\sup \{|A(p,s,0)-A(q,t,0)|:\\
&|p|,|q|,|s|,|t|\leq \rho,~~|p-q|,|s-t|\leq r \},~\forall ~r>0.
  \end{aligned}
\end{equation*}

\begin{theorem}\label{th1.1}
Let $0<\alpha<1$ and $u\in C(\bar{B}_1)$ be a viscosity solution of \cref{e1.1-2} where $F$ is $\rho$-uniformly elliptic and $A$ is continuous with modulus $\omega_A$. Suppose that for some $K_B,b_0\geq 0$,
\begin{equation}\label{e1.4}
|B(p,s,x)|\leq K_B,~\forall ~|p|,|s|<\rho,~x\in B_1,
\end{equation}
\begin{equation}\label{e1.5}
    -b_0|p-q|\leq B(p,s,x)-B(q,s,x)\leq b_0|p-q|,~\forall ~|p|,|q|,|s|\leq \rho,~x\in B_1
\end{equation}
and
\begin{equation*}
  \|u\|_{L^{\infty}(B_1)}\leq \delta,
\end{equation*}
where $\delta>0$ depends only on $n,\lambda,\Lambda,\rho,b_0,\alpha,\omega_A$ and $K_B$.

Then $u\in C^{1,\alpha}(0)$, i.e., there exists $P\in \mathcal{P}_1$ such that
\begin{equation*}
|u(x)-P(x)|\leq C|x|^{1+\alpha},~\forall ~x\in B_1
\end{equation*}
and
\begin{equation*}
\|P\|\leq \bar C\delta,
\end{equation*}
where $\bar C$ depends only on $n,\lambda,\Lambda,\rho,b_0,\alpha$, and $C$ depends also on
$\omega_A$.
\end{theorem}

\begin{remark}\label{re1.3}
We will apply \Cref{th1.1} to the prescribed mean curvature equation in \Cref{S.MS}. In fact, we can prove interior $C^{1,\alpha}$ regularity under the following more general condition:
\begin{equation*}
F(tM,p,s,x)=tF(M,p,s,x),~\forall ~|M|,|p|,|s|\leq \rho,~x\in B_1,~t>0.
\end{equation*}
Since we do not expect any application, we only consider the simpler equation \cref{e1.1-2} in this paper.
\end{remark}
~\\

%

Next, we consider the $C^{2,\alpha}$ regularity. If $F$ is $\rho$-uniformly elliptic and $D_MF$ is
continuous, we define the modulus of continuity
\begin{equation*}
  \begin{aligned}
\omega_F(r)=&\sup \{|D_MF(M,p,s,0)-D_MF(N,q,t,0)|:\\
&|M|,|N|,|p|,|q|,|s|,|t|\leq \rho,~~|M-N|,|p-q|,|s-t|\leq r \},~\forall ~r>0.
  \end{aligned}
\end{equation*}
We also need the following structure condition. For any $|M|,|p|,|q|,|s|,|t|\leq \rho$ and $x\in B_1$,
\begin{equation}\label{e.st.1}
    -b_0|p-q|-c_0|s-t|\leq F(M,p,s,x)-F(M,q,t,x)\leq b_0|p-q|+c_0|s-t|,
\end{equation}
where $b_0,c_0$ are two nonnegative constants.

The following is the interior pointwise $C^{2,\alpha}$ regularity.
\begin{theorem}\label{th1.2}
Let $0<\alpha<1$ and $u\in C(\bar{B}_1)$ be a viscosity solution of \cref{e1.1} where $F$ is $\rho$-uniformly elliptic, $F(0,0,0,x)\equiv 0$ and $D_MF$ is continuous with modulus
$\omega_F$. Suppose that \cref{e.st.1} holds,
\begin{equation}\label{e.st.2}
|F(M,p,s,x)-F(M,p,s,0)|\leq \delta|x|^{\alpha},~\forall ~|M|,|p|,|s|\leq \rho,~x\in B_1
\end{equation}
and
\begin{equation}\label{e11.4}
  \|u\|_{L^{\infty}(B_1)}\leq \delta,\quad \|f\|_{C^{\alpha}(0)}\leq \delta,
\end{equation}
where $\delta>0$ depends only on $n,\lambda,\Lambda,\rho,b_0,c_0,\alpha$ and $\omega_F$.

Then $u\in C^{2,\alpha}(0)$, i.e., there exists $P\in \mathcal{P}_2$ such that
\begin{equation*}
|u(x)-P(x)|\leq C|x|^{2+\alpha},~\forall ~x\in B_1
\end{equation*}
and
\begin{equation*}
\|P\|\leq \bar C\delta,\quad  F(D^2P,DP(0),P(0),0)=f(0),
\end{equation*}
where $\bar C$ depends only on $n,\lambda,\Lambda,\rho,b_0,c_0,\alpha$, and $C$ depends also on
$\omega_F$.
\end{theorem}
\begin{remark}\label{re1.7}
\Cref{th1.2} was first proved by Savin \cite{MR2334822} with $f\equiv 0$.
\end{remark}
~\\

For higher $C^{k,\alpha}$ ($k\geq 3$) regularity, we have
\begin{theorem}\label{th1.3}
Let $k\geq 3, 0<\alpha<1$ and $u\in C(\bar{B}_1)$ be a viscosity solution of \cref{e1.1}
where $F$ is $\rho$-uniformly elliptic. Let $F_0\in C^{k-1}$ be $\rho$-uniformly elliptic,
$F_0(0,0,0,x)\equiv 0$ and denote
\begin{equation*}
K_F=\|F_{0}\|_{C^{k-1}(\bar{\textbf{B}}_{\rho}\times \bar{B}_1)},\quad
\textbf{B}_{\rho}=\left\{(M,p,s):  |M|,|p|,|s|<\rho\right\}.
\end{equation*}
Suppose that \cref{e.st.1} holds,
\begin{equation}\label{e1.2}
|F(M,p,s,x)-F_0(M,p,s,x)|\leq \delta |x|^{k-2+\alpha},~\forall ~|M|,|p|,|s|\leq \rho,~x\in B_1
\end{equation}
and
\begin{equation*}
  \|u\|_{L^{\infty}(B_1)}\leq \delta,\quad \|f\|_{C^{k-2+\alpha}(0)}\leq \delta,
\end{equation*}
where $\delta>0$ depends only on $n,\lambda,\Lambda,\rho,b_0,c_0,k,\alpha$ and $K_F$.

Then $u\in C^{k,\alpha}(0)$, i.e., there exists $P\in \mathcal{P}_{k}$ such that
\begin{equation*}
|u(x)-P(x)|\leq C|x|^{k+\alpha},~\forall ~x\in B_1
\end{equation*}
and
\begin{equation*}
\|P\|\leq C\delta,\quad  |F_0(D^2P(x),DP(x),P(x),x)-P_f(x)|\leq C|x|^{k-1},~\forall ~x\in B_1,
\end{equation*}
where $C$ depends only on $n,\lambda,\Lambda,\rho,b_0,c_0,k,\alpha$ and $K_F$.
\end{theorem}

\begin{remark}\label{re1.8}
This higher pointwise regularity is new. Even for the uniformly elliptic equations, it was proved in \cite{lian2020pointwise} recently. In fact, the proof of \Cref{th1.3} is inspired by \cite{lian2020pointwise}. In addition, by similar arguments in \cite{lian2020pointwise}, we can obtain the pointwise $C^{k}(0)$ and $C^{k\mathrm{ln}L}(0)$ ($k\geq 2$) regularity as well.
\end{remark}
~\\

In fact, in above theorem, we only need that $\|f||_{C^{\alpha}(0)}$ is small.
\begin{corollary}\label{co1.0}
Let $k\geq 3, 0<\alpha<1$ and $u\in C(\bar{B}_1)$ be a viscosity solution of \cref{e1.1}
where $F$ is $\rho$-uniformly elliptic. Let $F_0\in C^{k-1}$ be $\rho$-uniformly elliptic and
$F_0(0,0,0,x)\equiv 0$. Suppose that \cref{e.st.1} holds,
\begin{equation*}
|F(M,p,s,x)-F_0(M,p,s,x)|\leq \delta|x|^{k-2+\alpha},~\forall ~|M|,|p|,|s|\leq \rho,~x\in B_1
\end{equation*}
and
\begin{equation*}
  \|u\|_{L^{\infty}(B_1)}\leq \delta,\quad f\in C^{k-2+\alpha}(0),\quad \|f\|_{C^{\alpha}(0)}\leq \delta,
\end{equation*}
where $\delta>0$ depends only on $n,\lambda,\Lambda,\rho,b_0,c_0,\alpha$ and $K_F$.

Then $u\in C^{k,\alpha}(0)$. That is, there exists $P\in \mathcal{P}_{k}$ such that
\begin{equation*}
|u(x)-P(x)|\leq C|x|^{k+\alpha},~\forall ~x\in B_1
\end{equation*}
and
\begin{equation*}
\|P\|\leq C\delta,\quad  |F_0(D^2P(x),DP(x),P(x),x)-P_f(x)|\leq C|x|^{k-1},~\forall ~x\in B_1,
\end{equation*}
where $C$ depends only on $n,\lambda,\Lambda,\rho,b_0,c_0,k,\alpha,K_F$ and $\|f\|_{C^{k-2+\alpha}(0)}$.
\end{corollary}

As a special case, we have
\begin{corollary}\label{co1.1}
Let $k\geq 3, 0<\alpha<1$ and $u\in C(\bar{B}_1)$ be a viscosity solution of \cref{e1.1}
where $F\in C^{k-1}(\bar{\textbf{B}}_{\rho}\times \bar{B}_1)$ is $\rho$-uniformly elliptic and $F(0,0,0,x)\equiv 0$. Suppose that
\begin{equation*}
  \|u\|_{L^{\infty}(B_1)}\leq \delta,\quad f\in C^{k-2+\alpha}(0),\quad \|f\|_{C^{\alpha}(0)}\leq \delta,
\end{equation*}
where $\delta>0$ depends only on $n,\lambda,\Lambda,\rho,b_0,c_0,\alpha$ and $K_F$. Then $u\in C^{k,\alpha}(0)$.
\end{corollary}

\begin{remark}\label{re1.2}
In this paper, a constant is called universal if it depends only on $n,\lambda,\Lambda,\rho, b_0$ and $c_0$.
\end{remark}
~\\

We use the perturbation technique as in \cite{MR2334822} to prove above theorems. The idea is the following. Take $C^{2,\alpha}$ regularity for example. If $u$ and $f$ are small and $F$ is smooth, the equation is close the Laplace equation (by the compactness method, see \cref{e11.3}). Hence, there exists $P_1\in \mathcal{P}_2$ such that
\begin{equation*}
\|u-P_1\|_{L^{\infty }(B _{\eta})}\leq \eta^{2+\alpha},
\end{equation*}
where $0<\eta<1$. By scaling, we have
\begin{equation*}
\|u-P_m\|_{L^{\infty }(B _{\eta^{m}})}\leq \eta ^{m(2+\alpha)},~~m\geq 1,
\end{equation*}
which means the $C^{2,\alpha}$ regularity. During the scaling, $\|P_m\|$ are kept small such that the scaled operators are always $\rho/2$-uniformly elliptic.

The main obstacle is to show the compactness of solutions. Savin \cite{MR2334822} proved a Harnack inequality by the technique of sliding paraboloids and then the H\"{o}lder regularity follows. In this paper, we follows the idea in \cite{MR1351007} to prove a weak Harnack inequality, which leads to the H\"{o}lder regularity as well. In fact, we just repeat the argument in \cite{MR1351007}.

The idea ``smallness implies regularity'' has been found and used many years ago (e.g. \cite[Proposition 2]{MR673830}, \cite{MR0179651}, \cite[(2.20)]{MR481490}, \cite[Theorem 1.3]{MR1139064}). Since the equation is regarded as a perturbation of the Laplace equation, we can obtain $C^{k,\alpha}$ ($k\geq 1$) regularity for \emph{any} $0<\alpha<1$. On the contrast, we can only obtain $C^{k,\alpha}$ regularity for \emph{some} $0<\alpha<\bar{\alpha}$ for a general fully nonlinear elliptic equation, where $0<\bar{\alpha}<1$ is a universal constant (see \cite[Chapter 8]{MR1351007}). In addition, the proofs of above theorems are relatively simpler than that for uniformly elliptic equations without the smallness assumptions (compare with \cite{lian2020pointwise}).

Note that even for uniformly elliptic equations, one usually start the proof by assuming that $f$ is small (see \cite[P. 75, Proof of Theorem 8.1]{MR1351007}). Then a scaling argument can transform a general $f$ to a small $f$. However, for locally uniformly elliptic equations, as pointed out in \Cref{re1.1}, the scaling argument is restricted. Hence, we have to make the assumption in the theorem that $f$ is small.

The paper is organized as follows. We first give the applications of above theorems in \Crefrange{S.MS}{S.SL}. Precisely, we shall prove a series of interior pointwise $C^{k,\alpha}$ regularity for the prescribed mean curvature equation in \Cref{S.MS}, the Monge-Amp\`{e}re equation in \Cref{S.MA}, the $k$-Hessian equation in \Cref{S.H}, the $k$-Hessian quotient equation in \Cref{S.HQ} and the Lagrangian mean curvature equation in \Cref{S.SL} respectively.

The proofs of above theorems are postponed to \Crefrange{S10}{S4}. In \Cref{S10}, we prepare some preliminaries, such as the Alexandrov-Bakel'man-Pucci maximum principle etc. We prove $C^{1,\alpha}$ regularity in \Cref{S3} and \Cref{S11} is devoted the $C^{2,\alpha}$ regularity. Finally, we give the proof of $C^{k,\alpha}$ ($k\geq 3$) regularity in \Cref{S4}.

\begin{notation}\label{no1.1}
\begin{enumerate}~~\\
\item $\mathbb{R}^n$: the $n$-dimensional Euclidean space; $\mathcal{S}^n$: the set of $n\times n$
    symmetric matrices with the standard order.
\item $\{e_i\}^{n}_{i=1}$: the standard basis of $\mathbb{R}^n$, i.e.,
    $e_i=(0,...0,\underset{i^{th}}{1},0,...0)$.
\item $x=(x_1,...,x_n)=(x',x_n)\in \mathbb{R}^n$.
\item $|x|=\left(\sum_{i=1}^{n} x_i^2\right)^{1/2}$ for $x\in \mathbb{R}^n$; $|M|=$ the spectrum
    radius of $M\in \mathcal{S}^n$.
\item $I$: the unit matrix in $\mathcal{S}^n$.
\item $a^+=\max(a,0)$, $a^-=\max(-a,0)$.
\item $B_r(x_0)=B(x_0,r)=\{x\in \mathbb{R}^{n}:  |x-x_0|<r\}$, $B_r=B_r(0)$,
\item $\Omega^c$: the complement of $\Omega$;\quad $\bar \Omega $: the closure of $\Omega$, where $ \Omega\subset \mathbb{R}^n$.
\item $\mathrm{diam}(\Omega)$: the diameter of $\Omega$.
\item Given $\varphi:\mathbb{R}^n\to \mathbb{R}$, define $\varphi _i=\partial \varphi/\partial x _{i}$, $\varphi
    _{ij}=\partial ^{2}\varphi/\partial x_{i}\partial x_{j}$ and we also use similar notations for higher
    order derivatives.
\item $D^0\varphi=\varphi$, $D \varphi= (\varphi_1 ,...,\varphi_{n} )$ and $D^2 \varphi= \left(\varphi _{ij}\right)_{n\times n}$ etc.
\item We also use the standard multi-index notation. Let $\sigma=(\sigma_1,...,\sigma_n)\in \mathbb{N}^n$, i.e., each component $\sigma_i$ is a nonnegative integer. Define
\begin{equation*}
|\sigma|= \sum_{i=1}^{n}\sigma_i,\quad\sigma!= \prod_{i=1}^{n}(\sigma_i!),\quad
x^{\sigma}= \prod_{i=1}^{n} x_i^{\sigma_i},\quad
D^{\sigma}\varphi = \frac{\partial^{|\sigma|} \varphi }{\partial x_1^{\sigma_1}\cdots \partial x_n^{\sigma_n}}.
\end{equation*}
\item $|D^k\varphi |= \left(\sum_{|\sigma|=k}|D^{\sigma}\varphi| ^2\right)^{1/2}$ for $k\geq 0$.
\item Given $F\colon \mathcal{S}^n\times \mathbb{R}^n\times \mathbb{R}\times \Omega\to \mathbb{R}$, define
\begin{equation*}
F_{M_{ij}}= \frac{\partial F}{\partial M_{ij}},\quad F_{p_i}= \frac{\partial F}{\partial p_{i}},\quad F_{s}= \frac{\partial F}{\partial s},\quad F_{x_i}= \frac{\partial F}{\partial x_{i}},
\end{equation*}
where $1\leq i,j\leq n$. Moreover, let $\xi\in \mathbb{N}^{n\times n}$ denote the matrix-valued multi-index. Then define
\begin{equation*}
D^{\xi}_MF= \frac{\partial^{|\xi|} F}{\partial M_{ij}^{\xi_{ij}}},\quad
D^k_MF=\left\{\frac{\partial^k F}{\partial M^{\xi}}\colon  |\xi|=k\right\},\quad
|D^k_MF|= \left(\sum_{|\xi|=k}\left|\frac{\partial^k F}{\partial M^{\xi}}\right| ^2\right)^{1/2}.
\end{equation*}
Similarly, we can define $D^k_{p}F$, $D^k_{s}F$, $D^k_{x}F$ and their norms etc.
\item $\mathcal{P}_k (k\geq 0):$ the set of polynomials of degree less than or equal to $k$. That is, any
    $P\in \mathcal{P}_k$ can be written as
\begin{equation*}
P(x)=\sum_{|\sigma|\leq k}\frac{a_{\sigma}}{\sigma!}x^{\sigma}
\end{equation*}
where $a_{\sigma}$ are constants. Define for $r>0$
\begin{equation*}
\|P\|_r= \sum_{|\sigma|\leq k}r^{|\sigma|}|a_{\sigma}|, \quad \|P\|=\|P\|_1= \sum_{|\sigma|\leq k}|a_{\sigma}|
\end{equation*}
\item $\mathcal{HP}_k (k\geq 0):$ the set of homogeneous polynomials of degree $k$. That is, any $P\in
    \mathcal{HP}_k$ can be written as
\begin{equation*}
P(x)=\sum_{|\sigma|= k}\frac{a_{\sigma}}{\sigma!}x^{\sigma}.
\end{equation*}
\end{enumerate}
\end{notation}

\section{Prescribed mean curvature equation}\label{S.MS}
In the following sections, we give the applications of the regularity of locally uniformly elliptic equations to several classical non-uniformly elliptic equations.

Since non-uniformly elliptic equations are much more difficult than uniformly elliptic ones, it is natural to
study them by assuming the data are good enough and then obtain a priori estimates, existence of smooth solutions and Liouville type theorems etc. This classical strategy has been used widely. For example,
\begin{itemize}
\item For the prescribed mean curvature equation, see \cite{MR1081183,MR0248647,MR66533,
MR222467,MR843597,MR265745,MR282058,MR0296832,MR1617971}.
  \item For the Monge-Amp\`{e}re equation, see \cite{MR739925,MR106487, MR437805,
      chu2023liouville,MR718679,MR815272,MR579472,MR62326,MR691980,MR2415390}.
  \item For the $k$-Hessian equations, see \cite{MR806416,MR3961212,MR815272,MR2038151,
      MR585774,MR992979,MR1368245,MR1275451,MR2500526}.
  \item For the $k$-Hessian quotient equations, see \cite{MR1963687,MR3342406,MR4525733,
      MR3635978,MR4217915,dai2023interior,MR4439494,lu2023interior,lu2024interior,MR1368245}
  \item For the Lagrangian mean curvature equation, see \cite{MR4314140,BMS2022,MR806416,MR2492708,MR4578557,MR3188067,
      MR2511754,MR2666907,MR1930884,Zhou2023_2}
\end{itemize}

In addition, there exist Pogorelov's type interior $C^{1,1}$ estimates for some equations. The smoothness requirements of the boundary can be relaxed. As a compensation, the boundary value must be an affine function or an admissible function in some sense. For example,
\begin{itemize}
  \item For the Monge-Amp\`{e}re equation, see \cite{MR293227},\cite[P. 73, Chapter 5.3]{MR478079}.
  \item For the $k$-Hessian equations, see \cite{MR1835381}.
  \item For the $k$-Hessian quotient equations, see \cite{MR4525733,MR4217915}.
\end{itemize}
For the prescribed mean curvature equation and the Lagrangian mean curvature equation, we have stronger pure interior $C^{0,1}$ and $C^{1,1}$ estimates (i.e., estimates independent of the boundary information).

In another direction, it is also natural to consider the regularity of solutions. One may introduce weak solutions in some sense and then prove the existence and regularity. In this respect, there are following examples:
\begin{itemize}
\item For the prescribed mean curvature equation, see \cite{MR2911884,MR437925,MR336532,MR487722,leonardi2023,MR352682,MR1318169}.
  \item For the Monge-Amp\`{e}re equation, see
\cite{MR2318817,MR1038359,MR1038360,MR1127042,MR1124980,MR1426885,MR4334976,MR3084340,
MR3032325,MR111943,MR2983006}.
  \item For the $k$-Hessian equations, see \cite{MR2133413,MR1466315,MR1634570,MR1726702,MR1923626,MR1777141,MR1840289}.
  \item For the Lagrangian mean curvature equation, see \cite{BS2020,MR4649187,MR4655360,
      MR1808031}.
\end{itemize}
As for the Schauder's type regularity, say $C^{2,\alpha}$ regularity under the assumption $f\in C^{\alpha}$, there are few results in this respect (see \cite{MR4093797} for the $2$-Hessian equation in dimension $3$ and \cite{BS2020} for the Lagrangian mean curvature equation) besides the Monge-Amp\`{e}re equation (see \cite{MR1038360,MR2338423,MR425854,MR643043,MR926846,
WXJ-1992} for the interior regularity; see \cite{MR2983006} and \cite{MR2415390} for the boundary regularity).

In the following sections, based on the regularity theory for locally uniformly elliptic equations, we will develop pointwise $C^{k,\alpha}$ ($k\geq 1$, $0<\alpha<1$) regularity for the prescribed mean curvature equation (\Cref{S.MS}), the Monge-Amp\`{e}re equation (\Cref{S.MA}), the $k$-Hessian equation (\Cref{S.H}), the $k$-Hessian quotient equation (\Cref{S.HQ}) and the Lagrangian mean curvature equation (\Cref{S.SL}).

We use the notion of viscosity solution for these equations as well. For the prescribed mean curvature equation and the Lagrangian mean curvature equation, the definition of viscosity solution is exactly the same as the usual (see \cite[Chapter 2]{MR1351007}). For the Monge-Amp\`{e}re equation, the test function $\varphi$ should be a convex function (see \cite[Definition 1.3.1]{MR1829162}). The $\varphi$ should be $k$-admissible if we define a viscosity solution for a $k$-Hessian equation or $k$-Hessian quotient equation (see \cite[Section 2]{MR1089043}, \cite[Definition 1.1]{li2020hessian}).

Note that above definitions of viscosity solution are different from \Cref{de10.2}. In fact, in the following sections, we will regard (transform) an equation as (to) a locally uniformly elliptic equation. Then the (transformed) solution will be a viscosity solution in the sense of \Cref{de10.2}.

%

We first give a general result:
\begin{theorem}\label{th2.1}
Let $k\geq 2$, $0<\alpha<1$ and $u$ be a viscosity solution of
\begin{equation*}
F(D^2u,Du,u,x)=f~~\mbox{ in}~B_1,
\end{equation*}
where $F$ is smooth. Suppose that $u\in C^{2}(0),f\in C^{k-2,\alpha}(0)$ and
\begin{equation}\label{e2.0}
A:=D_MF(D^2u(0),Du(0),u(0),0)>0.
\end{equation}
Then $u\in C^{k,\alpha}(0)$.
\end{theorem}
\proof By $u\in C^{2}(0)$, there exist $P\in\mathcal{P}_2$ and a modulus of continuity $\omega$ such that
\begin{equation*}
  |u(x)-P(x)|\leq \omega(|x|)|x|^2,~\forall ~x\in B_1.
\end{equation*}
Since $u$ is a viscosity solution,
\begin{equation*}
F(D^2u(0), Du(0), u(0),0)=F(D^2P(0), DP(0), P(0),0)=f(0).
\end{equation*}

For $r>0$, let
\begin{equation*}
y=\frac{x}{r}, \quad  \tilde{u}(y)=\frac{u(x)-P(x)}{r^2}, \quad \tilde{f}(y)=f(x)-F(D^2P,DP(x),P(x),x).
\end{equation*}
Then $\tilde{u}$ is a solution of
\begin{equation}\label{e2.0-2}
\tilde F(D^2\tilde{u},D\tilde{u},\tilde{u},y)=\tilde{f}\quad\mbox{in}~B_1,
\end{equation}
where
\begin{equation*}
\tilde F(M,p,s,y):=F(M+D^2P, rp+DP(x), r^2s+P(x),x)-F(D^2P,DP(x),P(x),x).
\end{equation*}
Then $\tilde{F}(0,0,0,y)\equiv 0$ and by \cref{e2.0},
\begin{equation*}
D_M\tilde{F}(0,0,0,0)=A>0.
\end{equation*}
Since $\tilde F$ is smooth, $\tilde{F}$ is $\rho$-uniformly elliptic with some $0<\lambda\leq \Lambda$ and these three constants depends only on $n,A,\|P\|$ and $\|F\|_{C^{1,1}}$.

In addition, by the definition of $\tilde{u}$ and $\tilde{f}$,
\begin{equation*}
\|\tilde{u}\|_{L^{\infty}(B_1)}\leq \omega(r), \quad  \tilde{f}(0)=0, \quad [f]_{C^{\alpha}(0)}\leq Cr^{\alpha}.
\end{equation*}
Thus, we can choose $r$ small enough such that
\begin{equation*}
\|\tilde u\|_{L^{\infty}(B_1)}\leq\delta, \quad \|\tilde f\|_{C^{\alpha}(0)}\leq \delta,
\end{equation*}
where $\delta$ is small enough such that $\tilde{u}$ is a viscosity solution of \cref{e2.0-2} in the sense of
\Cref{de10.2} and we can apply \Cref{co1.1} to \cref{e2.0-2}. Therefore, $\tilde{u}\in C^{k,\alpha}(0)$ and hence $u\in C^{k,\alpha}(0)$.~\qed~\\

\begin{remark}\label{re2.0}
If $u$ is an appropriate viscosity solution, \cref{e2.0} holds automatically. At least, it holds for all equations treated in the following sections. Hence, in general, if $u\in C^{2}(0)$ and $f\in C^{k-2,\alpha}(0)$, then $u\in C^{k,\alpha}(0)$. In conclusion, we have one rough but interesting assertion: For any elliptic equation with a smooth operator, the pointwise $C^{2,\alpha}$ regularity holds almost everywhere if $u$ is a convex viscosity solution and $f\in C^{\alpha}$ ($0<\alpha<1$).
\end{remark}
~\\

Now, we consider the prescribed mean curvature equation.
\begin{theorem}\label{th.MS}
Let $0<\alpha<1$ and $u\in C(\bar{B}_1)$ be a viscosity solution of
\begin{equation}\label{e.MS}
  \mathrm{div} \left(\frac{Du}{\sqrt{1+|Du|^2}}\right)=
\frac{1}{\sqrt{1+|Du|^2}}\left(\delta^{ij}-\frac{u_iu_j}{1+|Du|^2}\right)u_{ij}
=f~\mbox{ in}~B_1,
\end{equation}
where $f\in L^{\infty}(B_1)$. Then $u\in C^{1,\alpha}(0)$ provided one of the following conditions holds:\\
(i) there exists $P\in \mathcal{P}_1$ such that
\begin{equation*}
\|u-P\|_{L^{\infty}(B_1)}\leq \delta,
\end{equation*}
where $0<\delta<1$ depends only on $n,\alpha,|DP|$ and $\|f\|_{L^{\infty}(B_1)}$.\\
(ii) there exists a constant $|K|<n-1$ such that
\begin{equation*}
  \|f-K\|_{L^{\infty}(B_1)}\leq \delta,
\end{equation*}
where $0<\delta<1$ depends only on $n,\alpha,K$ and $\|u\|_{L^{\infty}(B_1)}$.\\
(iii) $u\in C^{0,1}(0)$.\\
(iv) there exists $P\in \mathcal{P}_1$ such that
\begin{equation*}
u=P\quad\mbox{on}~~\partial B_{\delta},
\end{equation*}
where $\delta>0$ depends only on $n,\alpha,|DP|$ and $\|f\|_{L^{\infty}(B_1)}$.
\end{theorem}

For higher regularity, we have
\begin{theorem}\label{th.MS-2}
Let $k\geq 2, 0<\alpha<1$, $u\in C(\bar{B}_1)$ be a viscosity solution of \cref{e.MS} and $f\in C^{k-2,\alpha}(0)$. Then $u\in C^{k,\alpha}(0)$ provided one of the four conditions in \Cref{th.MS} holds.
\end{theorem}

\begin{remark}\label{re2.2}
The mean curvature operator is given by
\begin{equation*}
F(M,p)= \frac{1}{\sqrt{1+|p|^2}}\left(\delta^{ij}-\frac{p_ip_j}{1+|p|^2}\right)M_{ij}.
\end{equation*}
Hence, $F$ is smooth and it is easy to verify (cf. \Cref{re1.1}) that $F$ is $1$-uniformly elliptic with $\lambda=\sqrt{2}/4$ and $\Lambda=\sqrt{2}/2$. Moreover, for any $P\in\mathcal{P}_1$, define
\begin{equation*}
G(M,p)=F(M,p+DP).
\end{equation*}
Then $G$ is smooth and $1$-uniformly elliptic with $\tilde{\lambda}$ and $\tilde{\Lambda}$, which depend only on $|DP|$. Hence, the regularity theory for locally uniformly elliptic equations are applicable to \cref{e.MS}.
\end{remark}

\begin{remark}\label{re1.0}
To investigate the regularity of solutions of \cref{e.MS}, one usually assumes $f\in C^{0,1}$ at least (e.g.
\cite[Theorem 4.2]{MR2836589},\cite{MR437925},\cite[Chapter 16]{MR1814364},\cite{MR487722}). In fact, if $f\in C^{0,1}$, one can prove the interior gradient bound, which was first obtained for the minimal surface equation by Bombieri, De Giorgi and Miranda \cite{MR0248647} (see also
\cite{MR1081183,MR66533,MR843597,MR265745,MR0296832,MR1617971}). Then the equation
becomes uniformly elliptic and the regularity follows.

On the other hand, if $f\notin C^{0,1}$, we cannot obtain the $C^{1,\alpha}$ regularity in general.
Consider the following counterexample borrowed from \cite[Section 8]{MR2911884}:
\begin{equation*}
u(x)=\left\{
  \begin{aligned}
&-(1-|x|)^{\theta},&&~\mbox{ if}~0\leq r\leq 1;\\
&(|x|-1)^{\theta},&&~\mbox{ if}~1<r<2,\\
  \end{aligned}
  \right.
\end{equation*}
where $0<\theta<1/2$. It can be checked directly that $u$ is a viscosity solution of \cref{e.MS} in $B_2$
with
\begin{equation*}
f\in C^{1-2\theta}(\bar{B}_2), \quad [f]_{C^{1-2\theta}(x_0)}\leq C\theta^{-2},~\forall ~x_0\in
\partial B_1.
\end{equation*}
However, we have only $u\in C^{\theta}$ at $\partial B_1$. Hence, the ``smallness'' assumptions in \Cref{th.MS} cannot be removed.

Moreover, if $\theta$ is smaller, $f$ is smoother but $[f]_{C^{1-2\theta}}$ is bigger. Correspondingly, $u$
has lower regularity. This phenomenon indicates that the smallness is more important than the smoothness
for the regularity in this case (i.e., non-uniformly elliptic equations with lower regularity on $f$).

Above observation may imply that regarding the prescribed mean curvature equation as
a locally uniformly elliptic equation to study the regularity is essential (see also \Cref{re2.1}, \Cref{re5.2} and \Cref{re6.2}).
\end{remark}

\begin{remark}\label{re2.4}
Roughly speaking, \Cref{th.MS} states that the interior $C^{1,\alpha}$ regularity holds for the prescribed mean curvature equation under some ``smallness'' assumption (except (iii)).

As pointed out in \Cref{re1.0}, to obtain $C^{1,\alpha}$ regularity, one usually prove the interior gradient bound first and then the equation becomes uniformly elliptic. The assumption (iii) can be understood in some sense that the equation is \emph{pointwise} uniformly elliptic at $0$. Then we obtain the pointwise $C^{1,\alpha}(0)$ regularity.
\end{remark}

\begin{remark}\label{re2.10}
The $C^{1,\alpha}$ regularity for the minimal surface equation under the assumption that $u$ is small, i.e.,
\begin{equation*}
\|u\|_{L^{\infty}(B_1)}\leq \delta,~f\equiv 0 \Longrightarrow u\in C^{1,\alpha}(\bar{B}_{1/2})
\end{equation*}
has been proved by De Giorgi \cite {MR0179651} as a special case (see also \cite[Chapters 6-8]{MR775682}, \cite[P. 676]{MR2334822} and \cite[P. 42]{MR2480601}). \Cref{th.MS} extends this result to
\begin{equation*}
\|u\|_{L^{\infty}(B_1)}\leq \delta,~f\in L^{\infty}(B_1)  \Longrightarrow u\in C^{1,\alpha}(\bar{B}_{1/2}).
\end{equation*}
\end{remark}
~\\

We first prove a lemma.
\begin{lemma}\label{le2.4}
Let $u$ be a viscosity solution of \cref{e.MS}. Then for any $\delta_1>0$, there exists $\delta>0$
depending only on $n$, $K$ and $\|u\|_{L^{\infty}(B_1)}$ such that if
\begin{equation*}
\|f-K\|_{L^{\infty}(B_1)}\leq \delta, \quad |K|<n-1,
\end{equation*}
we have for some $P\in \mathcal{P}_1$
\begin{equation*}
\|u-P\|_{L^{\infty}(B_{r})}\leq \delta_1 r
\end{equation*}
and
\begin{equation*}
|DP|\leq C,
\end{equation*}
where $0<r<1/2$ and $C$ depend only on $n,K$ and $\|u\|_{L^{\infty}(B_1)}$.
\end{lemma}
\proof Since $|K|<n-1$, there exist two solutions $v^{\pm}\in C^{\infty}(B_1)\cap C(\bar{B}_1)\cap W^{1,1}(B_1)$ of the following prescribed mean curvature equations (see \cite[Theorem 16.11]{MR1814364} and \cite[Theorem 1]{MR437925})
\begin{equation*}
  \left\{
  \begin{aligned}
\mathrm{div} A(Dv^{\pm})=&K\pm2\delta &&\quad\mbox{in}~~B_1;\\
v^{\pm}=&u &&\quad\mbox{on}~~\partial B_1,
  \end{aligned}
  \right.
\end{equation*}
where $\delta$ is taken small such that $|K\pm2\delta|\leq n-1$ and
\begin{equation}\label{e2.18}
A(p):=\frac{p}{\sqrt{1+|p|^2}},~\forall ~p\in \mathbb{R}^n.
\end{equation}
Since $v^{\pm}$ are smooth, by the definition of viscosity solution,
\begin{equation}\label{e2.9}
v^+\leq u\leq v^-\quad\mbox{in}~~B_1.
\end{equation}

We claim that for any $\delta_2>0$, if $\delta$ small enough (depending only on $\delta_2,n,K,\|u\|_{L^{\infty}(B_1)}$),
\begin{equation}\label{e2.3}
\|v^{\pm}-v\|_{L^{\infty}(B_{1/2})}\leq \delta_2,
\end{equation}
where $v$ is the solution of
\begin{equation*}
  \left\{
  \begin{aligned}
\mathrm{div} A(Dv)=&K\quad\mbox{in}~~B_1;\\
v=&u\quad\mbox{on}~~\partial B_1.
  \end{aligned}
  \right.
\end{equation*}

We prove the claim by contradiction and the proof is inspired by \cite{MR0301343}. Suppose not. Then there exist $\delta_2$ and sequences of $u_m,v^+_m,v_m$ satisfying $\|u_m\|_{L^{\infty}(B_1)}\leq K$,
\begin{equation}\label{e2.15}
  \left\{
  \begin{aligned}
\mathrm{div}  A(Dv_m^+)
=&K+\frac{1}{m} &&~~\mbox{in}~~B_1;\\
v_m^{+}=&u_m &&~~\mbox{on}~~\partial B_1,
  \end{aligned}
  \right.~~
   \left\{
  \begin{aligned}
\mathrm{div}  A(Dv_m)=&K &&~~\mbox{in}~~B_1;\\
v_m=&u_m &&~~\mbox{on}~~\partial B_1
  \end{aligned}
  \right.
\end{equation}
and
\begin{equation}\label{e2.10}
\|v_m^{+}-v_m\|_{L^{\infty}(B_{1/2})}> \delta_2.
\end{equation}

By the interior derivatives estimates (see \cite[Corollary 16.7]{MR1814364}),
\begin{equation*}
\|v_m^+\|_{C^{1,1}(\bar\Omega')},\quad \|v_m\|_{C^{1,1}(\bar\Omega')}\leq C, ~\forall ~\Omega'\subset\subset B_1,
\end{equation*}
where $C$ depends only on $n,\Omega',\|v_m^+\|_{L^{\infty}(B_1)}$ and $\|v_m\|_{L^{\infty}(B_1)}$. From the Alexandrov-Bakel'man-Pucci maximum principle (see \cite[Theorem 6]{MR2854672}),
\begin{equation}\label{e2.17}
\|v_m^+\|_{L^{\infty}(B_1)},~\|v_m\|_{L^{\infty}(B_1)}\leq C,
\end{equation}
where $C$ depends only on $n$ and $K$. Hence, there exist $\bar{v}^+$ and $\bar{v}$ such that (up to a subsequence)
\begin{equation}\label{e2.13}
v_m^+\to \bar{v}^+,~v_m\to \bar{v}\quad\mbox{in}~~C^1(\bar{\Omega}'),~\forall ~\Omega'\subset\subset B_1.
\end{equation}
In addition, since $\|v_m^+\|_{W^{1,1}(B_1)},\|v_m\|_{W^{1,1}(B_1)}$ are uniformly bounded (see \cite[(16), P. 319]{MR0301343}),
\begin{equation}\label{e2.14}
v_m^+\to \bar{v}^+,~v_m\to \bar{v}\quad\mbox{weakly in}~~W^{1,1}(B_1).
\end{equation}

Note that $v_m^+-v_m\in W^{1,1}_0(B_1)$. By using it as the test function in \cref{e2.15}, we have
\begin{equation*}
\int_{B_1} \left(A(Dv_m^+)-A(Dv_m)\right)\left(Dv_m^+-Dv_m\right)
=\frac{1}{m}\int_{B_1}\left(v_m^+-v_m\right).
\end{equation*}
By \cref{e2.17} and
\begin{equation*}
(A(p)-A(q))\cdot (p-q)\geq \frac{|p-q|^2}{(1+|p|^2+|q|^2)^{3/2}},~\forall ~p,q\in \mathbb{R}^n,
\end{equation*}
we have
\begin{equation*}
\int_{B_1} \frac{|Dv_m^+-Dv_m|^2}{(1+|Dv_m^+|^2+|Dv_m|^2)^{3/2}}
\leq \frac{C}{m}.
\end{equation*}
Let $m\to \infty$ and by the Fatou's lemma,
\begin{equation*}
\int_{B_1} \frac{|D\bar{v}^+-D\bar{v}|^2}{(1+|D\bar{v}^+|^2+|D\bar{v}|^2)^{3/2}}
=0.
\end{equation*}
By \cref{e2.14}, $\bar{v}^+-\bar{v}\in W^{1,1}_0(B_1)$. Hence, $\bar{v}^+\equiv \bar{v}$ in $B_1$, which contradicts with \cref{e2.10}. Therefore, \cref{e2.3} holds.

By combining \cref{e2.9} with \cref{e2.3},
\begin{equation}\label{e2.11}
  \|u-v\|_{L^{\infty}(B_1)}\leq \delta_2.
\end{equation}
Since $v$ is smooth, there exists $P\in \mathcal{P}_1$ such that
\begin{equation*}
|v(x)-P(x)|\leq C|x|^{2},~\forall ~x\in B_{1/2},
\end{equation*}
and
\begin{equation*}
|DP|\leq C,
\end{equation*}
where $C$ depends only on $n$ and $\|u\|_{L^{\infty}(B_1)}$.

Take $0<r<1/2$ small enough such that $Cr<\delta_1/2$. Then
\begin{equation}\label{e2.12}
\|v-P\|_{L^{\infty}(B_{r})}\leq \frac{\delta_1}{2}r.
\end{equation}
In addition, take $\delta_2$ small enough such that $\delta_2\leq \delta_1r/2$. Therefore, by \cref{e2.11} and \cref{e2.12},
\begin{equation*}
\|u-P\|_{L^{\infty}(B_{r})}\leq \delta_1r.
\end{equation*}
~\qed~\\

Now, we give the~\\
\noindent\textbf{Proof of \Cref{th.MS}.} For (i), let $\tilde{u}=u-P$ and $\tilde{u}$ is a solution of
\begin{equation*}
G(D^2v,Dv)=F(D^2v,Dv+DP)=f\quad\mbox{in}~B_1,
\end{equation*}
where $F,G$ are defined as in \Cref{re2.2}. Hence, $G$ is smooth and is $1$-uniformly elliptic with ellipticity constants depending only on $|DP|$. Then by \Cref{th1.1}, the conclusion follows.

Next, we prove (ii). Let $\delta_1>0$ to be specified later. By \Cref{le2.4}, there exist $0<r<1/2$ and
$P\in \mathcal{P}_1$ such that
\begin{equation*}
\|u-P\|_{L^{\infty}(B_{r})}\leq \delta_1 r.
\end{equation*}
Let
\begin{equation*}
y=\frac{x}{r}, \quad  \tilde{u}(y)=\frac{u(x)-P(x)}{r}, \quad \tilde{f}(y)=rf(x).
\end{equation*}
Then $\tilde{u}$ is a solution of
\begin{equation*}
F(D^2\tilde{u},D\tilde{u}+DP)=\tilde f\quad\mbox{in}~B_1
\end{equation*}
and
\begin{equation*}
\|\tilde{u}\|_{L^{\infty}(B_1)}\leq \delta_1, \quad \|\tilde{f}\|_{L^{\infty}(B_1)}\leq \|f\|_{L^{\infty}(B_1)}\leq n-1.
\end{equation*}
Then $\tilde u\in C^{1,\alpha}(0)$ by \Cref{th1.1} provided $\delta_1$ is small enough, which is guaranteed by taking $\delta$ small enough.

Next, we prove (iii). For $r>0$, let
\begin{equation}\label{e2.22}
y=\frac{x}{r}, \quad  \tilde{u}(y)=\frac{u(x)-u(0)}{r}, \quad \tilde{f}(y)=r f(x).
\end{equation}
Then $\tilde{u}$ is a solution of
\begin{equation*}
F(D^2\tilde{u},D\tilde{u})=\tilde{f}\quad\mbox{in}~B_1.
\end{equation*}
By the assumption, $\tilde{u}$ is bounded. We choose $r$ small enough (depending on $\|f\|_{L^{\infty}(B_1)}$) such that $\|\tilde f\|_{L^{\infty}}$ is small. Then the conclusion follows from (ii).

Finally, we prove (iv). Let
\begin{equation*}
y=\frac{x}{\delta}, \quad  \tilde{u}(y)=\frac{u(x)-P(x)}{\delta}, \quad \tilde{f}(y)=\delta f(x).
\end{equation*}
Then $\tilde{u}$ is a solution of
\begin{equation*}
F(D^2\tilde{u},D\tilde{u}+DP)=\tilde{f}\quad\mbox{in}~B_1.
\end{equation*}
Take $\delta$ small enough such that $\|\tilde f\|_{L^{\infty}(B_1)}$ is small. In addition, since $\tilde{u}=0$ on $\partial B_1$, by the Alexandrov-Bakel'man-Pucci maximum principle (see \cite[Theorem 6]{MR2854672}), $\|\tilde{u}\|_{L^{\infty}(B_1)}$ is small. Then the conclusion follows from \Cref{th1.1}.~\qed~\\

The next is the~\\
\noindent\textbf{Proof of \Cref{th.MS-2}.} By \Cref{th.MS}, $u\in C^{1,\alpha}(0)$. For $r>0$, let
\begin{equation*}
y=\frac{x}{r}, \quad  \tilde{u}(y)=\frac{u(x)-P_u(x)}{r}, \quad \tilde{f}(y)=rf(x).
\end{equation*}
Then $\tilde{u}$ is a solution of
\begin{equation*}
F(D^2\tilde{u},D\tilde{u}+DP_u)=\tilde{f}\quad\mbox{in}~B_1.
\end{equation*}
Since $u\in C^{1,\alpha}(0)$,
\begin{equation*}
\begin{aligned}
&\|\tilde{u}\|_{L^{\infty}(B_1)}= r^{-1}\|u-P_u\|_{L^{\infty}(B_r)}\leq Cr^{\alpha},\\
&\|\tilde f\|_{C^{k-2,\alpha}(0)}\leq \sum_{i=0}^{k-2}r^{i+1}|D^if(0)|+r^{k-1+\alpha}[f]_{C^{k-2,\alpha}(0)}.
\end{aligned}
\end{equation*}
Hence, we can choose $r$ small enough such that $\|\tilde{u}\|_{L^{\infty}}$ and $\|\tilde f\|_{C^{k-2,\alpha}(0)}$ are small. Then the conclusion follows from \Cref{th1.2} and \Cref{th1.3}.~\qed~\\

At the end of this section, we give two remarks.
\begin{remark}\label{re2.3}
Since the prescribed mean curvature equation has the divergence structure, one may consider a weak solution $u\in W^{1,1}(B_1)$ of \cref{e.MS} rather than a viscosity solution. We have the conclusion for weak solutions as well. For example, let us assume
\begin{equation*}
\|u\|_{L^{\infty}(B_1)}\leq \delta, \quad \|f\|_{L^{\infty}(B_1)}\leq \delta,
\end{equation*}
where $0<\delta<1$ depends only on $n$ and $\alpha$.

Then we can approach the regularity by an approximation argument similar to \cite[Proof of Theorem 1]{MR0301343}. Take sequences of smooth functions $u_m,f_m$ such that
\begin{equation*}
u_m\to u \quad\mbox{in}~~ W^{1,1}(B_1), \quad f_m\to f \quad\mbox{in}~~L^2(B_1)
\end{equation*}
and
\begin{equation*}
  \|u_m\|_{L^{\infty}(B_1)}\leq \delta, \quad \|f_m\|_{L^{\infty}(B_1)}\leq \delta,~\forall ~m\geq 1.
\end{equation*}

Let $v_m\in C^{\infty}(\bar{B}_1)$ be solutions of (see \cite[Theorem 16.10]{MR1814364})
\begin{equation*}
  \left\{
  \begin{aligned}
\mathrm{div} A(Dv_m)=&f_m\quad\mbox{in}~~B_1;\\
v_m=&u_m\quad\mbox{on}~~\partial B_1.
  \end{aligned}
  \right.
\end{equation*}
By the Alexandrov-Bakel'man-Pucci maximum principle,
\begin{equation*}
\|v_m\|_{L^{\infty}(B_1)}\leq C\delta,
\end{equation*}
where $C$ depends only on $n$. Then from (i) of \Cref{th.MS},
\begin{equation*}
\|v_m\|_{C^{1,\alpha}(\bar{B}_{1/2})}\leq C,
\end{equation*}
where $C$ depends only on $n$ and $\alpha$. Hence, there exists $\bar{v}\in C^{1,\alpha}(\bar{B}_{1/2})$ such that (up to a subsequence)
\begin{equation*}
v_m\to \bar{v}\quad\mbox{in}~~C^1(\bar{B}_{1/2}).
\end{equation*}

As in the proof of \Cref{le2.4},
\begin{equation*}
\int_{B_1} \left(A(Dv_m)-A(Du)\right)\left(Dv_m-Du_m\right)
=\int_{B_1}(f_m-f)\left(v_m-u_m\right).
\end{equation*}
Hence,
\begin{equation*}
  \begin{aligned}
\int_{B_1}& \frac{|Dv_m-Du|^2}{(1+|Dv_m|^2+|Du|^2)^{3/2}}\leq
\int_{B_1} \left(A(Dv_m)-A(Du)\right)\left(Dv_m-Du\right)\\
=&\int_{B_1}(f_m-f)\left(v_m-u_m\right)
+\int_{B_1} \left(A(Dv_m)-A(Du)\right)\left(Du_m-Du\right)\\
\leq& C\|f_m-f\|_{L^2(B_1)}+2\|Du_m-Du\|_{L^1(B_1)}.
  \end{aligned}
\end{equation*}
Let $m\to \infty$ and we have
\begin{equation*}
\int_{B_{1/2}} \frac{|D\bar{v}-Du|^2}{(1+|D\bar{v}|^2+|Du|^2)^{3/2}}=0.
\end{equation*}
Thus, for some constant $c_0$,
\begin{equation*}
u\equiv \bar{v}+c_0\quad\mbox{in}~~B_{1/2}.
\end{equation*}
Therefore, $u\in C^{1,\alpha}(\bar{B}_{1/2})$.
\end{remark}

\begin{remark}\label{re2.5}
The theory for locally uniformly elliptic equations is not applicable to the $p$-Laplace equations. Indeed, the $p$-Laplace equations are more like uniformly elliptic equations (in particular when $Du\neq 0$). Hence, one may use the uniformly elliptic equations techniques to study the $p$-Laplace equations (e.g. \cite{MR1351007,MR1264526}). In addition, the $p$-Laplace operator is not smooth.
\end{remark}
~\\

\section{Monge-Amp\`{e}re equation}\label{S.MA}
In this section, we consider the Monge-Amp\`{e}re equation:
\begin{equation*}
\det D^2u=f\quad\mbox{in}~~B_1.
\end{equation*}
We have the following observation. For any convex polynomial $P\in\mathcal{P}_2$ with $\det(D^2P)=1$, define
\begin{equation*}
F(M)=\det (M+D^2P)-1,~\forall ~M\in \mathcal{S}^n.
\end{equation*}
Then $F\in C^{\infty}, F(0)=0$. Moreover, $F$ is $\rho$-uniformly elliptic with $\tilde{\lambda},\tilde{\Lambda}$ and these three constants depend only on $n$ and $|D^2P|$. Therefore, the regularity theory for locally uniformly elliptic equations is applicable.

We have the following interior pointwise $C^{k,\alpha}$ regularity.
\begin{theorem}\label{th8.1}
Let $k\geq 2,0<\alpha<1$ and $u$ be a strictly convex viscosity solution of
\begin{equation*}
 \det D^2u=f\quad\mbox{in}~B_1.
\end{equation*}
Suppose that $0<\lambda\leq f\leq \Lambda$ and $f\in C^{k-2,\alpha}(0)$. Then $u\in C^{k,\alpha}(0)$.
\end{theorem}

If the dimension $n=2$, the solution is always strictly convex (see \cite{MR7626}, \cite{MR116293} \cite[Theorem 2.19]{MR3617963}, \cite[Remark 3.2]{MR2483373}). Hence, we have
\begin{corollary}\label{th.MA-2}
Let $k\geq 2,0<\alpha<1$ and $u$ be a convex viscosity solution of
\begin{equation*}
 \det D^2u=f\quad\mbox{in}~B_1\subset \mathbb{R}^2.
\end{equation*}
Suppose that $0<\lambda\leq f\leq \Lambda$ and $f\in C^{k-2,\alpha}(0)$. Then $u\in C^{k,\alpha}(0)$.
\end{corollary}

By applying above result to the prescribed Gaussian curvature equation, we have
\begin{corollary}\label{th.MA-3}
Let $k\geq 2,0<\alpha<1$ and $u$ be a strictly convex viscosity solution of
\begin{equation*}
 \det D^2u=K(x)\left(1+|Du|^2\right)^{\frac{n+2}{2}}\quad\mbox{in}~B_1.
\end{equation*}
Suppose that $0<\lambda\leq K\leq \Lambda$ and $K\in C^{k-2,\alpha}(0)$. Then $u\in C^{k,\alpha}(0)$.
\end{corollary}
\proof Since $u$ is convex,
\begin{equation*}
  \|Du\|_{L^{\infty}(B_{1/2})}\leq C\|u\|_{L^{\infty}(B_1)}.
\end{equation*}
Hence, the right-hand term is bounded between two positive constants. Then from the strict convexity of $u$,
we have $u\in C^{1,\beta}(\bar{B}_{1/2})$ for some $0<\beta<1$ (see \cite[Theorem 2]{MR1127042}, \cite[Corollary 4.21]{MR3617963} and \cite[Lemma 3.5]{MR2483373}). Therefore, the conclusion follows from \Cref{th8.1}.~\qed~\\

\begin{remark}\label{re2.9}
The $C^{2,\alpha}$ regularity for the Monge-Amp\`{e}re equation is well-known. Sabitov \cite{MR425854} and Schulz \cite{MR643043} proved the $C^{2,\alpha}$ regularity for $n=2$. If $f\in C^{0,1}$, the $C^{2,\alpha}$ regularity for general dimensions was derived by Urbas \cite{MR926846}. The $C^{2,\alpha}$ regularity for general dimensions was due to Caffarelli \cite{MR1038360} (see also \cite{MR2338423}). Trudinger and Wang obtained boundary $C^{2,\alpha}$ regularity \cite{MR2415390} and the pointwise version was proved by Savin \cite{MR2983006}. Of course, Savin's proof is also applicable to derive interior pointwise $C^{2,\alpha}$ regularity.

As pointed out in \Cref{re1.8}, we can obtain pointwise $C^{k}$ and $C^{k,\mathrm{Ln}L}$ ($k\geq 2$) regularity as well. In this respect, the $C^2$ regularity was proved by Wang \cite{WXJ-1992} and  $C^{2,\mathrm{Ln}L}$ regularity was proved by Jian and Wang \cite{MR2338423}.
\end{remark}

\begin{remark}\label{re2.1}
If the dimension $n\geq 3$, by the well-known Pogorelov's counterexample (see \cite[P. 81-83]{MR478079} and a good explanation on this counterexample in \cite[Chapter 3.2]{MR3617963}), the condition ``$u$ is strictly convex'' can not be dropped. In addition, ``$f\geq \lambda>0$'' cannot be replaced by ``$f\geq 0$''. Indeed, the best regularity for the latter is $C^{1,1}$ in general (see \cite[Theorem 2]{MR766792} and \cite{MR1687172} for $C^{1,1}$ regularity; see  \cite[Example 2]{MR864651}, \cite[(1.3), P. 88]{MR1687172} and \cite[Example 3 and Remark 1]{MR1223269} for counterexamples).

In conclusion, by transforming the Monge-Amp\`{e}re equation to a locally uniformly elliptic equation, we can obtain the best expected regularity, which may imply that this viewpoint is essential. That the Monge-Amp\`{e}re equation can be transformed to a locally uniformly elliptic equation has been noted before (e.g. \cite[P. 673, L. 1]{MR3158522}). We learned this from a note written by Prof. Chuanqiang Chen.
\end{remark}

\begin{remark}\label{re3.1}
There are several criterions to ensure that $u$ is strictly convex in a general domain $\Omega$, such as
\begin{itemize}
  \item $u$ is an affine function on $\partial \Omega$ (see \cite[Corollary 2]{MR1038359}, \cite[Theorem 4.10]{MR3617963} and \cite[Lemma 3.4]{MR2483373}).
  \item $u\in C^{1,\alpha}(\partial \Omega)$ ($\alpha>1-2/n$) (see \cite[Corollary 4]{MR1038359}, \cite[Corollary 4.11]{MR3617963} and \cite[Remark 3.1]{MR2483373}).
  \item $u\in W^{2,p}(\Omega)$ ($p>n(n-1)/2$) (see \cite{MR769173} and \cite{MR926846}).
\end{itemize}
\end{remark}

\begin{remark}\label{re2.7}
Another important notion of weak solution is the Alexandrov's generalized solution (see \cite{MR96903} and \cite[Definition 2.5]{MR3617963}). The right-hand term $f$ is not necessarily continuous if we use this notion. These two definitions are equivalent if $f$ is continuous (see \cite[Proposition 1.3.4, Proposition 1.7.1]{MR1829162}).

\Cref{th8.1} holds for Alexandrov's generalized solutions as well. In the following, we make an explanation. Since $u$ is strictly convex, without loss of generality, we can assume that $u$ is an Alexandrov's generalized solution of
\begin{equation*}
  \left\{\begin{aligned}
    \det D^2u&=f &&~\mbox{in}~\Omega;\\
    u&=0&&~\mbox{on}~\partial\Omega,
  \end{aligned}\right.
\end{equation*}
where $\Omega$ is a convex domain and $B_{1/n}\subset \Omega \subset B_1$. Choose a smooth function $\varphi$ with compact support in $B_1$ to mollify $f$ by convolution: (see \cite[Lemma 3.5.6]{MR1014685})
\begin{equation*}
f_{\varepsilon}:=f*\varphi_{\varepsilon}, \quad \varphi_{\varepsilon}(x):=\varepsilon^{-n}\varphi(x/\varepsilon), \quad \varepsilon>0
\end{equation*}
such that
\begin{equation*}
P=P*\varphi_{\varepsilon},~\forall ~P\in \mathcal{P}_{k-2},~\varepsilon>0.
\end{equation*}
Then $f_{\varepsilon}\in C^{k-2,\alpha}(0)$ with $P_{f_{\varepsilon}}\equiv P_{f}$.

Let $u_{\varepsilon}$ be Alexandrov's solutions (be viscosity solutions as well) of (see \cite[Theorem 2.13]{MR3617963})
\begin{equation*}
  \left\{\begin{aligned}
    \det D^2u_{\varepsilon}&=f_{\varepsilon} &&~\mbox{in}~\Omega;\\
    u_{\varepsilon}&=0&&~\mbox{on}~\partial\Omega.
  \end{aligned}\right.
\end{equation*}
Since $f_{\varepsilon}\to f$ in $L^1(\Omega)$, we have (see \cite[Proposition 2.16]{MR3617963})
\begin{equation*}
  u_{\varepsilon}\to u \quad\mbox{in}~~L^{\infty}(\Omega).
\end{equation*}
By \Cref{th8.1}, $u_{\varepsilon}\in C^{k,\alpha}(0)$. That is, there exist $P_{\varepsilon}\in \mathcal{P}_k$ such that
\begin{equation}\label{e2.16}
  |u_{\varepsilon}(x)-P_{\varepsilon}(x)|\leq C|x|^{k+\alpha},~\forall ~x\in \Omega
\end{equation}
and
\begin{equation*}
  \|P_{\varepsilon}\|\leq C,
\end{equation*}
where $C$ is independent of $\varepsilon$. Hence, up to a subsequence, there exists $P\in \mathcal{P}_k$ such that
$P_{\varepsilon}\to P$. Let $\varepsilon\to 0$ in \cref{e2.16} and then
\begin{equation*}
  |u(x)-P(x)|\leq C|x|^{k+\alpha},~\forall ~x\in \Omega.
\end{equation*}
That is, $u\in C^{k,\alpha}(0)$.
\end{remark}
~\\

We first prove a lemma.
\begin{lemma}\label{le8.3}
Let $B_{1/n}\subset \Omega \subset B_1$ and $u$ be a viscosity solution of
\begin{equation*}
  \left\{\begin{aligned}
    \det D^2u&=f &&~\mbox{in}~\Omega;\\
    u&=0&&~\mbox{on}~\partial\Omega.
  \end{aligned}\right.
\end{equation*}
Assume that $0$ is the minimum point of $u$.

Then for any $\delta_1>0$, there exist $\delta,r>0$ (depending only on $n$ and $\delta_1$) such that for some $P\in \mathcal{P}_2$ with $D^2P\geq 0$,
\begin{equation*}
\|u-P\|_{L^{\infty}(B_{r})}\leq \delta_1 r^2
\end{equation*}
and
\begin{equation}\label{e8.8}
  \det D^2P=1,\quad \|P\|\leq C,
\end{equation}
provided
\begin{equation*}
  \|f-1\|_{L^{\infty}(\Omega)}\leq \delta,
\end{equation*}
where $C$ depends only on $n$.
\end{lemma}
\proof Suppose not. Then there exist $\delta_1>0$, sequences of $u_m,\Omega_m$ and $f_m$ such that
$B_{1/n}\subset \Omega_m \subset B_1$,
\begin{equation*}
  \left\{\begin{aligned}
    \det D^2u_m&=f_m &&~\mbox{in}~\Omega_m;\\
    u_m&=0&&~\mbox{on}~\partial\Omega_m,
  \end{aligned}\right.
\end{equation*}
and
\begin{equation*}
  \|f_m-1\|_{L^{\infty}(\Omega_m)}\leq 1/m.
\end{equation*}
Moreover, for any convex polynomial $P\in\mathcal{P}_2$ satisfying \cref{e8.8}, we have
\begin{equation}\label{e8.5}
\|u_m-P\|_{L^{\infty}(B_{r})}\geq  \delta_1 r^2,
\end{equation}
where $r$ is to be specified later.

By the stability of solutions (see \cite[Corollary 2.12]{MR3617963}), up to a subsequence, there exist
$\bar{u},\tilde{\Omega}$ such that
\begin{equation*}
  u_m\to \bar{u}, \quad \Omega_m\to \tilde{\Omega}, \quad B_{1/n}\subset \tilde{\Omega}\subset B_1.
\end{equation*}
Moreover, $\bar{u}$ is the solution of
\begin{equation*}
  \left\{\begin{aligned}
    \det D^2\bar u&=1 &&~\mbox{in}~\tilde\Omega;\\
    \bar u&=0&&~\mbox{on}~\partial\tilde\Omega.
  \end{aligned}\right.
\end{equation*}

Since $0$ is the minimum point of $u_m$, we have (see \cite[Proposition 4.4]{MR3617963})
\begin{equation*}
d(0,\partial \Omega_m)\geq c_0>0,
\end{equation*}
where $c_0>0$ depends only on $n$. Hence,
\begin{equation*}
d(0, \partial \tilde\Omega)\geq c_0.
\end{equation*}
By the interior regularity (see \cite[Theorem
3.10]{MR3617963}), there exists a convex polynomial $P\in\mathcal{P}_2$ such that \cref{e8.8} holds and
\begin{equation*}
|\bar u(x)-P(x)|\leq C|x|^3,~\forall ~x\in B_{c_0/2}.
\end{equation*}
Take $r$ small such that
\begin{equation*}
  Cr=\frac{\delta_1}{2}.
\end{equation*}
Then
\begin{equation}\label{e8.7}
   \|\bar u-P\|_{L^{\infty}(B_{r})}\leq \frac{1}{2} \delta_1 r^2.
\end{equation}
By taking the limit in \cref{e8.5}, we have
\begin{equation*}
   \|\bar u-P\|_{L^{\infty}(B_{r_0})}\geq \delta_1 r^2,
\end{equation*}
which contradicts with \cref{e8.7}.~\qed~\\

Now, we give the~\\
\noindent\textbf{Proof of \Cref{th8.1}.} By subtracting an affine function, we may assume that $u\geq 0$ and $u(0)=0$. Since $u$ is strictly convex, for $h>0$ small (to be specified later), $S_h(0)\subset \subset B_1$, where $S_h(0)$ the section of $u$ at
$0$, i.e.
\begin{equation*}
S_h(0):=\left\{x\in B_1:~~u(x)<h\right\}.
\end{equation*}
With the aid of John's lemma (see \cite{MR30135} and \cite[A.3.2]{MR3617963}), we normalize the section $S_h(0)$ as follows:
\begin{equation}\label{e2.21}
y=T_hx, \quad \tilde{\Omega}=T_h(S_h(0)),
\quad B_{1/n}(\tilde y)\subset \tilde{\Omega}\subset B_1(\tilde y), \quad \tilde{u}(y)=\frac{u(x)}{h},
\end{equation}
where $T_h\in \mathcal{S}^n$ and $\tilde{y}\in B_1$. Then $\tilde{u}$ is a solution of
\begin{equation*}
  \left\{\begin{aligned}
    \det D^2\tilde{u}&=\tilde f &&~\mbox{in}~\tilde{\Omega};\\
    \tilde{u}&=1&&~\mbox{on}~\partial \tilde{\Omega},
  \end{aligned}\right.
\end{equation*}
where
\begin{equation*}
  \tilde{f}(y)=\frac{f(x)}{(\det T_h)^2h^n}.
\end{equation*}

By the uniform estimate (see \cite[Lemma 3.2]{MR2483373}),
\begin{equation}\label{e2.23}
C_1\leq (\det T_h)^2h^n\leq C_2,
\end{equation}
where $C_1,C_2$ depends only on $n,\lambda,\Lambda$. Without loss of generality, we assume $\tilde f(0)=1$. Since
$f\in C^{k-2,\alpha}(0)$,
\begin{equation*}
|\tilde{f}(y)-1|\leq C|f(x)-f(0)|\leq C|x|^{\alpha}\leq C|T_h^{-1}|\cdot |y|^{\alpha}.
\end{equation*}
Since $u$ is strictly convex,
\begin{equation}\label{e2.19}
|T_h^{-1}|\to 0\quad\mbox{as}~~h\to 0.
\end{equation}

Let $0<\delta_1<1$ to be determined later. By \cref{e2.19}, we can take $h$ small enough such that
\begin{equation}\label{e2.20}
C\|T_h^{-1}\|\leq \delta_1, \quad \|\tilde f-1\|_{L^{\infty}(\tilde\Omega)}\leq \delta,
\end{equation}
where $\delta\leq \delta_1$ is chosen small enough such that \Cref{le8.3} holds. Then by \Cref{le8.3}, there exist $r>0$ and a convex polynomial $P\in \mathcal{P}_2$ such that
\begin{equation}\label{e2.8}
\|\tilde{u}-P\|_{L^{\infty}(B_{r})}\leq \delta_1 r^2,
\end{equation}
and
\begin{equation*}
  \det D^2P=1,\quad \|P\|\leq C.
\end{equation*}

Let
\begin{equation*}
  z=\frac{y}{r}, \quad \hat{u}(z)=\frac{\tilde{u}(y)-P(y)}{r^2}.
\end{equation*}
Then $\hat{u}$ is a solution of
\begin{equation}\label{e2.7}
F(D^2\hat{u})=\hat{f}~~\mbox{ in}~B_1,
\end{equation}
where
\begin{equation*}
F(M)=\det (M+D^2P)-1,~\forall ~M\in \mathcal{S}^n, \quad \hat{f}(z)=\tilde{f}(y)-1.
\end{equation*}
Note that $F\in C^{\infty}, F(0)=0$ and $F$ is $\rho$-uniformly elliptic with ellipticity constants $\tilde{\lambda},\tilde{\Lambda}$ and they depend only on $n,\lambda,\Lambda$.

By \cref{e2.20}, \cref{e2.8} and the definition of $\hat{f}$,
\begin{equation*}
  \|\hat{u}\|_{L^{\infty}(B_1)}\leq \delta_1,\quad
  \|\hat{f}\|_{C^{\alpha}(0)}\leq \delta\leq \delta_1.
\end{equation*}
From \Cref{th1.2} and \Cref{co1.1} (choosing $\delta_1$ small enough), $\hat{u}\in C^{k,\alpha}(0)$ and hence $u\in C^{k,\alpha}(0)$. ~\qed~\\

\begin{remark}\label{re3.2}
Two cornerstone results for \Cref{th8.1} are a priori estimates/existence of classical solutions and Pogorelov's type estimate. They are implicitly used in \Cref{le8.3}. We refer to \cite[Theorem
3.10]{MR3617963} for details.
\end{remark}

\begin{remark}\label{re2.8}
For the Monge-Amp\`{e}re equation, we obtain almost the same $C^{k,\alpha}$ regularity as the Poisson's equation. However, for the prescribed mean curvature equation (and other types of equations below), we must assume that $u$ (or $f$) is small or $u$ is uniformly elliptic at $0$ (i.e. $u\in C^{0,1}$). The reason is that for the Monge-Amp\`{e}re equation, we can make an anisotropic scaling such that the equation is unchanged. This is the unique feature of the Monge-Amp\`{e}re equation.

For the prescribed mean curvature equation, to make $f$ small and the equation unchanged, we have to make the following scaling (see \cref{e2.22}):
\begin{equation}\label{e2.24}
  y=\frac{x}{r}, \quad  \tilde{u}(y)=\frac{u(x)-u(0)}{r}.
\end{equation}
Thus, to guarantee that $\tilde{u}$ is bounded, we must assume that $u\in C^{0,1}(0)$.

On the contrast, for the Monge-Amp\`{e}re equation, to make $f$ small and the equation unchanged, we can make an anisotropic scaling (see \cref{e2.21}):
\begin{equation}\label{e2.25}
y=T_hx, \quad \tilde{\Omega}=T_h(S_h(0)), \quad \tilde{u}(y)=\frac{u(x)}{h},
\end{equation}
where $T_h$ is anisotropic, i.e., its eigenvalues may be not comparable. Note that $h\simeq r^2$ in the sense of anisotropy (see \cref{e2.23}). Hence, the scaling transformation \cref{e2.25} is essentially the same to \cref{e2.24} except that \cref{e2.25} is anisotropic. Maybe this is the reason that there are few pointwise Schauder regularity for non-uniformly equations except the Monge-Amp\`{e}re equation.
\end{remark}
~\\

\section{$k$-Hessian equations}\label{S.H}
Next, we consider the $k$-Hessian equations. For $u\in C^2$, denote its eigenvalues by
\begin{equation*}
\lambda(D^2u)=(\lambda_1,\cdots,\lambda_n)\in{\mathbb{R}}^n.
\end{equation*}
Define for $1\leq k\leq n$
\begin{equation*}
\sigma_k(D^2u)=\sigma_k(\lambda(D^2u))=\sum_{1\leq i_1<\cdots<i_k\leq n}\lambda_{i_1}\cdots\lambda_{i_k}.
\end{equation*}
We say $u$ is $k$-admissible if $ \lambda(D^2u)\in\Gamma_k,$ where $\Gamma_k$ is the G{\aa}rding convex cone in ${\mathbb{R}}^n$ defined by
  \begin{equation*}
   \Gamma_k\equiv\big\{\lambda\in{\mathbb{R}}^n\big|\ \sigma_i(\lambda)>0, \quad \forall~ 1\leq i\leq k\big\}.
  \end{equation*}
For more basic knowledge of the $k$-Hessian equations, we refer to \cite{MR2500526}. We consider the $k$-Hessian equation:
\begin{equation*}
    \sigma_k(D^2u)=f ~~\mbox{ in}~B_1.
\end{equation*}
Similar to the Monge-Amp\`{e}re equation, for any $k$-admissible polynomial $P\in\mathcal{P}_2$ with $\sigma_k(D^2P)=1$, define
\begin{equation*}
F(M)=\sigma_k (M+D^2P)-1,~\forall ~M\in \mathcal{S}^n.
\end{equation*}
Then $F\in C^{\infty}, F(0)=0$ and $F$ is $\rho$-uniformly elliptic with $\tilde{\lambda},\tilde{\Lambda}$ which depend only on $n,k$ and $|D^2P|$.

Now, we state the main results in this section. Since the $k$-Hessian equation reduces to the Monge-Amp\`{e}re equation if the dimension $n=2$, we only consider $n\geq 3$ in this section.
\begin{theorem}\label{th.H}
Let $2\leq k\leq n-1$, $0<\alpha<1$ and $u\in C(\bar{B}_1)$ be a viscosity solution of
\begin{equation}\label{e.H}
    \sigma_k(D^2u)=f ~~\mbox{ in}~B_1,
\end{equation}
where $0<\lambda\leq f\leq \Lambda$. Then $u\in C^{l,\alpha}(0)$ ($l\geq 2$) provided one of the
following conditions holds :\\
(i) there exists a $k$-admissible $P\in \mathcal{P}_2$ such that
\begin{equation*}
\|u-P\|_{L^{\infty}(B_1)}\leq \delta, \quad \sigma_k(D^2P)=f(0), \quad  f\in C^{l-2,\alpha}(0), \quad
[f]_{C^{\alpha}(0)}\leq \delta,
\end{equation*}
where $0<\delta<1$ depends only on $n,k,\lambda,\Lambda,\alpha$ and $|D^2P|$.\\
(ii) there exists a $k$-admissible $P\in \mathcal{P}_2$ such that
\begin{equation*}
\|u-P\|_{L^{\infty}(B_1)}\leq \delta, \quad \sigma_k(D^2P)=f(0), \quad  f\in C^{l-2,\alpha}(0),
\end{equation*}
where $0<\delta<1$ depends only on $n,k,\lambda,\Lambda,\alpha,|D^2P|$ and $[f]_{C^{\alpha}(0)}$.\\
(iii) $u\in C^{2}(0)$ and $f\in C^{l-2,\alpha}(0)$.\\
(iv) there exists $p>k(n-1)/2$ such that
\begin{equation*}
u\in W^{2,p}(B_1), \quad f\in C^{l-2,\alpha}(0), \quad [f]_{C^{\alpha}(0)}\leq \delta,
\end{equation*}
where $\delta$ depends only on $n,k,\lambda,\Lambda,\alpha$ and $\|u\|_{W^{2,p}(B_1)}$.~\\
(v)  there exists a $k$-admissible $P\in \mathcal{P}_2$ such that
\begin{equation*}
u=P~~\mbox{ on}~\partial B_{\delta}, \quad \sigma_k(D^2P)=f(0), \quad  f\in C^{l-2,\alpha}(0),
\end{equation*}
where $\delta$ depends only on $n,k,\lambda,\Lambda,\alpha$ and $|D^2P|$.~\\

\end{theorem}

Chaudhuri and Trudinger \cite{MR2133413} proved that if $k>n/2$, $u\in C^2(x_0)$ for \emph{a.e.} $x_0$. Hence, we have the following corollary of (iii).
\begin{corollary}\label{co.4.1}
Let $k>n/2$, $ l\geq 2, 0<\alpha<1$ and $u\in C(\bar{B}_1)$ be a viscosity solution of
\begin{equation}\label{e.H-2}
    \sigma_k(D^2u)=f ~~\mbox{ in}~B_1,
\end{equation}
where $0<\lambda\leq f\leq \Lambda$ and $f\in C^{l-2,\alpha}(\bar{B}_1)$. Then $u\in C^{l,\alpha}(x_0)$ for \emph{a.e.} $x_0\in B_1$.
\end{corollary}

\begin{remark}\label{re4.2}
We point out that $\sigma_k(D^2P)=f(0)$ is not needed. In fact, with the aid of the definition of viscosity solution, we can modify $P$ such that this equality holds if $\delta$ is small enough.

We do not know whether the $C^{2,\alpha}$ regularity holds if we only assume $u\in C^{1,1}(0)$ instead of $u\in C^2(0)$ in (iii). On the other hand, $u\in C^2(0)$ can be relaxed to (by a similar proof to that of \Cref{th.H2})
\begin{equation*}
K_1|x|^2\leq u-P\leq K_1|x|^2\quad\mbox{in}~~B_1,
\end{equation*}
where $K_1,K_2>0$ and $P\in\mathcal{P}_1$. A similar condition has been used by Wang and Bao \cite[Theorem 1.1]{MR4438235} to prove the rigidity in $\mathbb{R}^n$ for the $k$-Hessian equation.

Urbas \cite{MR1840289} (see also \cite{MR769173,MR1777141} for similar estimates) proved the interior $C^{1,1}$ estimate if $p>k(n-1)/2$ and $f\in C^{1,1}$. This is the cornerstone of the regularity under the assumption (iv).

The condition (v) concerns the regularity in a small domain. In this direction, Urbas \cite[Theorem 3]{MR1089043} proved global $C^1$ estimate. Tian, Wang and Wang \cite{MR3536033} studied the local solvability of the $k$-Hessian equation in $B_{r_0}$ for $r_0$ small. Maybe \Cref{th.H} can be applied to study the local solvability.
\end{remark}
~\\

Since the $k$-Hessian equations do not possess pure interior $C^{1,1}$ estimate, we do not have a similar regularity by assuming only $f-f(0)$ small as in (ii) \Cref{th.MS}. Instead, we have the following regularity based on a priori estimates and the Pogorelov's type estimate.
\begin{theorem}\label{th.H2}
Let $2\leq k\leq n-1$, $0<\alpha<1$ and $u\in C(\bar{\Omega})$ be a viscosity solution of
\begin{equation*}
  \left\{\begin{aligned}
    \sigma_k(D^2u)&=f &&~\mbox{in}~\Omega;\\
    u&=g&&~\mbox{on}~\partial \Omega,
  \end{aligned}\right.
\end{equation*}
where $0<\lambda\leq f\leq \Lambda$. Then $u\in C^{l,\alpha}(0)$ ($l\geq 2$) provided one of the
following conditions holds :\\
(i) $\partial \Omega\in C^{3,1}$, $\Omega$ is $(k-1)$-convex and
\begin{equation*}
\|g-g_0\|_{L^{\infty}(\partial \Omega)}\leq
\delta,~g_0\in C^{3,1}(\partial \Omega),~ f\in C^{l-2,\alpha}(0),~[f]_{C^{\alpha}(0)}\leq \delta,
\end{equation*}
where $0<\delta<1$ depends only on $n,k,\lambda,\Lambda,\alpha,\Omega$ and
$\|g_0\|_{C^{3,1}(\partial \Omega)}$.\\
(ii) $\Omega$ is a $(k-1)$-convex domain and
\begin{equation*}
\|g-P\|_{L^{\infty}(\partial \Omega)}\leq
\delta,~P\in\mathcal{P}_1,~f\in C^{l-2,\alpha}(0),~[f]_{C^{\alpha}(0)}\leq \delta,
\end{equation*}
where $0<\delta<1$ depends only on $n,k,\lambda,\Lambda,\alpha$ and
$\Omega$.\\
(iii)
\begin{equation*}
\Omega=B_{\delta},~~\|g-P\|_{L^{\infty}(\partial B_{\delta})}\leq \delta_ 1\delta^2,~~P\in\mathcal{P}_1,~~
 f\in C^{l-2,\alpha}(0),
\end{equation*}
where $\delta_1$ depends only on $n,k,\lambda,\Lambda,\alpha$ and $\delta$ depends also on $[f]_{C^{\alpha}(0)}$.~\\

\end{theorem}

\begin{remark}\label{re4.4}
The regularity under assumption (i) is based on a priori estimates for sufficient smooth data, which was first proved by Caffarelli, Nirenberg and Spruck \cite{MR806416} (see also \cite{MR1368245} and \cite[Theorem 3.4]{MR2500526}). The regularity under assumption (ii) is based on the Pogorelov's type estimate, established by Chou and Wang \cite[Theorem 1.5]{MR1835381} (see also \cite[Theorem 4.3]{MR2500526}). In (ii), $\partial \Omega$ is not necessary to be smooth. The $(k-1)$-convexity can be defined in some weak sense (e.g. by an approximation, see \cite[P. 226, L. 4]{MR1634570}; in the viscosity sense, see \cite[P. 580]{MR1726702}).
\end{remark}

\begin{remark}\label{re4.3}
There are other types of weak solutions, e.g., the weak solution defined by an approximation (see \cite{MR1466315}, \cite{MR1777141}, \cite{MR1840289}) and the weak solution based on the Borel measure (see \cite{MR1634570,MR1726702,MR1923626}). As before, \Cref{th.H} and \Cref{th.H2} hold with these weak solutions by an approximation argument.
\end{remark}

\begin{remark}\label{re4.5}
We do not know any counterexample for the $k$-Hessian equation to demonstrate the necessity of the ``smallness'' as in the prescribed mean curvature equation (see \Cref{re1.0}), the $k$-Hessian quotient equations (see \Cref{re5.2}) and the Lagrangian mean curvature equation (see \Cref{re6.2}). Hence, it is interesting to construct a counterexample, e.g., a strictly convex viscosity solution $u$ of \cref{e.H} with $f\in C^{\alpha}$ but $u\notin C^{2,\alpha}$. Or, can we prove \Cref{th.H} and \Cref{th.H2} without the smallness assumptions? Note that the counterexample constructed by Pogorelov and its extensions (see \cite{Caf_1,MR4468909,lu2024interior}, \cite[P. 6]{MR1777141}, \cite{MR1089043} etc.) are not strictly convex.
\end{remark}
~\\

We first prove a lemma.
\begin{lemma}\label{le2.5}
Let $u$ be a viscosity solution of
\begin{equation*}
\begin{aligned}
    \sigma_k(D^2u)&=f &&~\mbox{in}~B_1.
\end{aligned}
\end{equation*}
Assume that $u\in W^{2,p}(B_1)$ $(p>k(n-1)/2)$.

Then for any $\delta_1>0$ and $K>0$, there exist $\delta,r>0$ (depending only on $n,\delta_1$ and $K$) such that for some $k$-admissible $P\in \mathcal{P}_2$,
\begin{equation*}
\|u-P\|_{L^{\infty}(B_{r})}\leq \delta_1 r^2
\end{equation*}
and
\begin{equation}\label{e2.36}
  \sigma_k(D^2P)=1,\quad \|P\|\leq C,
\end{equation}
provided
\begin{equation*}
\|u\|_{W^{2,p}(B_1)}\leq K, \quad  \|f-1\|_{L^{\infty}(B_1)}\leq \delta,
\end{equation*}
where $C$ depends only on $n,k,p$ and $K$.
\end{lemma}
\proof Suppose not. Then there exist $\delta_1,K>0$, sequences of $u_m,f_m$ such that
\begin{equation*}
    \sigma_k(D^2u_m)=f_m~~\mbox{in}~B_1
\end{equation*}
and
\begin{equation*}
 \|u_m\|_{W^{2,p}(B_1)}\leq K, \quad  \|f_m-1\|_{L^{\infty}(B_1)}\leq \frac{1}{m}.
\end{equation*}
Moreover, for any $k$-admissible $P\in\mathcal{P}_2$ satisfying \cref{e2.36}, we have
\begin{equation}\label{e2.37}
\|u_m-P\|_{L^{\infty}(B_{r})}\geq  \delta_1 r^2,
\end{equation}
where $r$ is to be specified later.

Since $n\geq 3$ and $k\geq 2$, we have $p>k(n-1)/2>n/2$ and thus for some $\beta>0$,
\begin{equation*}
\|u_m\|_{C^{\beta}(\bar B_1)}\leq C\|u_m\|_{W^{2,p}(B_1)}\leq CK.
\end{equation*}
Then up to a subsequence, there exist $\bar{u}$ such that
\begin{equation*}
  u_m\to \bar{u}\quad\mbox{in}~~L^{\infty}(B_1), \quad \|\bar u\|_{W^{2,p}(B_1)}\leq K.
\end{equation*}
Moreover, $\bar{u}$ is a viscosity solution of
\begin{equation*}
    \sigma_k(D^2\bar u)=1~~\mbox{in}~B_1.
\end{equation*}
Hence, $\bar{u}$ is smooth (see \cite[Theorem 1.1]{MR1840289} and there exists $k$-admissible $P\in \mathcal{P}_2$ such that \cref{e2.36} holds and
\begin{equation*}
|\bar u(x) -P(x)|\leq C|x|^3,~\forall ~x\in B_{1/2}.
\end{equation*}
Take $r$ small such that $Cr=\delta_1/2$. Then
\begin{equation*}
   \|\bar u-P\|_{L^{\infty}(B_{r})}\leq \frac{1}{2} \delta_1 r^2.
\end{equation*}
By taking the limit in \cref{e2.37}, we have
\begin{equation*}
   \|\bar u-P\|_{L^{\infty}(B_{r})}\geq \delta_1 r^2,
\end{equation*}
which is contradiction.~\qed~\\

Now, we give the~\\
\noindent\textbf{Proof of \Cref{th.H}.} The proof is similar to that of the prescribed mean curvature equation. We assume $f(0)=1$ without loss of generality. For (i), let $\tilde{u}=u-P$ and then $\tilde{u}$ is a solution of
\begin{equation*}
F(D^2\tilde{u})=\tilde {f}\quad\mbox{in}~B_1,
\end{equation*}
where
\begin{equation}\label{e4.5}
F(M):=\sigma_k(M+D^2P)-1,~\forall ~M\in \mathcal{S}^n, \quad \tilde{f}:=f-1.
\end{equation}
Thus, $F\in C^{\infty}$, $F(0)=0$ and $F$ is $\rho$-uniformly elliptic with ellipticity constants $\tilde{\lambda},\tilde{\Lambda}$ and they depend only on $n,k,\lambda,\Lambda$ and $|D^2P|$.
In addition, by assumption (i),
\begin{equation*}
  \|\tilde{u}\|_{L^{\infty}(B_1)}\leq \delta,\quad
  \|\tilde{f}\|_{C^{\alpha}(0)}\leq \delta.
\end{equation*}
Therefore, by \Cref{co1.1}, $\tilde{u}\in C^{l,\alpha}(0)$ and hence $u\in C^{l,\alpha}(0)$.

Next, we prove (ii). For $r>0$, let
\begin{equation}\label{e4.19}
y=\frac{x}{r}, \quad  \tilde{u}(y)=\frac{u(x)-P(x)}{r^2}, \quad \tilde{f}(y)=f(x)-1.
\end{equation}
Then $\tilde{u}$ is a solution of
\begin{equation*}
F(D^2\tilde{u})=\tilde{f}\quad\mbox{in}~B_1.
\end{equation*}
We first choose $r$ small enough such that $\|\tilde f\|_{C^{\alpha}}$ is small. Next, take $\delta$ small
such that $\|\tilde{u}\|_{L^{\infty}}$ is small. Then the conclusion follows from \Cref{co1.1}.

Note that (iii) is a direct consequence of \Cref{th2.1}.

As regards (iv), by \Cref{le2.5}, for any $\delta_1>0$, there exist $k$-admissible $P\in \mathcal{P}_2$ and $r>0$ such that,
\begin{equation*}
\|u-P\|_{L^{\infty}(B_{r})}\leq \delta_1 r^2.
\end{equation*}
Then the conclusion follows as above through a transformation like \cref{e4.19}.

Finally, we prove (v). By the comparison estimate for $k$-Hessian estimate (see \cite[Lemma 2.1]{MR1466315}),
\begin{equation*}
\|u-P\|_{L^{\infty}(B_{\delta})}\leq C\delta^2\|f^{1/k}-1\|_{L^{\infty}(B_{\delta})}
\leq C\delta^{2+\alpha},
\end{equation*}
where $C$ depends only on $n,k$. Then the conclusion follows as above.
~\qed~\\

Next, we give the~\\
\noindent\textbf{Proof of \Cref{th.H2}.} We assume $f(0)=1$ as before. For (i), since $\partial \Omega$ and $g_0$ are smooth, there exists a solution $v\in C^{3}(\bar{\Omega})$ of (see \cite[Theorem 1.1]{MR1368245}, \cite[Theorem 3.4]{MR2500526})
\begin{equation*}
  \left\{\begin{aligned}
    \sigma_k(D^2v)&=1 &&~\mbox{in}~\Omega;\\
    v&=g_0&&~\mbox{on}~\partial\Omega.
  \end{aligned}\right.
\end{equation*}
By the comparison estimate for $k$-Hessian estimate (see \cite[Lemma 2.1]{MR1466315}),
\begin{equation}\label{e2.26}
\|u-v\|_{L^{\infty}(\Omega)}\leq \|g-g_0\|_{L^{\infty}(\partial\Omega)}+
C\mathrm{diam}(\Omega)^2\|f^{1/k}-1\|_{L^{\infty}(\Omega)}\leq C\delta^{\alpha},
\end{equation}
where $C$ depends only on $n,k$ and $\Omega$.

In addition, by the interior $C^3$ estimate for $v$, there exists a $k$-admissible $P\in \mathcal{P}_2$ such that
\begin{equation}\label{e2.27}
|v(x)-P(x)|\leq C|x|^3,~\forall ~x\in B_{r_0}
\end{equation}
and
\begin{equation}\label{e2.28}
  \sigma_k (D^2P)=1,\quad \|P\|\leq C
\end{equation}
where $r_0$ depends only on $\Omega$ and $C$ depends only on $n,k,\Omega, \|g_0\|_{C^{3,1}(\partial \Omega)}$.

Let $0<\delta_1<1$ to be specified later. Take $\delta$ and $r$ small enough such that
\begin{equation*}
Cr\leq \frac{\delta_1}{2}, \quad C\delta^{\alpha}\leq \frac{\delta_1}{2} r^2.
\end{equation*}
Then by \cref{e2.26} and \cref{e2.27},
\begin{equation}\label{e2.29}
\|u-P\|_{L^{\infty}(B_{r})}\leq \delta_1r^2.
\end{equation}
Therefore, the conclusion follows as before (cf. the proof of \Cref{th.H}).


The proof for (ii) is similar to the above. By subtracting $P$, we may assume that $|u|\leq \delta$ on $\partial \Omega$. Let $v\in C^{3}(\Omega)\cap C(\bar{\Omega})$ be the solution of (see \cite[Theorem 4.4]{MR2500526})
\begin{equation*}
  \left\{\begin{aligned}
    \sigma_k(D^2v)&=1 &&~\mbox{in}~\Omega;\\
    v&=0&&~\mbox{on}~\partial\Omega.
  \end{aligned}\right.
\end{equation*}
By the comparison estimate,
\begin{equation}\label{e2.31}
\|u-v\|_{L^{\infty}(\Omega)}\leq \|u\|_{L^{\infty}(\partial\Omega)}
+C\|f^{1/k}-1\|_{L^{\infty}(\partial\Omega)}
\leq C\delta^{\alpha}.
\end{equation}
By the interior $C^3$ estimate for $v$, there exist $r_0>0$ and a $k$-admissible $P\in \mathcal{P}_2$ such that
\begin{equation*}
|v(x)-P(x)|\leq C|x|^3,~\forall ~x\in B_{r_0}
\end{equation*}
and
\begin{equation*}
  \sigma_k (D^2P)=1,\quad \|P\|\leq C.
\end{equation*}
Then the conclusion follows as above.

For (iii), with the aid of the transformation
\begin{equation*}
  y=\frac{x}{\delta}, \quad \tilde{u}(y)=\frac{u(x)-P(x)}{\delta^2},
\end{equation*}
the conclusion follows from (ii).~\qed~\\

Since the $2$-Hessian equation has pure interior $C^{1,1}$ estimate in some circumstances, we have the following additional regularity for this equation.
\begin{theorem}\label{th.H3}
Let $0<\alpha<1$ and $u\in C(\bar{B}_1)$ be a viscosity solution of
\begin{equation}\label{e.H-3}
    \sigma_2(D^2u)=f ~~\mbox{ in}~B_1\subset \mathbb{R}^n,
\end{equation}
where $0<\lambda\leq f\leq \Lambda$. Then $u\in C^{l,\alpha}(0)$ ($l\geq 2$) provided one of the
following conditions holds :\\
(i)
\begin{equation*}
n=3,4, \quad f\in C^{l-2,\alpha}(0), \quad [f]_{C^{\alpha}(0)}\leq \delta,
\end{equation*}
where $0<\delta<1$ depends only on $n,\lambda,\Lambda,\alpha$ and $\|u\|_{L^{\infty}(B_1)}$.\\
(ii)
\begin{equation*}
n=3,4, \quad u\in C^{1,1}(0), \quad  f\in C^{l-2,\alpha}(0).
\end{equation*}
(iii) there exists $\varepsilon> 0$ such that $u+\left((n(n-1)/2)^{-1/2}-\varepsilon\right)|x|^2$ is convex and
\begin{equation*}
f\in C^{l-2,\alpha}(0), \quad [f]_{C^{\alpha}(0)}\leq \delta,
\end{equation*}
where $0<\delta<1$ depends only on $n,\lambda,\Lambda,\alpha,\varepsilon$ and $\|u\|_{L^{\infty}(B_1)}$.\\
(iv) there exists $\varepsilon> 0$ such that $u+\left((n(n-1)/2)^{-1/2}-\varepsilon\right)|x|^2$ is convex and
\begin{equation*}
u\in C^{1,1}(0), \quad  f\in C^{l-2,\alpha}(0).
\end{equation*}
\end{theorem}

\begin{remark}\label{re4.6}
The regularity under assumptions (i) and (ii) is based on the pure interior $C^{1,1}$ estimate for the $2$-Hessian equation when $f\equiv 1$. For $n=3$, it was derived by Warren and Yuan \cite{MR2487850}. For $n=4$, it was derived by Shankar and Yuan \cite{shankar2023hessian}.

The regularity under assumptions (iii) and (iv) is based on the pure interior $C^{1,1}$ estimates for almost convex viscosity solutions. Mooney \cite{MR4246798} proved the convex case; Shankar and Yuan \cite{MR4210289} proved the almost convex case. We also note that the interior $C^{1,1}$ estimate holds as well for \emph{smooth} semi-convex solutions (see \cite{MR3961212}, \cite{MR4054864}). However, we can not obtain a regularity result by an approximation argument based on these estimates. This has been pointed out in \cite[P. 2474, L. 2]{MR4246798} and \cite[P. 2]{MR4210289}.
\end{remark}

\begin{remark}\label{re4.7}
The regularity under assumptions (i), (iii) (resp. (ii), (iv)) is analogous to (ii) (resp. (iii)) in \Cref{th.MS} since we have pure interior $C^{1,1}$ estimates for constant $f$. A regularity result similar to (i) was proved by Xu \cite{MR4093797} in dimension $3$. The interior $C^{2,\alpha}$ regularity in dimension $3$ with $f\in C^{0,1}$ was obtained by Zhou \cite[Theorem 1.4]{Zhou2023}.
\end{remark}
~\\

We first give a lemma.
\begin{lemma}\label{le4.2}
Let $u\in C(\bar{B}_1)$ be a viscosity solution of
\begin{equation*}
    \sigma_2(D^2u)=f ~~\mbox{ in}~B_1.
\end{equation*}
and $\|u\|_{L^{\infty}(B_1)}\leq K$. Suppose that one of the following two conditions holds:~\\
(i) $n=3$ or $4$;\\
(ii) there exists $\varepsilon> 0$ such that $u+\left((n(n-1)/2)^{-1/2}-\varepsilon\right)|x|^2$ is convex.

Then for any $\delta_1>0$, there exist $\delta, r>0$ depending only on $n,K,\varepsilon$ and $\delta_1$ such that if $\|f-1\|_{L^{\infty}(B_1)}\leq \delta$, we have
\begin{equation*}
  \|u-P\|_{L^{\infty}(B_{r})}\leq \delta_1 r^2
\end{equation*}
and
\begin{equation}\label{e4.17}
\sigma_2(D^2P)=1, \quad \|P\|\leq C,
\end{equation}
where $P\in\mathcal{P}_2$ is $2$-admissible and $C$ depends only on $n,K,\varepsilon$.
\end{lemma}
\proof We only give the proof under the condition (ii). Suppose not. Then there exist $K,\delta_1>0$ and sequences of $u_m,f_m$ such that
\begin{equation*}
    \sigma_2(D^2u_m)=f_m ~~\mbox{ in}~B_1,
\end{equation*}
\begin{equation*}
  \|u_m\|_{L^{\infty}(B_1)}\leq K, \quad \|f_m-1\|_{L^{\infty}(B_1)}\leq 1/m
\end{equation*}
and $u_m+\left((n(n-1)/2)^{-1/2}-\varepsilon\right)|x|^2$ is convex. However, for any $2$-admissible $P\in \mathcal{P}_2$ satisfying \cref{e4.17}, we have
\begin{equation}\label{e4.16}
\|u_m-P\|_{L^{\infty}(B_{r})}> \delta_1 r^2,
\end{equation}
where $r$ is to be specified later.

By the interior H\"{o}lder regularity (see \cite[Theorem 4.1]{MR1466315}, \cite[corollary 9.1]{MR2500526}), $u_m$ are uniformly bounded and equicontinuous in any compact subset of $B_1$. Then up a subsequence, there exists $\bar u$ such that
\begin{equation*}
u_m\to \bar u\quad\mbox{in}~~L^{\infty}_{loc}(B_1).
\end{equation*}
Then $\bar u+\left((n(n-1)/2)^{-1/2}-\varepsilon\right)|x|^2$ is convex and $\bar u$ is a viscosity solution of
\begin{equation*}
  \sigma_2(D^2\bar u)=1 ~~\mbox{ in}~B_1.
\end{equation*}
Thus, $\bar u$ is smooth (see \cite[Theorem 1.1]{MR4210289}) and there exists a $2$-admissible $P\in\mathcal{P}_2$ such that \cref{e4.17} holds and
\begin{equation*}
|\bar u(x) -P(x)|\leq C|x|^3,~\forall ~x\in B_{1/2},
\end{equation*}
where $C$ depends only on $n$ and $K$. Take $r$ small such that $Cr=\delta_1/2$. Then
\begin{equation*}
   \|\bar u-P\|_{L^{\infty}(B_{r})}\leq \frac{1}{2} \delta_1 r^2.
\end{equation*}
By taking the limit in \cref{e4.16}, we have
\begin{equation*}
   \|\bar u-P\|_{L^{\infty}(B_{r})}\geq \delta_1 r^2,
\end{equation*}
which is a contradiction.~\qed~\\

\noindent\textbf{Proof of \Cref{th.H3}.} With the aid of \Cref{le4.2}, the theorem can be proved in a similar way as before and we omit it.~\qed~\\

\section{$k$-Hessian quotient equation}\label{S.HQ}
In this section, we consider the $k$-Hessian quotient equations:
\begin{equation*}
 S_{k,l}(D^2u):=\frac{\sigma_k(D^2u)}{\sigma_l(D^2u)} =f ~~\mbox{ in}~B_1,
\end{equation*}
where $1\leq l<k\leq n$. Similar to the $k$-Hessian equation, for any $k$-admissible polynomial $P\in\mathcal{P}_2$ with $S_{k,l}(D^2P)=1$, define
\begin{equation*}
F(M)=S_{k,l}(M+D^2P)-1,~\forall ~M\in \mathcal{S}^n.
\end{equation*}
Then $F\in C^{\infty}, F(0)=0$ and $F$ is $\rho$-uniformly elliptic with $\tilde{\lambda},\tilde{\Lambda}$ which depend only on $n,k,l$ and $|D^2P|$.

The main results are the following.
\begin{theorem}\label{th.HQ}
Let $1\leq l< k\leq n$, $0<\alpha<1$ and $u\in C(\bar{B}_1)$ be a viscosity solution of
\begin{equation}\label{e.HQ}
 S_{k,l}(D^2u)=f ~~\mbox{ in}~B_1,
\end{equation}
where $0<\lambda\leq f\leq \Lambda$. Then $u\in C^{m,\alpha}(0)$ ($m\geq 2$) provided one of the
following conditions holds :\\
(i) there exists a $k$-admissible $P\in \mathcal{P}_2$ such that
\begin{equation*}
\|u-P\|_{L^{\infty}(B_1)}\leq \delta, \quad S_{k,l}(D^2P)=f(0), \quad  f\in C^{m-2,\alpha}(0), \quad
[f]_{C^{\alpha}(0)}\leq \delta,
\end{equation*}
where $0<\delta<1$ depends only on $n,k,l,\lambda,\Lambda,\alpha$ and $|D^2P|$.\\
(ii) there exists a $k$-admissible $P\in \mathcal{P}_2$ such that
\begin{equation*}
\|u-P\|_{L^{\infty}(B_1)}\leq \delta, \quad S_{k,l}(D^2P)=f(0), \quad  f\in C^{m-2,\alpha}(0),
\end{equation*}
where $0<\delta<1$ depends only on $n,k,l,\lambda,\Lambda,\alpha,|D^2P|$ and $[f]_{C^{\alpha}(0)}$.\\
(iii) $u\in C^{2}(0)$ and $f\in C^{m-2,\alpha}(0)$.\\
(iv)  there exists a $k$-admissible $P\in \mathcal{P}_2$ such that
\begin{equation*}
u=P~~\mbox{ on}~\partial B_{\delta}, \quad \sigma_k(D^2P)=f(0), \quad  f\in C^{m-2,\alpha}(0),
\end{equation*}
where $\delta$ depends only on $n,k,l,\lambda,\Lambda,\alpha$ and $|D^2P|$.~\\

\end{theorem}

\begin{remark}\label{re5.1}
Since the proof is quite similar to the that for the $k$-Hessian equation, we omit it. We point out that we need to apply a comparison estimate when proving (iv) (see the proof of (v) in \Cref{th.H}). There is no existing literature to cite. Indeed, the comparison estimate for $k$-Hessian equation (see \cite[Lemma 2.2]{MR1466315}) can be extended to the $k$-Hessian quotient equation with the $L^p$ norm in the estimate replaced by the $L^{\infty}$ norm. The key is that the operator $(S_{k,l})^{1/(k-l)}$ is $1$-homogenous and concave. Then an inequality similar to \cite[(2.4)]{MR1466315} holds and the comparison estimate can be proved similarly.
\end{remark}
~\\

For the $k$-Hessian quotient equations, there are few pure interior $C^{1,1}$ estimates and Pogorelov's type estimates for \emph{smooth convex} solutions (see \cite{dai2023interior,
lu2023interior,lu2024interior}) until now. Unfortunately, as before, we can not use these to build regularity. Instead, we have the following regularity based on a priori estimates established by Trudinger \cite[Theorem 1.1]{MR1368245}. Since the proof is similar to that for the $k$-Hessian equation and we omit it.
\begin{theorem}\label{th.HQ2}
Let $1\leq l< k\leq n$, $0<\alpha<1$ and $u\in C(\bar{\Omega})$ be a viscosity solution of
\begin{equation*}
  \left\{\begin{aligned}
   S_{k,l}(D^2u)&=f &&~\mbox{in}~\Omega;\\
    u&=g&&~\mbox{on}~\partial \Omega,
  \end{aligned}\right.
\end{equation*}
where $0<\lambda\leq f\leq \Lambda$. Then $u\in C^{m,\alpha}(0)$ ($m\geq 2$) provided one of the
following conditions holds :\\
(i) $\partial \Omega\in C^{3,1}$, $\Omega$ is $(k-1)$-convex and
\begin{equation*}
\|g-g_0\|_{L^{\infty}(\partial \Omega)}\leq
\delta,~g_0\in C^{3,1}(\partial \Omega),~ f\in C^{m-2,\alpha}(0),~[f]_{C^{\alpha}(0)}\leq \delta,
\end{equation*}
where $0<\delta<1$ depends only on $n,k,l\lambda,\Lambda,\alpha,\Omega$ and
$\|g_0\|_{C^{3,1}(\partial \Omega)}$.\\
(ii)
\begin{equation*}
\Omega=B_{\delta},~~\|g-P\|_{L^{\infty}(\partial B_{\delta})}\leq \delta_ 1\delta^2,~~P\in\mathcal{P}_1,~~
 f\in C^{m-2,\alpha}(0),
\end{equation*}
where $\delta_1$ depends only on $n,k,l,\lambda,\Lambda,\alpha$ and $\delta$ depends also on $[f]_{C^{\alpha}(0)}$.~\\
\end{theorem}

\begin{remark}\label{re5.2}
Similar to the prescribed mean curvature equation, the $C^{2,\alpha}$ regularity for $k$-Hessian quotient equation can not hold unconditionally. Consider the following counterexample borrowed from \cite[P. 2]{Zhou2023}. For any $1\leq l<k\leq n$ and $0<\theta<1$, define
\begin{equation*}
u(x)=\frac{1}{2}|x'|^2+\frac{1}{1+\theta}|x_n|^{1+\theta}.
\end{equation*}
Then $u$ is a viscosity solution of \cref{e.HQ} with $f\in C^{1-\theta}$. Similar to the prescribed mean curvature equation, $[f]_{C^{1-\theta}(0)}\simeq \theta^{-1}$. Hence, if $\theta$ is smaller, $f$ is smoother but $[f]_{C^{1-\theta}(0)}$ is bigger. Correspondingly, $u$ has lower regularity.
\end{remark}
~\\

\section{Lagrangian mean curvature equation}\label{S.SL}
In this section, we consider the Lagrangian mean curvature equation (called special Lagrangian equation for constant $f$):
\begin{equation}\label{e.SL}
F(D^2u)=\sum_{i=1}^{n}\arctan\lambda_i=f\quad\mbox{in}~~B_1,
\end{equation}
where $\lambda_i$ are the eigenvalues of $D^2u$. The $f$ is called phase function. Obviously, there must hold
\begin{equation*}
-n\frac{\pi}{2}<f<n\frac{\pi}{2}.
\end{equation*}
We define
\begin{equation}\label{e6.1}
\varepsilon_f=\inf_{x\in B_1}\min (n\frac{\pi}{2}-f(x), f(x)+n\frac{\pi}{2}).
\end{equation}

Similar to the prescribed mean curvature equation, $F\in C^{\infty}$ is $1$-uniformly elliptic with $\lambda=1/5,\Lambda=1$. Moreover, for any $P\in \mathcal{P}_2$, define
\begin{equation*}
G(M)=F(M+D^2P)-F(D^2P).
\end{equation*}
Then $G(0)=0$ is $\rho$-uniformly elliptic with $\tilde{\lambda},\tilde{\Lambda}$ and they depend only on $n$ and $|D^2P|$.

The phase $f$ is divided into three categories by Yuan \cite{MR2199179}: critical ($|f|=(n-2)\pi/2$), subcritical ($|f|<(n-2)\pi/2$) and supercritical ($|f|>(n-2)\pi/2$). For critical/supercritical phases, the level set $\left\{(\lambda_1,..,\lambda_n): \sum_{i=1}^{n} \arctan \lambda_i=c\right\}$ is convex and then an extended Evans-Krylov theorem \cite{MR1793687} can be applied if we assume $u\in C^{1,1}$ a priori. For more knowledge and historic literature with respect to the special Lagrangian equation, we refer to \cite{MR2711859} and \cite{MR4181013}.

Now, we state the main results in this section. The following theorem is similar to the previous ones and we omit its proof.
\begin{theorem}\label{th.SL}
Let $0<\alpha<1$ and $u\in C(\bar{B}_1)$ be a viscosity solution of \cref{e.SL}. Then $u\in C^{k,\alpha}(0)$ ($k\geq 2$) provided one of the following conditions holds :\\
(i) there exists $P\in \mathcal{P}_2$ such that
\begin{equation*}
\|u-P\|_{L^{\infty}(B_1)}\leq \delta, \quad F(D^2P)=f(0), \quad  f\in C^{k-2,\alpha}(0), \quad
[f]_{C^{\alpha}(0)}\leq \delta,
\end{equation*}
where $0<\delta<1$ depends only on $n,\alpha$ and $|D^2P|$.\\
(ii) there exists $P\in \mathcal{P}_2$ such that
\begin{equation*}
\|u-P\|_{L^{\infty}(B_1)}\leq \delta, \quad F(D^2P)=f(0), \quad  f\in C^{k-2,\alpha}(0),
\end{equation*}
where $0<\delta<1$ depends only on $n,\alpha,|D^2P|$ and $[f]_{C^{\alpha}(0)}$.\\
(iii) $u\in C^{2}(0)$ and $f\in C^{k-2,\alpha}(0)$.\\
(iv)  there exists $P\in \mathcal{P}_2$ such that
\begin{equation*}
u=P~~\mbox{ on}~\partial B_{\delta}, \quad F(D^2P)=f(0), \quad  f\in C^{k-2,\alpha}(0),
\end{equation*}
where $\delta$ depends only on $n,\alpha$ and $|D^2P|$.~\\
\end{theorem}

\begin{remark}\label{re6.3}
We can not obtain $C^{2,\alpha}$ regularity for the special Lagrangian equation without additional assumptions besides $f\in C^{\alpha}$. If $f$ is subcritical, even a constant $f$, $C^{2,\alpha}$ regularity can not be expected. Nadirashvili and Vl\u{a}du\c{t}\cite{MR2683755} constructed counterexamples for $n\geq 3$ with solutions only belonging to $C^{1,1/3}$ for any constant subcritical phase. In fact, for any $\delta>0$, there exists a viscosity solution $u\notin C^{1,\delta}$, which was given by Wang and Yuan \cite{MR3117304}. In \cite{Mooney2023}, Mooney and Savin provided a counterexample such that the solution $u\in C^{0,1}$ but $u\notin C^1$.
\end{remark}
~\\

If $u$ is a convex viscosity solution, we have the following regularity.
\begin{theorem}\label{th.SL2}
Let $0<\alpha<1$ and $u\in C(\bar{B}_1)$ be a convex viscosity solution of \cref{e.SL}. Suppose that $f\in C^{k-2,\alpha}(0)$ ($k\geq 2$). Then $u\in C^{k,\alpha}(0)$ provided $u\in C^{1,1}(0)$ or $[f]_{C^{\alpha}(0)}\leq \delta$, where $0<\delta<1$ depends only on $n,\alpha,\varepsilon_f$ and $\|u\|_{L^{\infty}(B_1)}$.
\end{theorem}
\begin{remark}\label{re6.1}
If $f$ is a constant, \Cref{th.SL2} has been proved by Chen, Warren and Yuan \cite{MR2492708} (for smooth solutions) and Chen, Shankar and Yuan \cite{MR4655360} (for viscosity solutions). For nonconstant $f$ and viscosity solutions, Bhattacharya and Shankar proved the following regularity results:
\begin{itemize}
  \item $f\in C^{2,\alpha}\Longrightarrow u\in C^{4,\alpha}$ (see \cite{MR4649187}).
  \item $f\in C^{2}\Longrightarrow u\in C^{3,\alpha}$ (see \cite{BS2020}).
  \item $u\in C^{1,\beta},f\in C^{\alpha}\Longrightarrow u\in C^{2,\alpha}$, where $\beta>(1+\alpha)^{-1}$ (see \cite{BS2020}).
\end{itemize}
\end{remark}

\begin{remark}\label{re6.2}
Similar to the prescribed mean curvature equation, the smallness condition in \Cref{th.SL2} can not be removed. Consider the following counterexample borrowed from \cite[Remark 1.3]{MR4314140} (see also \cite[Remark 1.2]{BS2020}). Take $n=2$ and $0<\theta<1$. Define
\begin{equation*}
u(x_1,x_2)=\frac{1}{1+\theta}|x_1|^{1+\theta}+\frac{1}{2}x_2^2, \quad f(x)=\frac{3\pi}{4}-\arctan(\theta^{-1}|x_1|^{1-\theta}).
\end{equation*}
Then $u$ is a strictly convex viscosity solution and $f\in C^{1-\theta}$ is supercritical. However, $u\in C^{1,\theta}(0)$ only.

On the other hand, note that $[f]_{C^{1-\theta}(0)}\simeq \theta^{-1}$. Thus, if $\theta$ is smaller, $f$ has more smoothness but $[f]_{C^{1-\theta}(0)}$ is bigger. Correspondingly, $u$ has less regularity. This demonstrate the assertion ``smallness is more important than smoothness'' again.
\end{remark}
~\\

If the phase $f$ is critical and supercritical, we can also deduce regularity. We first consider the supercritical case.
\begin{theorem}\label{th.SL3}
Let $0<\alpha<1$ and $u\in C(\bar{B}_1)$ be a viscosity solution of \cref{e.SL}. Suppose that for some $\varepsilon>0$,
\begin{equation*}
|f|\geq (n-2)\frac{\pi}{2}+\varepsilon, \quad f\in C^{k-2,\alpha}(0)~(k\geq 2).
\end{equation*}
Then $u\in C^{k,\alpha}(0)$ provided $u\in C^{1,1}(0)$ or $[f]_{C^{\alpha}(0)}\leq \delta$, where $0<\delta<1$ depends only on $n,\alpha,\varepsilon,\varepsilon_f$ and $\|u\|_{L^{\infty}(B_1)}$.
\end{theorem}

For the critical case, we have
\begin{theorem}\label{th.SL4}
Let $0<\alpha<1$ and $u\in C(\bar{B}_1)$ be a viscosity solution of \cref{e.SL}. Suppose that
\begin{equation*}
|f|= (n-2)\frac{\pi}{2}, \quad f\in C^{k-2,\alpha}(0)~(k\geq 2).
\end{equation*}
Then $u\in C^{k,\alpha}(0)$ provided $[f]_{C^{\alpha}(0)}\leq \delta$, where $0<\delta<1$ depends only on $n,\alpha,\varepsilon_f,\|u\|_{L^{\infty}(B_1)}$ and the modulus of continuity of $u$.
\end{theorem}
\begin{remark}\label{re6.4}
If $f$ is a constant, \Cref{th.SL3} and \Cref{th.SL4} have been proved by Warren and Yuan \cite{MR2511754} ($n=2$), \cite{MR2666907} ($n=3$) and Wang and Yuan \cite{MR3188067} ($n\geq 3$).

For nonconstant $f\in C^{1,1}$, Bhattacharya \cite{MR4314140} proved the supercritical case and Lu \cite{MR4578557} proved the critical and supercritical cases. Recently, Zhou \cite{Zhou2023_2} extended above results to $f\in C^{0,1}$.

Note that the $f$ is supercritical in the counterexample in \Cref{re6.2}. Hence, the smallness condition can not be removed in \Cref{th.SL3}.
\end{remark}
~\\

Without the convexity assumption on $u$ or the critical/supercritical assumption on $f$, Yuan \cite{MR1808031} proved $C^{2,\alpha}$ regularity for $C^{1,1}$ viscosity solutions and constant phase in dimension $3$. Based on this result, we have the following regularity for a general phase $f$.
\begin{theorem}\label{th.SL5}
Let $0<\alpha<1$ and $u\in C^{1,1}(\bar{B}_1)$ be a viscosity solution of \cref{e.SL}. Suppose that $f\in C^{k-2,\alpha}(0)~(k\geq 2)$. Then $u\in C^{k,\alpha}(0)$.
\end{theorem}

We first give the~\\
\noindent\textbf{Proof of \Cref{th.SL2}.} We only give the proof for the case $[f]_{C^{\alpha}(0)}\leq \delta$ since the case $u\in C^{1,1}(0)$ can be transformed to the former case.

\textbf{Claim:} For any $\delta_1>0$, if $\delta$ small enough, there exist $r>0$ and $P\in \mathcal{P}_2$ such that
\begin{equation}\label{e6.2-1}
F(D^2P)=f(0)
\end{equation}
and
\begin{equation}\label{e6.2}
\|u-P\|_{L^{\infty}(B_r)}\leq \delta_1 r^2.
\end{equation}
The proof is similar to the previous (e.g. \Cref{le4.2}). Suppose not. Then there exist $\varepsilon_0,K>0$ and sequences of $u_m,f_m$ such that $u_m,f_m$ satisfy \cref{e.SL} and
\begin{equation*}
\|u_m\|_{L^{\infty}(B_1)}\leq K, \quad [f_m]_{C^{\alpha}(0)}\leq 1/m, \quad \varepsilon_{f_m}\geq \varepsilon_0.
\end{equation*}
In addition, for any $P\in \mathcal{P}_2$ with \cref{e6.2-1} holding for $f_m$, we have
\begin{equation}\label{e6.4}
\|u_m-P\|_{L^{\infty}(B_r)}> \delta_1 r^2,
\end{equation}
where $0<r<1$ is to be specified later.

Since $u_m$ are convex,
\begin{equation*}
\|u_m\|_{C^{0,1}(\bar{B}_{1/2})}\leq C\|u_m\|_{L^{\infty}(B_1)}\leq K.
\end{equation*}
Hence, $u_m$ are uniformly bounded and equicontinuous. Up to a subsequence (similarly in the following argument), there exists $\bar{u}$ such that
\begin{equation*}
u_m\to \bar{u}\quad\mbox{in}~~L^{\infty}(B_{1/2}).
\end{equation*}
In addition, since $\varepsilon_{f_m}\geq \varepsilon_0$,
\begin{equation*}
f_m(0)\to f_0\in \left[-n\pi/2+\varepsilon_0, n\pi/2-\varepsilon_0\right].
\end{equation*}
By combining with $[f_m]_{C^{\alpha}(0)}\to 0$, we conclude that $\bar{u}$ is a viscosity solution of
\begin{equation*}
F(D^2\bar{u})=f_0\quad\mbox{in}~~B_{1/2}.
\end{equation*}

From the regularity for constant phases (see \cite{MR4655360}), there exists $\bar P\in\mathcal{P}_2$ such that $F(D^2\bar P)=f_0$ and
\begin{equation*}
|\bar u(x) -\bar P(x)|\leq C|x|^3,~\forall ~x\in B_{1/2},
\end{equation*}
where $C$ depends only on $n$ and $K$. Take $r$ small such that $Cr=\delta_1/2$. Then
\begin{equation*}
   \|\bar u-\bar P\|_{L^{\infty}(B_{r})}\leq \frac{1}{2} \delta_1 r^2.
\end{equation*}
In addition, since $f_m(0)\to f_0$, we can choose
\begin{equation*}
P_m(x)=\bar P(x)+\sum_{i=1}^{n} c_{i,m}x_i^2
\end{equation*}
such that $c_{i,m}\to 0$ and \cref{e6.2-1} holds for $P_m$ and $f_m(0)$. Thus, \cref{e6.4} holds for $P_m$. By taking the limit in \cref{e6.4}, we have
\begin{equation*}
   \|\bar u-\bar P\|_{L^{\infty}(B_{r})}\geq \delta_1 r^2,
\end{equation*}
which is a contradiction. Therefore, the claim holds.

Once the claim is proved, the regularity for $u$ follows as before.~\qed~\\

Next, we give the~\\
\noindent\textbf{Proof of \Cref{th.SL3}.} We only consider the case
\begin{equation*}
f>(n-2)\frac{\pi}{2}+\varepsilon, \quad [f]_{C^{\alpha}(0)}\leq \delta.
\end{equation*}

We need the following claim as before:~\\
\textbf{Claim:} For any $\delta_1>0$, if $\delta$ small enough, there exist $r>0$ and $P\in \mathcal{P}_2$ such that
\begin{equation*}
F(D^2P)=f(0)
\end{equation*}
and
\begin{equation*}
\|u-P\|_{L^{\infty}(B_r)}\leq \delta_1 r^2.
\end{equation*}

The proof is similar to the proof of \Cref{th.SL3}. We only need to take care of the compactness of solutions. Indeed, since $f>(n-2)\frac{\pi}{2}+\varepsilon$, the $u+K|x|^2$ is convex where $K=\tan(\pi/2 -\varepsilon)$. This assertion can be proved directly by the definition of viscosity solution. Then
\begin{equation*}
\|u\|_{C^{0,1}(\bar{B}_{1/2})}\leq C_1\|u\|_{L^{\infty}(B_1)}+C_2,
\end{equation*}
where $C_1$ depends only on $n$ and $C_2$ depends also on $\varepsilon$.

Once the compactness of solutions is built, with the aid of regularity for constant supercritical phase (see \cite{MR3188067}), the rest proof is quite similar to that of \Cref{th.SL3} and we omit it. ~\qed~\\

Next, we give the~\\
\noindent\textbf{Proof of \Cref{th.SL4}.} The difficulty lies in that we can not derive the compactness of solutions and we have to rely on the modulus of continuity of $u$ directly. Hence, the constant $\delta$ depends on the modulus of continuity of $u$.

Then the following claim can be proved as above and we omit it.~\\
\textbf{Claim:} For any $\delta_1>0$, if $\delta$ small enough, there exist $r>0$ and $P\in \mathcal{P}_2$ such that
\begin{equation*}
F(D^2P)=f(0)
\end{equation*}
and
\begin{equation*}
\|u-P\|_{L^{\infty}(B_r)}\leq \delta_1 r^2.
\end{equation*}
Then with the aid of regularity for constant critical phase (see \cite{MR3188067}), the regularity of $u$ follows as before. ~\qed~\\

Next, we give the~\\
\noindent\textbf{Proof of \Cref{th.SL5}.} Since $u\in C^{1,1}(\bar{B}_1)$, we have the compactness of solutions a priori. Then with the aid of regularity for constant critical phase (see \cite{MR1808031}), the regularity of $u$ follows as before. ~\qed~\\

\begin{remark}\label{re6.5}
Note we have assumed $u\in C^{1,1}$. Thus, the equation is in fact uniformly elliptic. In addition, the regularity for constant critical phase has been proved in \cite{MR1808031}. Therefore, the $C^{2,\alpha}$ regularity can be proved directly by the theory for uniformly elliptic equations (cf. \cite[Chapter 8]{MR1351007}). However, for higher order pointwise $C^{k,\alpha}$ regularity, we must apply the regularity theory presented in the introduction.
\end{remark}
~\\

%

%

\section{ABP maximum principle, weak Harnack inequality, H\"{o}lder
regularity}\label{S10}
In this section, we develop the basic theory for fully nonlinear locally uniformly elliptic equations. We follows
almost exactly the strategy of \cite[Chapters 3, 4]{MR1391943}. The main difficulty is that we can not make the scaling
argument arbitrarily.

First, we introduce some notions.

We also introduce the Pucci's class as follows.
\begin{definition}\label{de10.3}
We say that $u\in \bar{S}_{\rho}(\lambda,\Lambda,b_0,f)$ if for any $\varphi\in C^2(\Omega)$ with
\cref{e1.3},
\begin{equation*}
\mathcal{M}^{-}(D^2\varphi(x_0),\lambda,\Lambda)-b_0|D\varphi(x_0)|\leq f(x_0).
\end{equation*}
Similarly, we denote
$u\in \underline{S}_{\rho}(\lambda,\Lambda,b_0,f)$ if for any $\varphi\in C^2(\Omega)$ with
\begin{equation*}
\begin{aligned}
&\|\varphi\|_{C^{1,1}(\bar\Omega)}\leq \rho, \quad  \varphi(x_0)=u(x_0),\quad \varphi\geq u~~\mbox{ in}~\Omega,\\
\end{aligned}
\end{equation*}
we have
\begin{equation*}
\mathcal{M}^{+}(D^2\varphi(x_0),\lambda,\Lambda)+b_0|D\varphi(x_0)|\geq f(x_0).
\end{equation*}

We also define
\begin{equation*}
  S^*_{\rho}(\lambda,\Lambda,b_0, f)=\underline{S}_{\rho}(\lambda,\Lambda,b_0,-|f|)\cap
  \bar{S}_{\rho}(\lambda,\Lambda,b_0,|f|).
\end{equation*}
 We will denote $\underline{S}_{\rho}(\lambda,\Lambda,b_0,f)$
 ($\bar{S}_{\rho}(\lambda,\Lambda,b_0,f)$, $S^*_{\rho}(\lambda,\Lambda,f)$) by
 $\underline{S}_{\rho}(f)$ ($\bar{S}_{\rho}(f)$, $S^*_{\rho}(f)$) for short if $\lambda,\Lambda,b_0$ are
 understood well.
\end{definition}
\begin{remark}\label{re10.1}
Note that if $\rho_1\geq \rho$,
\begin{equation*}
\bar{S}_{\rho_1}(\lambda,\Lambda,b_0,f)\subset \bar{S}_{\rho}(\lambda,\Lambda,b_0,f).
\end{equation*}
Hence, for any $\rho>0$,
\begin{equation*}
\bar{S}(\lambda,\Lambda,b_0,f)\subset \bar{S}_{\rho}(\lambda,\Lambda,b_0,f),
\end{equation*}
where $\bar{S}$ denotes the usual Pucci's class.
\end{remark}

\begin{remark}\label{re1.5}
If $u\in \underline{S}_{\rho}(\lambda,\Lambda,b_0,f)$, then $-u\in
\bar{S}_{\rho}(\lambda,\Lambda,b_0,f)$. Hence, we only consider the supersolution in the following
argument.
\end{remark}
~\\

As usual, any viscosity solution belongs to the Pucci's class.

\begin{proposition}\label{pr10.1}
Let $u$ be a viscosity supersolution of
\begin{equation*}
  F(D^2u,Du,u,x)=f~~\mbox{ in}~\Omega,
\end{equation*}
where $F$ is $5\rho$-uniformly elliptic and $\|u\|_{L^{\infty}(\Omega)}\leq \rho$. Let $\phi\in
C^2(\Omega)$ with \cref{e1.3}. Then\\
(i) if \cref{e.st.1} holds,
\begin{equation*}
u-\phi\in \bar{S}_{4\rho}(\lambda,\Lambda,b_0,f(x)+c_0|u(x)-\phi(x)|-F(D^2\phi(x),D\phi(x),\phi(x), x));
\end{equation*}
(ii) if for any $|M|,|p|,|q|,|s|<\rho$ and $x\in B_1$,
\begin{equation}\label{e.st.3}
    -b_0|p-q|\leq F(M,p,s,x)-F(M,q,s,x)\leq b_0|p-q|, ~
\end{equation}
then
\begin{equation*}
u-\phi\in \bar{S}_{4\rho}(\lambda,\Lambda,b_0,f(x)-F(D^2\phi(x),D\phi(x),u(x), x));
\end{equation*}
(iii) if $F(0,p,s,x)\equiv 0$,
\begin{equation*}
  u\in \bar{S}_{4\rho}(\lambda,\Lambda,0,f-\mathcal{M}^-(D^2\phi)).
\end{equation*}
\end{proposition}
\proof Clearly, $\|u-\phi\|_{L^{\infty}(\Omega)}\leq 2\rho$. Given $x_0\in \Omega$, for any $\varphi\in
C^2(\Omega)$ with \cref{e1.3} (replacing $u$ by $u-\phi$ and $\rho$ by $4\rho$ there), $\varphi+\phi$
satisfies \cref{e1.3} (replacing $\rho$ by $5\rho$). By the definition of viscosity solution
\begin{equation*}
  \begin{aligned}
f(x_0)\geq& F(D^2\varphi(x_0)+D^2\phi(x_0),D\varphi(x_0)+D\phi(x_0), \varphi(x_0)+\phi(x_0), x_0).
  \end{aligned}
\end{equation*}
If \cref{e.st.1} holds,
\begin{equation*}
  \begin{aligned}
f(x_0)\geq &\mathcal{M}^-(D^2\varphi(x_0))-b_0|D\varphi(x_0)|-c_0|\varphi(x_0)|
+F(D^2\phi(x_0),D\phi(x_0),\phi(x_0), x_0),\\
  \end{aligned}
\end{equation*}
which means
\begin{equation*}
u-\phi\in \overline{S}_{4\rho}(\lambda,\Lambda,b_0,f(x)+c_0|u(x)-\phi(x)|-F(D^2\phi(x),D\phi(x),\phi(x),
x)).
\end{equation*}
If \cref{e.st.3} holds,
\begin{equation*}
  \begin{aligned}
f(x_0)\geq
&\mathcal{M}^-(D^2\varphi(x_0))-b_0|D\varphi(x_0)|+F(D^2\phi(x_0),D\phi(x_0),\varphi(x_0)+\phi(x_0),
x_0).
  \end{aligned}
\end{equation*}
By combining with $u(x_0)=\varphi(x_0)+\phi(x_0)$,
\begin{equation*}
u-\phi\in \overline{S}_{4\rho}(\lambda,\Lambda,b_0,f(x)-F(D^2\phi(x),D\phi(x),u(x), x)).
\end{equation*}
If $F(0,p,s,x)\equiv 0$,
\begin{equation*}
  \begin{aligned}
f(x_0)
\geq &
\mathcal{M}^-(D^2\varphi(x_0))+F(D^2\phi(x_0),D\varphi(x_0)+D\phi(x_0),\varphi(x_0)+\phi(x_0),x_0)\\
\geq &\mathcal{M}^-(D^2\varphi(x_0))+\mathcal{M}^-(D^2\phi(x_0)).
  \end{aligned}
\end{equation*}
Hence,
\begin{equation*}
u\in \bar{S}_{4\rho}(\lambda,\Lambda,0,f-\mathcal{M}^-(D^2\phi)).
\end{equation*}
~\qed~\\

As usual, we have the following maximum principle.
\begin{lemma}\label{le2.1}
If $u\in C(\bar{\Omega})$ satisfies
\begin{equation*}
u\in \bar{S}_{\rho}(\lambda,\Lambda,b_0,0)~~\mbox{ in}~\Omega
\end{equation*}
and $u\geq 0$ on $\partial \Omega$, then
\begin{equation*}
  u\geq 0~~\mbox{ in}~\Omega.
\end{equation*}
\end{lemma}
\proof Suppose not. Then by choosing $\alpha$ large enough and $\varepsilon$ small enough, $\varphi:=
\varepsilon e^{-\alpha x_1}+c$ will touch $u$ by below at some $x_0\in \Omega$ and
\begin{equation*}
\mathcal{M}^{-}(D^2\varphi(x_0))-b_0|D\varphi(x_0)|>0,
\end{equation*}
which is a contradiction.~\qed~\\

We also have
\begin{lemma}\label{le2.2}
Suppose that $\rho\geq \rho_0$,
\begin{equation*}
u\in \bar{S}_{\rho}(\lambda,\Lambda,b_0,f)~~\mbox{ in}~B_1
\end{equation*}
and
\begin{equation*}
u\geq 0~~\mbox{ on}~\partial B_1,\quad
\|f\|_{L^{\infty}(B_1)}\leq 1,
\end{equation*}
where $\rho_0>1$ is universal. Then
\begin{equation*}
  \sup_{B_1} u^-\leq C \|f\|_{L^{\infty}(B_1)},
\end{equation*}
where $C$ is universal.
\end{lemma}
\proof Take $\alpha>0$ (universal) large enough such that $v:=\|f\|_{L^{\infty}(B_1)}e^{\alpha x_1}$ is a
classical solution of
\begin{equation}\label{e2.1}
  \mathcal{M}^{-}(D^2v)-b_0|Dv|\geq 2\|f\|_{L^{\infty}(B_1)}.
\end{equation}
If
\begin{equation*}
\sup_{B_1} u^-\geq 2\|f\|_{L^{\infty}(B_1)}e^{\alpha},
\end{equation*}
$v-u$ has a local maximum at some $x_0\in B_1$. Then by the definition of viscosity solution,
\begin{equation*}
  \mathcal{M}^{-}(D^2v(x_0))-b_0|Dv(x_0)|\leq f(x_0)\leq \|f\|_{L^{\infty}(B_1)},
\end{equation*}
which contradicts to \cref{e2.1}.~\qed~\\

Next, we prove the fundamental Alexandrov-Bakel'man-Pucci maximum principle by the same way as in
\cite[Chapter 3]{MR1351007}. First, we prove a lemma analogous to \cite[Lemma 3.3]{MR1351007}:
\begin{lemma}\label{le10.3}
Let $0<\delta<1$, $\rho\geq \rho_1$ and
\begin{equation*}
u\in \overline{S}_{\rho}(\lambda,\Lambda,b_0,f)~~\mbox{ in}~B_{\delta},
\end{equation*}
where $\rho_1>1$ is universal. Suppose that
\begin{equation*}
\|f\|_{L^{\infty}(B_{\delta})}\leq 1.
\end{equation*}
Assume that $\varphi$ is a convex function in $B_{\delta}$ such that $0\leq \varphi\leq u$ in $B_{\delta}$
and $0=\varphi(0)=u(0)$. Then
\begin{equation*}
\varphi(x)\leq C(\sup_{B_{\delta}} f^+)|x|^2,~\forall ~x\in B_{\nu\delta},
\end{equation*}
where $0<\nu<1$ and $C$ are universal.
\end{lemma}
\proof We prove the lemma in the same way as that of Lemma 3.3 in \cite[Chapter 3]{MR1351007}. The
constant $\rho_1$ is to be specified later. As in the proof Lemma 3.3 of \cite{MR1351007}, for
$0<r<\delta/4$, define
\begin{equation*}
\bar{C}=\frac{1}{r^2}\sup_{B_r}\varphi.
\end{equation*}
We aim to prove
\begin{equation*}
\bar{C}\leq \frac{17}{\lambda} \sup_{B_{\delta}}f^+.
\end{equation*}
Suppose not. Let
\begin{equation}\label{e10.1}
\tilde{C}=\frac{17}{\lambda} \sup_{B_{\delta}}f^+\leq \frac{17}{\lambda}.
\end{equation}
Recall (3.3) in \cite{MR1351007}, i.e.,
\begin{equation*}
\varphi\geq \bar C r^2~~\mbox{  in}~H\cap B_{\delta},
\end{equation*}
where $H$ is a hyperplane tangent to $B_r$ at some $x_0\in \partial B_r$. Since $\tilde{C}< \bar{C}$,
\begin{equation*}
\varphi\geq \tilde C r^2~~\mbox{  in}~H\cap B_{\delta}.
\end{equation*}

Construct the polynomial $P$ as in the proof of \cite[Lemma 3.3]{MR1351007}, i.e.,
\begin{equation*}
  P(x):=\frac{\tilde{C}}{8}(x_n+r)^2-4\tilde{C}\frac{r^2}{\delta^2}|x'|^2.
\end{equation*}
Then $P+c$ for an appropriate constant $c$ will touch $u$ by below at some point $x_0$. Note that
\begin{equation*}
\|D^2P\|_{L^{\infty}}\leq \tilde{C},\quad \|DP\|_{L^{\infty}}\leq \tilde{C},\quad
\|P+c\|_{L^{\infty}}\leq \tilde{C}+\rho/2.
\end{equation*}
Choose $\rho_1\geq 2\tilde{C}$. Then by the definition of viscosity solution,
\begin{equation*}
\mathcal{M}^-(D^2P)-b_0|DP(x_0)|\leq f(x_0)\leq \sup_{B_{\delta}}f^+.
\end{equation*}
By choosing
\begin{equation*}
\nu=\min\left(\frac{1}{8}\sqrt{\frac{\lambda}{(n-1)\Lambda}},\frac{\lambda}{12b_0} \right),
\end{equation*}
we have
\begin{equation*}
  \frac{\lambda \tilde{C}}{16}\leq \mathcal{M}^-(D^2P)-b_0|DP(x_0)|\leq \sup_{B_{\delta}}f^+.
\end{equation*}
Hence,
\begin{equation*}
  \tilde{C}\leq \frac{16}{\lambda} \sup_{B_{\delta}}f^+,
\end{equation*}
which contradicts to \cref{e10.1}.~\qed~\\

We have another lemma.
\begin{lemma}\label{le2.3}
Let $\rho\geq \rho_2$,
\begin{equation*}
u\in \bar{S}_{\rho}(\lambda,\Lambda,b_0,f)~~\mbox{ in}~B_1
\end{equation*}
and
\begin{equation*}
u\geq 0~~\mbox{ on}~\partial B_1,\quad  \|u\|_{L^{\infty}(B_1)}\leq \frac{\rho}{4},\quad
\|f\|_{L^{\infty}(B_1)}\leq 1,
\end{equation*}
where $\rho_2>\max(\rho_0,\rho_1)$ is universal. Then $\Gamma_u\in C^{1,1}(\bar{B}_1)$ and
\begin{equation*}
 \mathcal{M}^-(D^2\Gamma_u)-b_0|D\Gamma_u|\leq f~~a.e. \mbox{ in}~\left\{u=\Gamma_u\right\},
\end{equation*}
where we have extended $u$ by zero outside $B_1$ and $\Gamma_u$ is the convex envelop of $-u^-$ in
$B_2$.
\end{lemma}
\proof For any $x_0\in \left\{u=\Gamma_u\right\}\subset B_1$, let $L_{x_0}$ be the supporting affine
function of $\Gamma_{u}$ at $x_0$. Then by \Cref{le2.2} and noting $x_0\in B_1$,
\begin{equation*}
\|\Gamma_u\|_{L^{\infty}(B_{1})}=\sup_{B_1}u^-\leq C_0,\quad   |DL_{x_0}|\leq
C_0,\quad \|L_{x_0}\|_{L^{\infty}(B_1)}
\leq C_0,
\end{equation*}
where $C_0$ is universal. Set $v=u-L_{x_0}$. Then $v(x_0)=(\Gamma_u-L_{x_0})(x_0)=0 $ and by
choosing $\rho_2$ large enough,
\begin{equation*}
  v\in \bar{S}_{3\rho/4}(\lambda,\Lambda,b_0,f+b_0C_0).
\end{equation*}
Let $w=v/(1+b_0C_0)$. Then
\begin{equation*}
  v\in \bar{S}_{3\rho/(4+4b_0C_0)}(\lambda,\Lambda,b_0,1).
\end{equation*}
By taking $\rho_2$ large enough, we can apply \Cref{le10.3} for any $0<\delta<1$. Then
\begin{equation*}
L_{x_0}\leq \Gamma_u\leq L_{x_0}+C|x-x_0|^2,
\end{equation*}
where $C$ is universal. Therefore, by the same argument in \cite{MR1351007} (see the proof of Lemma
3.5), $\Gamma_u\in C^{1,1}(\bar{B}_1)$. Hence, $\Gamma_u$ is second order differentiable almost
everywhere.

Next, take $x_0\in \left\{u=\Gamma_u\right\}$ such that $\Gamma_u$ is second order differentiable at
$x_0$. Let $P$ denote the second order polynomial corresponding to $\Gamma_u$ at $x_0$. Then for any
$\varepsilon>0$, $P-
\varepsilon|x-x_0|^2$ will touch locally $\Gamma_u$ and hence $u$ by below at $x_0$. Note that $\|P\|\leq C$ for some universal $C$. Then by choosing $\rho_2$ large enough, $P-\varepsilon|x-x_0|^2$ is an admissible
test function. Hence,
\begin{equation*}
\mathcal{M}^-(D^2P-2\varepsilon I)-b_0|DP(x_0)-2\varepsilon (x-x_0)|\leq f(x_0).
\end{equation*}
By letting $\varepsilon\to 0$, we arrive at the conclusion.~\qed~\\

Based on above lemma, we have the following Alexandrov-Bakel'man-Pucci maximum principle analogous to
\cite[Theorem 3.2]{MR1351007}:
\begin{theorem}[\textbf{ABP}]\label{th10.2}
Let $\rho\geq \rho_2$,
\begin{equation*}
u\in \bar{S}_{\rho}(\lambda,\Lambda,b_0,f)~~\mbox{ in}~B_1
\end{equation*}
and
\begin{equation*}
u\geq 0~~\mbox{ on}~\partial B_1,\quad   \|u\|_{L^{\infty}(B_1)}\leq \frac{\rho}{4},\quad
\|f\|_{L^{\infty}(B_1)}\leq 1,
\end{equation*}
where $\rho_2$ is as in \Cref{le2.3}. Then
\begin{equation*}
  \sup_{B_1} u^-\leq C \|f^+\|_{L^n(B_1\cap \{u=\Gamma_u\})},
\end{equation*}
where $C$ is universal.
\end{theorem}
\proof The proof is standard and we omit it (see \cite[Chapt. 9.1]{MR1814364}).  ~\qed~\\

By scaling, we have
\begin{corollary}\label{co10.1}
Let $\rho\geq \rho_2$, $K>0$,
\begin{equation*}
u\in \bar{S}_{K\rho}(\lambda,\Lambda,b_0,f)~~\mbox{ in}~B_1
\end{equation*}
and
\begin{equation*}
u\geq 0~~\mbox{ on}~\partial B_1,\quad  \|u\|_{L^{\infty}(B_1)}\leq
\frac{K\rho}{4},\quad \|f\|_{L^{\infty}(B_1)}\leq K,
\end{equation*}
where $\rho_2$ is as in \Cref{le2.3}. Then
\begin{equation*}
  \sup_{B_1} u^-\leq C \|f^+\|_{L^n(B_1\cap \{u=\Gamma_u\})},
\end{equation*}
where $C$ is universal.
\end{corollary}
\proof Let
\begin{equation*}
\tilde{u}=\frac{u}{K},\quad \tilde{f}=\frac{f}{K}.
\end{equation*}
Then
\begin{equation*}
\tilde u\in \bar{S}_{\rho}(\lambda,\Lambda,b_0,\tilde f)~~\mbox{ in}~B_1,\quad
 \|\tilde u\|_{L^{\infty}(B_1)}\leq \frac{\rho}{4},\quad
\|\tilde f\|_{L^{\infty}(B_1)}\leq 1.
\end{equation*}
By \Cref{th10.2},
\begin{equation*}
  \sup_{B_1} \tilde{u}^-\leq C \|\tilde{f}^+\|_{L^n(B_1\cap \{\tilde u=\Gamma_{\tilde u}\})}.
\end{equation*}
By transforming to $u$, we obtain the conclusion.~\qed~\\

With the aid of the ABP maximum principle, we can prove the following lemma analogous to \cite[Lemma
4.5]{MR1351007}:
\begin{lemma}\label{le10.1}
Let $\rho\geq \rho_3$,
\begin{equation*}
u\in \bar{S}_{\rho}(\lambda,\Lambda,b_0,|f|)~~\mbox{ in}~Q_{4\sqrt{n}}
\end{equation*}
and
\begin{equation*}
u\geq 0~~\mbox{ in}~Q_{4\sqrt{n}},\quad \inf_{Q_3}u\leq 1,\quad
\|u\|_{L^{\infty}(Q_{4\sqrt{n}})}\leq \frac{\rho}{8},\quad
\|f\|_{L^{\infty}(Q_{4\sqrt{n}})}\leq \varepsilon_0,
\end{equation*}
where $\rho_3\geq \rho_2$ is universal. Then
\begin{equation}\label{e10.3}
|\{u<M\}\cap Q_1|>\mu,
\end{equation}
where $0<\varepsilon_0,\mu<1$ and $M>1$ are universal.
\end{lemma}
\proof The proof is the same as that of \cite[Lemma 4.5]{MR1351007}. Construct the auxiliary function
$\varphi$ as in \cite[Lemma 4.1]{MR1351007} such that the conclusion of \cite[Lemma 4.1]{MR1351007}
holds with
\begin{equation*}
\mathcal{M}^+(D^2\varphi,\lambda,\Lambda)\leq C\xi
\end{equation*}
replaced by
\begin{equation*}
\mathcal{M}^+(D^2\varphi,\lambda,\Lambda)+b_0|D\varphi|\leq C\xi.
\end{equation*}
Since
\begin{equation*}
\|D^2\varphi\|_{L^{\infty}}\leq C_0,\quad \|D\varphi\|_{L^{\infty}}\leq
C_0,\quad \|D^2\varphi\|_{L^{\infty}}\leq C_0,
\end{equation*}
where $C_0$ is universal. By taking $\rho_3$ large enough,
\begin{equation*}
  w:=u+\varphi\in \bar{S}_{3\rho/4}(\lambda,\Lambda,b_0,|f|+C\xi),\quad  \|w\|_{L^{\infty}}\leq
  \frac{3\rho}{16},\quad
\end{equation*}
Now, we can apply the ABP maximum principle \Cref{co10.1} to obtain \cref{e10.3} as in
\cite{MR1351007}. ~\qed~\\

By iteration, we have the following lemma analogous to \cite[Lemma 4.6]{MR1351007}:
\begin{lemma}\label{le10.2}
Let $\rho\geq \rho_3$,
\begin{equation*}
u\in \bar{S}_{\rho}(\lambda,\Lambda,b_0,f)~~\mbox{ in}~Q_{4\sqrt{n}}
\end{equation*}
and
\begin{equation*}
u\geq 0~~\mbox{ in}~Q_{4\sqrt{n}},\quad \inf_{Q_3}u\leq 1,\quad
\|u\|_{L^{\infty}(Q_{4\sqrt{n}})}\leq \frac{\rho}{8},\quad
\|f\|_{L^{\infty}(Q_{4\sqrt{n}})}\leq \varepsilon_0.
\end{equation*}
Then
\begin{equation}\label{e10.4}
|\{u\geq t\}\cap Q_1|\leq Ct^{-\varepsilon},~\forall ~t<\frac{\rho}{\rho_3},
\end{equation}
where $C$ and $0<\varepsilon<1$ are universal.
\end{lemma}
\proof The proof is the same as in \cite[Lemma 4.6]{MR1351007}. We only need to prove
\begin{equation}\label{e10.2}
|\{u>M^k\}\cap Q_1|\leq (1-\mu)^k,
~\forall ~1\leq k\leq \frac{1}{\ln M}\ln \frac{\rho}{\rho_3}
\end{equation}
since \cref{e10.4} follows from \cref{e10.2} by choosing $\varepsilon>0$ with
\begin{equation*}
  1-\mu=M^{-\varepsilon}.
\end{equation*}

For $k=1$, we have just proved in \Cref{le10.1}. Suppose that \cref{e10.2} holds for $k-1$. We use the
same scaling argument as in \cite{MR1351007}:
\begin{equation*}
  \tilde{u}(y):=\frac{u(x)}{M^{k-1}}.
\end{equation*}
Since $u\in \bar{S}_{\rho}(\lambda,\Lambda,b_0,f)$,
\begin{equation*}
\tilde{u}\in \bar{S}_{\rho/M^{k-1}}(\lambda,\Lambda,b_0,\tilde f),\quad
\tilde{f}(y)=\frac{f(x)}{2^{2i}M^{k-1}}.
\end{equation*}
Note that
\begin{equation*}
k\leq \frac{1}{\ln M}\ln \frac{\rho}{\rho_3}
\end{equation*}
implies
\begin{equation*}
 \frac{\rho}{M^{k-1}}\geq \rho_3.
\end{equation*}
Therefore, by \Cref{le10.1},
\begin{equation*}
|\{\tilde u>M\}\cap Q_1|\leq 1-\mu.
\end{equation*}
By transforming to $u$, we obtain \cref{e10.2}. ~\qed~\\

Next, we prove the ``weak Harnack inequality'' analogous to \cite[Theorem 4.8]{MR1351007}:
\begin{theorem}\label{th10.3}
Let $\rho\geq \rho_3$,
\begin{equation*}
u\in \bar{S}_{\rho}(\lambda,\Lambda,b_0,f)~~\mbox{ in}~B_1
\end{equation*}
and
\begin{equation*}
0\leq u\leq 1~~\mbox{ in}~B_1,\quad
\|u\|_{L^{\infty}(B_1)}\leq \frac{\rho}{8},\quad
\|f\|_{L^{\infty}(B_1)}\leq \varepsilon_0.
\end{equation*}
Then
\begin{equation*}
\|u_{\rho}\|_{L^{p_0}(B_{1/2})}\leq C\left(\inf_{B_{1/2}} u+\|f\|_{L^{\infty}(B_1)}\right),
\end{equation*}
where $C$ and $0<p_0<1$ are universal, and
\begin{equation*}
u_{\rho}(x)=\left\{
  \begin{aligned}
&u(x),&&~~\mbox{ if}~u(x)\leq \frac{\rho}{\rho_3};\\
&0,&&~~\mbox{ if}~u(x)>\frac{\rho}{\rho_3}.
  \end{aligned}
  \right.
\end{equation*}
\end{theorem}
\proof If
\begin{equation*}
\inf_{B_{1/2}} u+\varepsilon_0^{-1}\|f\|_{L^{\infty}(B_1)}\geq 1,
\end{equation*}
by \Cref{le10.2} and taking $p_0=\varepsilon/2$,
\begin{equation*}
  \begin{aligned}
\|u_{\rho}\|_{L^{p_0}(B_{1/2})}
=&p_0\int_{0}^{+\infty} t^{p_0-1}|\{x\in B_{1/2}: u_{\rho}(x)>t\}|dt\\
\leq & |B_{1/2}|+p_0\int_{1}^{\rho/\rho_3} t^{p_0-1}|\{x\in B_{1/2}: u(x)>t\}|dt\\
\leq &C+C\int_{1}^{\rho/\rho_3} t^{p_0-1}t^{-\varepsilon}dt\\
\leq &C\\
\leq &C\left(\inf_{B_{1/2}} u+\|f\|_{L^{\infty}(B_1)}\right).
  \end{aligned}
\end{equation*}

If
\begin{equation*}
K^{-1}:=\inf_{B_{1/2}} u+\varepsilon_0^{-1}\|f\|_{L^{\infty}(B_1)}< 1,
\end{equation*}
consider
\begin{equation*}
  \tilde{u}=Ku.
\end{equation*}
Then
\begin{equation*}
\tilde{u}\in \bar{S}_{K\rho}(\lambda,\Lambda,b_0,\tilde f)
\end{equation*}
and
\begin{equation*}
\inf_{B_{1/2}}\tilde u\leq 1,\quad
\|\tilde u\|_{L^{\infty}(B_1)}\leq \frac{\rho}{8},\quad
\|\tilde f\|_{L^{\infty}(B_1)}\leq \varepsilon_0.
\end{equation*}
Then by \Cref{le10.2},
\begin{equation*}
\|\tilde u_{K\rho}\|_{L^{p_0}(B_{1/2})}\leq C.
\end{equation*}
By transferring to $u$, we obtain the conclusion. \qed~\\

Next, we prove the H\"{o}lder regularity.
\begin{theorem}\label{th10.4}
Let $\rho\geq 2\rho_3$,
\begin{equation*}
u\in S^{*}_{\rho}(\lambda,\Lambda,b_0,f)~~\mbox{ in}~B_1
\end{equation*}
and
\begin{equation*}
\|u\|_{L^{\infty}(B_1)}\leq 1,\quad
\|f\|_{L^{\infty}(B_1)}\leq \varepsilon_1,
\end{equation*}
where $0<\varepsilon_1<1$ is universal. Then
\begin{equation*}
\underset{B_r}{\mathrm{osc}}~u\leq Cr^{\alpha},~\forall ~\sqrt{\frac{2\rho_3}{\rho}}\leq r\leq 1,
\end{equation*}
where  $0<\alpha<1$ and $C$ are universal.
\end{theorem}
\proof We only need to prove
\begin{equation}\label{e10.5}
\underset{B_{1/2^{k}}}{\mathrm{osc}}~u\leq 2(1-\mu)^{k},~\forall ~0\leq k\leq \frac{1}{2}\log_2\frac{\rho}{2\rho_3},
\end{equation}
where $0<\mu<1$ is universal. We prove it by induction. For $k=0$, \cref{e10.5} holds clearly. Suppose
that it holds for $k$. For $r>0$, denote
\begin{equation*}
M_r=\sup_{B_r} u,\quad  m_r=\inf_{B_r} u.
\end{equation*}
Set $r_0=1/2^k$. Note that
\begin{equation*}
  |\{x\in B_{r_0}: u(x)\geq \frac{M_{r_0}+m_{r_0}}{2}\}|\geq \frac{|B_{r_0}|}{2}~~\mbox{ or }~
  |\{x\in B_{r_0}: u(x)\leq \frac{M_{r_0}+m_{r_0}}{2}\}|\geq \frac{|B_{r_0}|}{2}.
\end{equation*}
Without loss of generality, we assume that the former holds. Let
\begin{equation*}
y=\frac{x}{r_0},\quad \tilde{u}(y)=\frac{u(x)-m_{r_0}}{M_{r_0}-m_{r_0}},\quad
\tilde{f}(y)=\frac{r_0^2f(x)}{M_{r_0}-m_{r_0}}.
\end{equation*}
Without loss of generality, we assume that
\begin{equation*}
  \underset{B_{r_0}}{\mathrm{osc}}~u=M_{r_0}-m_{r_0}\geq r_0^2.
\end{equation*}
Note that $M_{r_0}-m_{r_0}\leq 2$. Then
\begin{equation*}
  \tilde{u}\in S^{*}_{r_0^2\rho/2}(\lambda,\Lambda,b_0,\tilde f)~~\mbox{ in}~B_1
\end{equation*}
and
\begin{equation*}
0\leq \tilde u\leq 1~~\mbox{ in}~B_1,\quad \|\tilde f\|_{L^{\infty}(B_1)}\leq \varepsilon_1.
\end{equation*}
Moreover,
\begin{equation*}
  |\{y\in B_{1}: \tilde u(y)\geq \frac{1}{2}\}|\geq \frac{|B_{1}|}{2}.
\end{equation*}
Hence, by combining with \Cref{th10.3} and noting
\begin{equation*}
\frac{r_0^2\rho}{2}\geq 8,\quad \frac{r_0^2\rho}{2\rho_3}\geq 1,
\end{equation*}
we have
\begin{equation*}
c_0\leq \|\tilde u\|_{L^{p_0}(B_{1/2})}= \|\tilde u_{r_0^2\rho/2}\|_{L^{p_0}(B_{1/2})}\leq
C\left(\inf_{B_{1/2}} \tilde u+\|\tilde f\|_{L^{\infty}(B_1)}\right)\leq C\inf_{B_{1/2}} \tilde
u+C\varepsilon_1,
\end{equation*}
where $c_0>0,C$ are universal constants.

Choose $\varepsilon_1$ small enough such that
\begin{equation*}
  C\varepsilon_1\leq \frac{c_0}{2}.
\end{equation*}
Then for some universal constant $\mu>0$,
\begin{equation*}
\inf_{B_{1/2}} \tilde u\geq \mu.
\end{equation*}
By transforming back to $u$,
\begin{equation*}
  m_{r_0/2}-m_{r_0}\geq \mu\left(M_{r_0}-m_{r_0}\right).
\end{equation*}
Hence,
\begin{equation*}
\underset{B_{1/2^{k+1}}}{\mathrm{osc}}~u=M_{r_0/2}-m_{r_0/2}
\leq M_{r_0}-m_{r_0}-\mu\left(M_{r_0}-m_{r_0}\right)=(1-\mu)\left(M_{r_0}-m_{r_0}\right)
\leq 2(1-\mu)^{k+1}.
\end{equation*}
By induction, the proof is completed.~\qed~\\

\begin{remark}\label{re10.2}
Since $\rho$ is finite, we can make the scaling argument only finite times. Hence, we can't obtain the real
H\"{o}lder regularity. However, it can provide necessary compactness when we use the compactness method
to prove higher regularity (see \Cref{le3.1} and \Cref{le11.1}).
\end{remark}
~\\

By applying above theorem to each $x_0\in B_{1/2}$ in $B(x_0, 1/2)$, we have
\begin{corollary}\label{co10.2}
Let $\rho\geq 2\rho_3$,
\begin{equation*}
u\in S^{*}_{\rho}(\lambda,\Lambda,b_0,f)~~\mbox{ in}~B_1
\end{equation*}
and
\begin{equation*}
\|u\|_{L^{\infty}(B_1)}\leq 1,\quad \|f\|_{L^{\infty}(B_1)}\leq \varepsilon_1,
\end{equation*}
where $0<\varepsilon_1<1$ is universal. Then
\begin{equation*}
\underset{B_r(x_0)}{\mathrm{osc}}~u\leq \tilde Cr^{\alpha_0},~\forall ~x_0\in B_{1/2},~
\forall~\sqrt{\frac{2\rho_3}{\rho}}\leq r\leq \frac{1}{2},
\end{equation*}
where $0<\alpha_0<1$ and $\tilde C$ are universal.

Therefore,
\begin{equation*}
  |u(x_1)-u(x_2)|\leq \tilde C|x_1-x_2|^{\alpha_0},~\forall ~x_1,x_2\in B_{1/2},~
  \sqrt{\frac{2\rho_3}{\rho}}\leq |x_1-x_2|\leq \frac{1}{2}.
\end{equation*}
\end{corollary}

\section{Interior \texorpdfstring{$C^{1,\alpha}$}{C1,a} regularity}\label{S3}
In this section, we prove the interior pointwise $C^{1,\alpha}$ regularity by the classical technique of
perturbation. First, we prove the key step by the compactness method.
\begin{lemma}\label{le3.1}
Let $0<\alpha<1$ and $u\in C(\bar{B}_1)$ be a viscosity solution of
\begin{equation*}
  F(D^2u,Du,u,x)= A^{ij}(Du,u,x)u_{ij}+B(Du,u,x)=0~~\mbox{ in}~B_1,
\end{equation*}
where $F$ is $\rho$-uniformly elliptic and $A$ is continuous with modulus $\omega_A$. Suppose that
\cref{e1.4} and \cref{e1.5} hold. Let $r\leq r_0$ and assume that for some $P_0\in \mathcal{P}_1$,
\begin{equation*}
\|u-P_0\|_{L^{\infty}(B_{r})}\leq r^{1+\alpha},\quad \|P_0\|\leq \bar{C}r_0^{1+\alpha},
\end{equation*}
where $\bar{C}$ depends only on $n,\lambda,\Lambda,\rho,b_0$ and $\alpha$ and $0<r_0,\delta_0<1$
(small) depends also on $\omega_A$.

Then there exists $P\in \mathcal{P}_1$ such that
\begin{equation*}
  \|u-P\|_{L^{\infty}(B_{\eta r})}\leq  (\eta r) ^{1+\alpha},\quad \|P-P_0\|_r\leq \bar{C}(\eta r)^{1+\alpha},
\end{equation*}
where $0<\eta<1/2$ depends only on $n,\lambda,\Lambda,\rho,b_0,c_0$ and $\alpha$.
\end{lemma}
\proof We prove the lemma by contradiction. Suppose that the conclusion is false. Then there exist sequences
of $A_m,B_m,u_m,P_m,r_m$ such that $r_m\leq 1/m$,
\begin{equation*}
F_m(D^2u_m,Du_m,u_m,x)=A^{ij}_m(Du_m,u_m,x)u_{m,ij}+B_m(Du_m,u_m,x)
=0~~\mbox{ in}~B_{r_m}
\end{equation*}
and
\begin{equation*}
\|u_m-P_m\|_{L^{\infty}(B_{r_m})}\leq r_m^{1+\alpha},\quad \|P_m\|\leq \frac{\bar{C}}{m},
\end{equation*}
where $F_m$ are $\rho$-uniformly elliptic and $A_m$ are continuous (with the same modulus $\omega_A$).
In addition, \cref{e1.4} and \cref{e1.5} hold for $B_m$ (with $K_B$ and $b_0$). Moreover, for any $P\in
\mathcal{P}_1$ with
\begin{equation*}
\|P-P_m\|\leq \bar{C}(\eta r_m)^{1+\alpha},
\end{equation*}
we have
\begin{equation}\label{e2.2}
\|u_m-P\|_{L^{\infty}(B_{\eta r_m})}> (\eta r_m)^{1+\alpha},
\end{equation}
where $\bar{C}$ and $0<\eta<1/2$ are to be specified later.

Let
\begin{equation*}
\tilde{x}=\frac{x}{r_m},\quad \tilde{u}_m(\tilde{x})=\frac{u_m(x)-P_m(x)}{r_m^{1+\alpha}}.
\end{equation*}
Then $\tilde{u}_m$ are viscosity solutions of
\begin{equation}\label{e2.4}
    \tilde F_m(D^2\tilde u_m,D\tilde u_m,\tilde u_m,\tilde{x})=0~~\mbox{ in}~B_1,
\end{equation}
where
\begin{equation*}
  \begin{aligned}
\tilde{F}_m(M,p,s,\tilde{x})=&r_m^{1-\alpha}
F_m(r_m^{\alpha-1}M,r_m^{\alpha}p+DP_m(x),r_m^{1+\alpha}s+P_m(x),x),~\\
=&A^{ij}_m(r_m^{\alpha}p+DP_m(x),r_m^{1+\alpha}s+P_m(x),x)M_{ij}\\
&+r_m^{1-\alpha}B_m(r_m^{\alpha}p+DP_m(x),r_m^{1+\alpha}s+P_m(x),x).
  \end{aligned}
\end{equation*}
Indeed, since $F_m$ are $\rho$-uniformly elliptic and $\|P_m\|\to 0$, $\tilde{F}_m$ are $\frac{1}{2}
r_m^{-\alpha}\rho$-uniformly elliptic for $m$ large enough. Note that
\begin{equation*}
\|\tilde{u}_m\|_{L^{\infty}(B_1)}\leq 1\leq \frac{1}{4} r_m^{-\alpha}\rho.
\end{equation*}
Then it can be verified directly that $\tilde{u}_m$ are viscosity solutions of \cref{e2.4}.

Next, by \Cref{pr10.1}, for $m$ large enough,
\begin{equation*}
  \tilde u_m\in S^*_{\frac{1}{2}r_m^{-\alpha}\rho}(\lambda,\Lambda,r_mb_0,\bar{f}_m),
\end{equation*}
where
\begin{equation*}
\bar{f}_m(\tilde{x})=r_m^{1-\alpha}B_m(DP_m(x),r_m^{1+\alpha}\tilde u(\tilde{x})+ P_m(x),x).
\end{equation*}
Hence, for any $\varepsilon>0$, we can take $m_0$ large enough such that for any $m\geq m_0$,
\begin{equation*}
  \|\bar{f}_m\|_{L^{\infty}(B_1)}\leq \varepsilon_1,\quad
  \tilde C\left(\frac{4\rho_3}{r_m^{-\alpha}\rho}\right)^{\alpha_0/2}\leq \frac{\varepsilon}{2},
\end{equation*}
where $\alpha_0,\tilde C,\rho_3$ and $\varepsilon_1$ are as in \Cref{co10.2}. Set
\begin{equation*}
\delta=\left(\frac{4\rho_3}{r_m^{-\alpha}\rho}\right)^{1/2}.
\end{equation*}
By \Cref{co10.2}, for any $x_1,x_2\in B_{1/2}$ with $|x_1-x_2|\leq \delta$, by choosing $x_3\in B_{1/2}$
with $|x_1-x_3|=|x_2-x_3|=\delta$, we have for any $m\geq m_0$,
\begin{equation*}
|\tilde{u}_m(x_1)-\tilde{u}_m(x_2)|\leq |\tilde{u}_m(x_1)-\tilde{u}_m(x_3)|
+|\tilde{u}_m(x_2)-\tilde{u}_m(x_3)|\leq 2\tilde C\delta^{\alpha_0}\leq \varepsilon.
\end{equation*}
Thus, $\tilde{u}_m$ are equicontinuous. By Arzel\`{a}-Ascoli theorem, there exists $\tilde{u}\in
C(\bar{B}_{1/2})$ such that
$\tilde{u}_m\to \tilde{u}$ in $L^{\infty}(B_{1/2})$ (up to a subsequence and similarly hereinafter).

Since
\begin{equation*}
\lambda I\leq A_m(0,0,0)\leq \Lambda I,~\forall ~m\geq 1,
\end{equation*}
there exists a constant symmetric matrix $A=(A^{ij})_{n\times n}$ such that
\begin{equation*}
A_m(0,0,0)\to A.
\end{equation*}

Now, we show that $\tilde{u}$ is a viscosity solution of
\begin{equation}\label{e2.5}
A^{ij}\tilde{u}_{ij}=0~~\mbox{ in}~B_{1/2}.
\end{equation}
Given $\tilde x_0\in B_{1/2}$ and $\varphi\in C^2$ touching $\tilde{u}$ strictly by above at $\tilde x_0$.
Then there exist a sequence of $\tilde x_m\to \tilde x_0$ such that $\varphi+c_m$ touch $\tilde{u}_m$ by
above at $\tilde x_m$ and $c_m\to 0$. By the definition of viscosity solution, for $m$ large enough (e.g.
$r_m^{-\alpha}\rho>2\|\varphi\|_{C^{2}(\bar{B}_{1/2})}$),
\begin{equation*}
  \tilde{F}_m(D^2\varphi(\tilde x_m),D\varphi(\tilde x_m),\varphi(\tilde x_m)+c_m, \tilde x_m)\geq 0.
\end{equation*}
Since
\begin{equation*}
  \begin{aligned}
\tilde{F}_m&(D^2\varphi(\tilde x_m),D\varphi(\tilde x_m),\varphi(\tilde x_m)+c_m, \tilde x_m)\\
=&A^{ij}_m(r_m^{\alpha}D\varphi(\tilde x_m)+DP_m(x),r_m^{1+\alpha}\varphi(\tilde
x_m)+P_m(x),x)M_{ij}\\
&+r_m^{1-\alpha}B_m(r_m^{\alpha}D\varphi(\tilde x_m)+DP_m(x),r_m^{1+\alpha}\varphi(\tilde
x_m)+P_m(x),x),\\
\leq &A^{ij}_m(0,0,0)M_{ij}+\omega_A(\frac{C}{m})\\
&+r_m^{1-\alpha}B_m(r_m^{\alpha}D\varphi(\tilde x_m)+DP_m(x),r_m^{1+\alpha}\varphi(\tilde
x_m)+P_m(x),x),
  \end{aligned}
\end{equation*}
by letting $m\to \infty$, we have
\begin{equation*}
A^{ij}\varphi_{ij}(\tilde x_0)\geq 0.
\end{equation*}
Hence, $\tilde{u}$ is a subsolution of \cref{e2.5}. Similarly, we can prove that it is a viscosity supersolution
as well. That is, $\tilde{u}$ is a viscosity solution.

Since \cref{e2.5} is a linear equation with constant coefficients, $\tilde{u}\in C^{\infty}(B_{1/2})$. Then
there exists $\tilde{P}\in \mathcal{P}_1$ such that for any $0<\eta<1/4$,
\begin{equation*}
  \|\tilde{u}-\tilde P\|_{L^{\infty}(B_{\eta})}\leq C_1\eta^{2}\|\tilde{u}\|_{L^{\infty}(B_{1/2})}\leq
  C_1\eta^{2},
\end{equation*}
and
\begin{equation*}
  \|\tilde P\|\leq C_2\|\tilde{u}\|_{L^{\infty}(B_{1/2})}\leq C_2,
\end{equation*}
where $C_1$ and $C_2$ are universal. By taking $\eta$ small and $\bar{C}$ large such that
\begin{equation*}
C\eta^{1-\alpha}\leq \frac{1}{2},\quad C_2\leq \bar{C}\eta^{1+\alpha}.
\end{equation*}
Then
\begin{equation}\label{e2.6}
  \|\tilde{u}-\tilde P\|_{L^{\infty}(B_{\eta})}\leq \frac{1}{2}\eta^{1+\alpha},\quad
   \|\tilde P\|\leq \bar{C}\eta^{1+\alpha}.
\end{equation}

Let
\begin{equation*}
Q_m(x)=P_m(x)+r_m^{1+\alpha}\tilde P(\tilde x).
\end{equation*}
Then
\begin{equation*}
\|Q_m-P_m\|\leq \bar{C}(\eta r_m)^{1+\alpha}.
\end{equation*}
Hence, \cref{e2.2} holds for $Q_m$. That is,
\begin{equation*}
\|u_m-Q_m\|_{L^{\infty}(B_{\eta}r_m)}> (\eta r_m)^{1+\alpha}.
\end{equation*}
Equivalently,
\begin{equation*}
\|\tilde u_m-\tilde P\|_{L^{\infty}(B_{\eta})}> \eta^{1+\alpha}.
\end{equation*}
Let $m\to \infty$, we have
\begin{equation*}
\|\tilde u-\tilde P\|_{L^{\infty}(B_{\eta})}\geq \eta^{1+\alpha},
\end{equation*}
which contradicts with \cref{e2.6}.~\qed~\\

Now, we give the~\\
\noindent\textbf{Proof of \Cref{th1.1}.} To prove \Cref{th1.1}, we only to prove the following. There exist a
sequence of $P_m\in \mathcal{P}_1$ ($m\geq -1$) such that for all $m\geq 0$,

\begin{equation}\label{e3.2-1}
\|u-P_m\|_{L^{\infty }(B _{\eta^{m}r_0})}\leq (\eta ^{m}r_0)^{1+\alpha},
\end{equation}
and
\begin{equation}\label{e3.3-1}
\|P_m-P_{m-1}\|_{\eta^{m}r_0}\leq \bar{C}(\eta ^{m}r_0)^{1+\alpha},
\end{equation}
where $\eta,r_0$ and $\bar{C}$ are as in \Cref{le3.1}.

We prove above by induction. For $m=0$, by setting $P_0\equiv P_{-1}\equiv 0$, \cref{e3.2-1} and
\cref{e3.3-1} hold clearly. Suppose that the conclusion holds for $m\leq m_0$. By \cref{e3.3-1},
\begin{equation*}
\|P_{m_0}\|\leq \sum_{i=1}^{m_0}\|P_i-P_{i-1}\|\leq
\bar{C}r_0^{1+\alpha}\frac{\eta^{1+\alpha}}{1-\eta^{1+\alpha}}\leq \bar{C}r_0^{1+\alpha}.
\end{equation*}
By \Cref{le3.1}, the conclusion holds for $m=m_0+1$. By induction, the proof of \Cref{th1.1} is completed.
~\qed~\\

\section{Interior \texorpdfstring{$C^{2,\alpha}$}{C2,a} regularity}\label{S11}
In this section, we prove the interior pointwise $C^{2,\alpha}$ regularity. As before, we first prove the key step by the compactness method.
\begin{lemma}\label{le11.1}
Let $0<\alpha<1$ and $u\in C(\bar{B}_1)$ be a viscosity solution of
\begin{equation*}
  F(D^2u,Du,u,x)=f~~\mbox{ in}~B_1,
\end{equation*}
where $F$ is $\rho$-uniformly elliptic, $F(0,0,0,x)\equiv 0$ and $D_MF$ is continuous with modulus
$\omega_F$. Suppose that \cref{e.st.1} holds. Let $r\leq r_0$ and assume that for some $P_0\in
\mathcal{P}_2$,
\begin{equation*}
\|F(M,p,s,x)-F(M,p,s,0)\|_{L^{\infty}(B_{r})}\leq \delta_0r^{\alpha},~\forall ~|M|,|p|,|s|\leq \rho,
\end{equation*}
\begin{equation*}
\|u-P_0\|_{L^{\infty}(B_{r})}\leq r^{2+\alpha},\quad
\|f-f(0)\|_{L^{\infty}(B_{r})}\leq \delta_0r^{\alpha}
\end{equation*}
and
\begin{equation*}
\|P_0\|\leq \bar{C}r_0^{2+\alpha},\quad   F(D^2P_0,DP_0(0),P_0(0),0)=f(0),\quad
\end{equation*}
where $\bar{C}$ depends only on $n,\lambda,\Lambda,\rho,b_0,c_0$ and $\alpha$ and $0<r_0,\delta_0<1$
(small) depends also on $\omega_F$.

Then there exists $P\in \mathcal{P}_2$ such that
\begin{equation*}
  \|u-P\|_{L^{\infty}(B_{\eta r})}\leq  (\eta r) ^{2+\alpha}
\end{equation*}
and
\begin{equation*}
 \|P-P_0\|_r\leq \bar{C}(\eta r)^{2+\alpha},\quad   F(D^2P,DP(0),P(0),0)=f(0),\quad
\end{equation*}
where $0<\eta<1/2$ depends only on $n,\lambda,\Lambda,\rho,b_0,c_0$ and $\alpha$.
\end{lemma}
\proof We prove the lemma by contradiction. Suppose that the conclusion is false. Then there exist sequences
of $F_m,u_m,f_m,P_m,r_m$ such that $r_m\leq 1/m$,
\begin{equation*}
  F_m(D^2u_m,Du_m,u_m,x)=f_m~~\mbox{ in}~B_{r_m},
\end{equation*}
\begin{equation*}
\|F_m(M,p,s,x)-F_m(M,p,s,0)\|_{L^{\infty}(B_{r_m})}\leq \frac{r_m^{\alpha}}{m},~\forall
~|M|,|p|,|s|\leq \rho,
\end{equation*}
\begin{equation*}
\|u_m-P_m\|_{L^{\infty}(B_{r_m})}\leq r_m^{2+\alpha},\quad \|f_m-f_m(0)\|_{L^{\infty}(B_{r_m})}\leq
\frac{r_m^{\alpha}}{m}
\end{equation*}
and
\begin{equation*}
 F_m(D^2P_m,DP_m(0),P_m(0),0)=f_m(0),\quad \|P_m\|\leq \frac{\bar{C}}{m},\quad
\end{equation*}
where $F_m$ are $\rho$-uniformly elliptic, $F_m(0,0,0,x)\equiv 0$ and $D_MF_m$ are continuous (with the
same modulus $\omega_F$). In addition, \cref{e.st.1} holds for $F_m$ with $b_0$ and $c_0$. Moreover, for
any $P\in \mathcal{P}_2$ with
\begin{equation*}
\|P-P_m\|\leq \bar{C}(\eta r_m)^{2+\alpha},\quad F_m(D^2P,DP(0),P(0),0)=f_m(0),\quad
\end{equation*}
we have
\begin{equation}\label{e11.1}
\|u_m-P\|_{L^{\infty}(B_{\eta r_m})}> (\eta r_m)^{2+\alpha},
\end{equation}
where $\bar{C}$ and $0<\eta<1/2$ are to be specified later.

Let
\begin{equation}\label{e3.5}
\tilde{x}=\frac{x}{r_m},\quad \tilde{u}_m(\tilde{x})=\frac{u_m(x)-P_m(x)}{r_m^{2+\alpha}}.
\end{equation}
Then $\tilde{u}_m$ are viscosity solutions of
\begin{equation}\label{e3.1}
    \tilde F_m(D^2\tilde u_m,D\tilde u_m,\tilde u_m,\tilde{x})=\tilde f_m~~\mbox{ in}~B_1,
\end{equation}
where
\begin{equation*}
  \begin{aligned}
&\tilde{F}_m(M,p,s,\tilde{x})=r_m^{-\alpha}
F_m(r_m^{\alpha}M+D^2P_m,r_m^{1+\alpha}p+DP_m(x),r_m^{2+\alpha}s+P_m(x),x),~\\
&\tilde{f}_m(\tilde{x})=r_m^{-\alpha}f_m(x).\\
  \end{aligned}
\end{equation*}
Indeed, since $F_m$ are $\rho$-uniformly elliptic and $\|P_m\|\to 0$, $\tilde{F}_m$ are $\frac{1}{2}
r_m^{-\alpha}\rho$-uniformly elliptic for $m$ large enough. Note that
\begin{equation*}
\|\tilde{u}_m\|_{L^{\infty}(B_1)}\leq 1\leq \frac{1}{4} r_m^{-\alpha}\rho.
\end{equation*}
Then it can be verified directly that $\tilde{u}_m$ are viscosity solutions of \cref{e3.1}.

Next, by \Cref{pr10.1}, for $m$ large enough,
\begin{equation*}
  \tilde u_m\in S^*_{\frac{1}{2}r_m^{-\alpha}\rho}(\lambda,\Lambda,r_mb_0,\bar{f}_m),
\end{equation*}
where
\begin{equation*}
\bar{f}_m(\tilde{x})=|r_m^{-\alpha}
\left(f_m(x)-F_m(D^2P_m,DP_m(x),P_m(x),x)\right)|+r_m^2c_0|u_m(x)|.
\end{equation*}
With the aid of $F_m(D^2P_m,DP_m(0),P_m(0),0)=f_m(0)$,
\begin{equation*}
  \begin{aligned}
|f_m&-F_m(D^2P_m,DP_m,P_m,x)|\\
=&\left|f_m-f_m(0)-\left(F_m(D^2P_m,DP_m,P_m,x)-F_m(D^2P_m,DP_m,P_m,0)\right)\right.\\
&\left.-\left(F_m(D^2P_m,DP_m,P_m,0)-F_m(D^2P_m,DP_m(0),P_m(0),0)\right)\right|\\
\leq& \frac{r_m^{\alpha}}{m}+\frac{r_m^{\alpha}}{m}+b_0Cr_m+c_0Cr_m.
  \end{aligned}
\end{equation*}
Hence, for any $\varepsilon>0$, we can take $m_0$ large enough such that for any $m\geq m_0$,
\begin{equation*}
  \|\bar{f}_m\|_{L^{\infty}(B_1)}\leq \varepsilon_1,\quad
  \tilde C\left(\frac{4\rho_3}{r_m^{-\alpha}\rho}\right)^{\alpha_0/2}\leq \frac{\varepsilon}{2},
\end{equation*}
where $\alpha_0,\tilde C,\rho_3$ and $\varepsilon_1$ are as in \Cref{co10.2}. Set
\begin{equation*}
\delta=\left(\frac{4\rho_3}{r_m^{-\alpha}\rho}\right)^{1/2}.
\end{equation*}
By \Cref{co10.2}, for any $x_1,x_2\in B_{1/2}$ with $|x_1-x_2|\leq \delta$, by choosing $x_3\in B_{1/2}$
with $|x_1-x_3|=|x_2-x_3|=\delta$, we have for any $m\geq m_0$,
\begin{equation*}
|\tilde{u}_m(x_1)-\tilde{u}_m(x_2)|\leq |\tilde{u}_m(x_1)-\tilde{u}_m(x_3)|
+|\tilde{u}_m(x_2)-\tilde{u}_m(x_3)|\leq 2\tilde C\delta^{\alpha_0}\leq \varepsilon.
\end{equation*}
Thus, $\tilde{u}_m$ are equicontinuous. By Arzel\`{a}-Ascoli theorem, there exists $\tilde{u}\in
C(\bar{B}_{1/2})$ such that
$\tilde{u}_m\to \tilde{u}$ in $L^{\infty}(B_{1/2})$ (up to a subsequence and similarly hereinafter).

Since $F_m$ are $\rho$-uniformly elliptic,
\begin{equation*}
\lambda I\leq D_MF_m(0,0,0,0)\leq \Lambda I,~\forall ~m\geq 1.
\end{equation*}
Hence, there exists a constant symmetric matrix $A=(A^{ij})_{n\times n}$ such that
\begin{equation*}
D_MF_m(0,0,0,0)\to A.
\end{equation*}

Now, we show that $\tilde{u}$ is a viscosity solution of
\begin{equation}\label{e11.3}
A^{ij}\tilde{u}_{ij}=0~~\mbox{ in}~B_{1/2}.
\end{equation}
Given $\tilde x_0\in B_{1/2}$ and $\varphi\in C^2$ touching $\tilde{u}$ strictly by above at $\tilde x_0$.
Then there exist a sequence of $\tilde x_m\to \tilde x_0$ such that $\varphi+c_m$ touch $\tilde{u}_m$ by
above at $\tilde x_m$ and $c_m\to 0$. By the definition of viscosity solution, for $m$ large enough (e.g.
$r_m^{-\alpha}\rho>2\|\varphi\|_{C^{2}(\bar{B}_{1/2})}$),
\begin{equation*}
  \tilde{F}_m(D^2\varphi(\tilde x_m),D\varphi(\tilde x_m),\varphi(\tilde x_m)+c_m, \tilde x_m)\geq
  \tilde{f}_m(\tilde x_m).
\end{equation*}
We compute
\begin{equation*}
  \begin{aligned}
\tilde{F}_m&(D^2\varphi(\tilde x_m),D\varphi(\tilde x_m),\varphi(\tilde x_m)+c_m, \tilde x_m)\\
=&\frac{1}{r_m^{\alpha}}F_m(r_m^{\alpha}D^2\varphi+D^2P_m,
r_m^{1+\alpha}D\varphi+DP_m,r_m^{2+\alpha}(\varphi+c_m)+P_m,r_m\tilde x_m)\\
&-\frac{1}{r_m^{\alpha}}F_m(r_m^{\alpha}D^2\varphi+D^2P_m,
r_m^{1+\alpha}D\varphi+DP_m,r_m^{2+\alpha}(\varphi+c_m)+P_m,0)\\
&+\frac{1}{r_m^{\alpha}}F_m(r_m^{\alpha}D^2\varphi+D^2P_m,
r_m^{1+\alpha}D\varphi+DP_m,r_m^{2+\alpha}(\varphi+c_m)+P_m,0)\\
&-\frac{1}{r_m^{\alpha}}F_m(r_m^{\alpha}D^2\varphi+D^2P_m,DP_m(0),P_m(0),0)\\
&+\frac{1}{r_m^{\alpha}}F_m(r_m^{\alpha}D^2\varphi+D^2P_m,DP_m(0),P_m(0),0)\\
&-\frac{1}{r_m^{\alpha}}F_m(D^2P_m,DP_m(0),P_m(0),0)+\frac{1}{r_m^{\alpha}}f_m(0)\\
\leq&\frac{1}{m}+Cb_0r_m^{1-\alpha}+Cc_0r_m^{1-\alpha}\\
&+F_{m,M_{ij}}(\theta r_m^{\alpha}D^2\varphi(\tilde
x_m)+D^2P_m,DP_m(0),P_m(0),0)\varphi_{ij}(\tilde x_m)+\frac{1}{r_m^{\alpha}}f_m(0),
  \end{aligned}
\end{equation*}
where $0<\theta<1$. Note that in above inequality, the variable of $\varphi$ is $\tilde{x}_m$ and the
variable of $P_m$ is $x_m:=r_m\tilde{x}_m$. Hence,
\begin{equation*}
  \begin{aligned}
0\leq&\tilde{F}_m(D^2\varphi(\tilde x_m),D\varphi(\tilde x_m),\varphi(\tilde x_m)+c_m, \tilde x_m)
-\tilde{f}_m(\tilde x_m)\\
\leq & \frac{1}{m}+Cb_0r_m^{1-\alpha}+Cc_0r_m^{1-\alpha}\\
&+F_{m,M_{ij}}(0,0,0,0)\varphi_{ij}(\tilde x_m)+\frac{1}{r_m^{\alpha}}|f_m(r_m\tilde
x_m)-f_m(0)|+\omega_F(\frac{C}{m}).
  \end{aligned}
\end{equation*}
Let $m\to \infty$, we have
\begin{equation}\label{e3.4}
A^{ij}\varphi_{ij}(\tilde x_0)\geq 0.
\end{equation}
Hence, $\tilde{u}$ is a subsolution of \cref{e11.3}. Similarly, we can prove that it is a viscosity supersolution
as well. That is, $\tilde{u}$ is a viscosity solution.
%
bounded and equicontinuous in any compact set of $\mathcal{S}^n$. Then there exists a linear operator
$\tilde F$ with ellipticity constants $\lambda$ and $\Lambda$ such that $\tilde{F}_m\to \tilde{F}$ and
$D\tilde{F}_m\to D\tilde{F}$ uniformly in any compact set of $\mathcal{S}^n$ .
%

Since \cref{e11.3} is a linear equation with constant coefficients, $\tilde{u}\in C^{\infty}(B_{1/2})$. Then
there exists $\tilde{P}\in \mathcal{P}_2$ such that for any $0<\eta<1/4$,
\begin{equation*}
  \|\tilde{u}-\tilde P\|_{L^{\infty}(B_{\eta})}\leq C_1\eta^{3}\|\tilde{u}\|_{L^{\infty}(B_{1/2})}\leq
  C_1\eta^{3},
\end{equation*}
\begin{equation*}
A^{ij}\tilde P_{ij}=0
\end{equation*}
and
\begin{equation*}
  \|\tilde P\|\leq C_2\|\tilde{u}\|_{L^{\infty}(B_{1/2})}\leq C_2,
\end{equation*}
where $C_1$ and $C_2$ are universal. By taking $\eta$ small and $\bar{C}$ large such that
\begin{equation*}
C\eta^{1-\alpha}\leq \frac{1}{2},\quad C_2\leq (\bar{C}-1)\eta^{2+\alpha}.
\end{equation*}
Then
\begin{equation}\label{e11.2}
  \|\tilde{u}-\tilde P\|_{L^{\infty}(B_{\eta})}\leq \frac{1}{2}\eta^{2+\alpha},\quad
   \|\tilde P\|\leq (\bar{C}-1)\eta^{2+\alpha}.
\end{equation}
By a similar argument to prove \cref{e3.4},
\begin{equation*}
  \tilde{F}_m(D^2\tilde P,D\tilde P(0),\tilde P(0),0)\to A^{ij}\tilde P_{ij}=0.
\end{equation*}
Note that $\tilde f_m(0)\to 0$ as well. Hence, there exist a sequence of constants $t_m\to 0$ such that
\begin{equation*}
\tilde{F}_m(D^2\tilde P+t_mI,D\tilde P(0),\tilde P(0),0)=\tilde{f}_m(0).
\end{equation*}

Let
\begin{equation*}
Q_m(x)=P_m(x)+r_m^{2+\alpha}\left(\tilde P(\tilde x)+\frac{t_m}{2}|\tilde x|^2\right).
\end{equation*}
Then with the aid of \cref{e11.2},
\begin{equation*}
F_m(D^2Q_m,DQ_m(0),Q_m(0),0)=f_m(0),\quad \|Q_m-P_m\|\leq \bar{C}(\eta r_m)^{2+\alpha}.
\end{equation*}
Hence, \cref{e11.1} holds for $Q_m$. That is,
\begin{equation*}
\|u_m-Q_m\|_{L^{\infty}(B_{\eta}r_m)}> (\eta r_m)^{2+\alpha}.
\end{equation*}
Equivalently,
\begin{equation*}
\|\tilde u_m-\tilde P-\frac{t_m}{2}|\tilde x|^2\|_{L^{\infty}(B_{\eta})}> \eta^{2+\alpha}.
\end{equation*}
Let $m\to \infty$, we have
\begin{equation*}
\|\tilde u-\tilde P\|_{L^{\infty}(B_{\eta})}\geq \eta^{2+\alpha},
\end{equation*}
which contradicts with \cref{e11.2}.~\qed~\\

Now, we give the~\\
\noindent\textbf{Proof of \Cref{th1.2}.} To prove \Cref{th1.2}, we only to prove the following. There exist a
sequence of $P_m\in \mathcal{P}_2$ ($m\geq -1$) such that for all $m\geq 0$,

\begin{equation}\label{e3.2}
\|u-P_m\|_{L^{\infty }(B _{\eta^{m}r_0})}\leq (\eta ^{m}r_0)^{2+\alpha},
\end{equation}
and
\begin{equation}\label{e3.3}
F(D^2P_m,DP_m(0),P_m(0),0)=f(0),\quad \|P_m-P_{m-1}\|_{\eta^{m}r_0}\leq \bar{C}(\eta
^{m}r_0)^{2+\alpha},
\end{equation}
where $\eta,r_0$ and $\bar{C}$ are as in \Cref{le11.1}.

We prove above by induction. Set $P_{-1}\equiv 0$. For $m=0$, since $F(0,0,0,0)=0$, there exists $t\in
\mathbb{R}$ such that
\begin{equation*}
F(tI,0,0,0)=f(0),\quad |t|\leq \frac{|f(0)|}{n\lambda}\leq \frac{\delta}{n\lambda}.
\end{equation*}
Then by choosing $P_0=t|x|^2/2$ and $\delta=(1+1/n\lambda)r_0^{2+\alpha}$,
\begin{equation*}
\|u-P_0\|_{L^{\infty }(B_{r_0})}\leq \|u\|_{L^{\infty }(B_{r_0})}+\|P_0\|_{L^{\infty }(B_{r_0})}\leq
\delta +\frac{\delta}{2n\lambda}\leq r_0^{2+\alpha}.
\end{equation*}
Hence, \cref{e3.2} and \cref{e3.3} hold for $m=0$. Suppose that the conclusion holds for $m\leq m_0$. By
\cref{e3.3},
\begin{equation*}
\|P_{m_0}\|\leq \sum_{i=1}^{m_0}\|P_i-P_{i-1}\|+\|P_0\|\leq
\bar{C}r_0^{2+\alpha}\frac{\eta^{2+\alpha}}{1-\eta^{2+\alpha}}+r_0^{2+\alpha}\leq
\bar{C}r_0^{2+\alpha}.
\end{equation*}
By \Cref{le11.1}, the conclusion holds for $m=m_0+1$. By induction, the proof of \Cref{th1.2} is completed.
~\qed~\\

\begin{remark}\label{re1.4}
Note that we can not prove the pointwise $C^{1,\alpha}$ regularity for a general operator $F$. The reason is the following. For the $C^{1,\alpha}$ regularity, we will consider for
some $P_m\in \mathcal{P}_1$
\begin{equation*}
  \tilde{x}=\frac{x}{r_m},\quad \tilde u_m(\tilde{x})=\frac{u_m(x)-P_n(x)}{r_m^{1+\alpha}}
\end{equation*}
instead of \cref{e3.5} in the proof of \Cref{le11.1}. Then $\tilde{u}_m$ are solutions of
\begin{equation*}
  \begin{aligned}
\tilde{F}_m(D^2\tilde{u}_m,D\tilde{u}_m,\tilde{u}_m,\tilde{x})
:=r_m^{1-\alpha}F_m(r_m^{\alpha-1}D^2u_m,r_m^{\alpha}Du_m,r_m^{\alpha+1}u_m,x)
=r^{1-\alpha}f_m.\\
  \end{aligned}
\end{equation*}
Note that $\tilde{F}_m$ are only $r_m^{1-\alpha}\rho$-uniformly elliptic and $r_m^{1-\alpha}\rho\to 0$.
Hence, we can not proceed the scaling argument.
\end{remark}
~\\

\section{Interior \texorpdfstring{$C^{k,\alpha}$}{Ck,a} regularity}\label{S4}
In this section, we prove the interior pointwise $C^{k,\alpha}$ regularity \Cref{th1.3}. Since $k\geq 2$, the
assumption of \Cref{th1.3} is stronger than that of \Cref{th1.2}. Hence, $u\in C^{2,\alpha}(0)$. Then we
only need to prove:\\

\textbf{Claim I: If $u\in C^{k-1,\alpha}(0)$ under the assumptions of \Cref{th1.3},
 then $u\in C^{k,\alpha}(0)$}.\\

Indeed, we can prove \Cref{th1.3} by induction if \textbf{Claim I} has been proved. For $k=3$, by
\Cref{th1.2}, $u\in C^{2,\alpha}(0)$. Then \textbf{Claim I} implies $u\in C^{3,\alpha}(0)$. Hence,
\Cref{th1.3} holds for $k=3$. We assume that \Cref{th1.3} holds for $k=k_0$ and we only need to prove it
for $k_0+1$. Since \Cref{th1.3} holds for $k_0$ and the assumptions of \Cref{th1.3} with $k_0+1$ is
stronger than that with $k_0$, $u\in C^{k_0,\alpha}(0)$. Then \textbf{Claim I} implies $u\in
C^{k_0+1,\alpha}(0)$. Therefore, \Cref{th1.3} holds for $k=k_0+1$. By induction, the proof of
\Cref{th1.3} is completed. Thus, in this section, we only need to prove \textbf{Claim I} and assume that
$u\in C^{k-1,\alpha}(0)$ throughout this section.

We first prove a special result.
\begin{lemma}\label{le4.0}
Let $0<\alpha<1$ and $u\in C(\bar{B}_1)$ be a viscosity solution of
\begin{equation*}
  F(D^2u,Du,u,x)=f~~\mbox{ in}~B_1,
\end{equation*}
where $F$ is $\rho$-uniformly elliptic. Suppose that \cref{e.st.1} and \cref{e1.2} hold (with $\delta_0$) and
\begin{equation*}
u(0)=\cdots=|D^{k-1}u(0)|=0,\quad
  \|u\|_{L^{\infty}(B_1)}\leq \delta_0,\quad |f(x)|\leq \delta_0|x|^{k-2+\alpha},~\forall ~x\in B_1,
\end{equation*}
where $\delta_0>0$ depends only on $n,\lambda,\Lambda,\rho,b_0,c_0,\alpha,k$ and $K_F$. Assume for
some $\bar P\in \mathcal{HP}_{k}$,
\begin{equation}\label{e4.14}
\|\bar P\|\leq \delta_0,\quad  |F_0(D^2\bar P(x),D\bar P(x),\bar P(x),x)|\leq \bar C|x|^{k-1},~\forall ~x\in B_1,
\end{equation}
where $\bar{C}$ depending only on $n,\lambda,\Lambda,\rho,b_0,c_0,\alpha,k$ and $K_F$, is to specified
in \Cref{le4.1}.

Then $u\in C^{k,\alpha}(0)$. That is, there exists $P\in \mathcal{HP}_{k}$ such that
\begin{equation*}
|u(x)-P(x)|\leq C|x|^{k+\alpha},~\forall ~x\in B_1
\end{equation*}
and
\begin{equation*}
\|P\|\leq C\delta_0,\quad  |F_0(D^2P(x),DP(x),P(x),x)|\leq \bar C|x|^{k-1},~\forall ~x\in B_1,
\end{equation*}
where $C$ depends only on $n,\lambda,\Lambda,\rho,b_0,c_0,\alpha,k$ and $K_F$.
\end{lemma}

\begin{remark}\label{re4.1}
Since we have assumed $u\in C^{k-1,\alpha}(0)$, $D^{i}u(0)$ ($1\leq i\leq k-1$) are well defined.
\end{remark}
~\\

To prove \Cref{le4.0}, we prove the following key step by the compactness method as before.
\begin{lemma}\label{le4.1}
Let $0<\alpha<1$ and $u\in C(\bar{B}_1)$ be a viscosity solution of
\begin{equation*}
  F(D^2u,Du,u,x)=f~~\mbox{ in}~B_1,
\end{equation*}
where $F$ is $\rho$-uniformly elliptic. Suppose that \cref{e.st.1} holds. Let $r\leq r_0$ and assume that for
some $P_0\in \mathcal{HP}_{k}$, \cref{e1.2} holds (with $\delta_1$),
\begin{equation*}
u(0)=\cdots=|D^{k-1}u(0)|=0,\quad
\|u-P_0\|_{L^{\infty}(B_{r})}\leq r^{k+\alpha},\quad
\end{equation*}
\begin{equation*}
|f(x)|\leq \delta_1|x|^{k-2+\alpha},~\forall ~x\in B_1
\end{equation*}
and
\begin{equation*}
\|P_0\|\leq \bar{C}r_0^{k+\alpha},\quad   |F_0(D^2P_0(x),DP_0(x),P_0(x),x)|\leq \bar{C}|x|^{k-1},~\forall
~x\in B_1,~
\end{equation*}
where $\bar{C}$ depends only on $n,\lambda,\Lambda,\rho,b_0,c_0,\alpha$ and $k$, and
$0<r_0,\delta_1<1$ (small) depend also on $K_F$.

Then there exists $P\in \mathcal{HP}_{k}$ such that
\begin{equation*}
  \|u-P\|_{L^{\infty}(B_{\eta r})}\leq  (\eta r) ^{k+\alpha}
\end{equation*}
and
\begin{equation*}
\|P-P_0\|_r\leq \bar{C}(\eta r)^{2+\alpha},\quad  |F_0(D^2P(x),DP(x),P(x),x)|\leq \bar{C}|x|^{k-1},\quad
\end{equation*}
where $0<\eta<1/2$ depends only on $n,\lambda,\Lambda,\rho,b_0,c_0,\alpha$ and $k$.
\end{lemma}
\proof We prove the lemma by contradiction. Suppose that the conclusion is false. Then there exist sequences
of $F_m,F_{0m},u_m,f_m,P_m,r_m$ such that $r_m\leq 1/m$,
\begin{equation*}
  F_m(D^2u_m,Du_m,u_m,x)=f_m~~\mbox{ in}~B_{r_m},
\end{equation*}
\begin{equation*}
u_m(0)=\cdots=|D^{k-1}u_m(0)|=0,\quad
\|u_m-P_m\|_{L^{\infty}(B_{r_m})}\leq r_m^{k+\alpha},
\end{equation*}
\begin{equation*}
|f_m(x)|\leq \frac{1}{m}|x|^{k-2+\alpha},~\forall ~x\in B_1
\end{equation*}
and
\begin{equation}\label{e4.3}
\|P_m\|\leq \frac{\bar{C}}{m},\quad   |F_{0m}(D^2P_m(x),DP_m(x),P_m(x),x)|\leq
\frac{1}{m}|x|^{k-1},~\forall ~x\in B_1,~
\end{equation}
where $F_m,F_{0m}$ are $\rho$-uniformly elliptic and \cref{e1.2} (with $\delta$ replaced by $1/m$) holds
for $F_m$ and $F_{0m}$. Furthermore, $F_{0m}\in C^{k-1}$ (with the same bound $K_F$). Moreover, for
any $P\in \mathcal{HP}_{k}$ with
\begin{equation}\label{e4.7}
|F_{0m}(D^2P(x),DP(x),P(x),x)|\leq \bar{C}|x|^{k-1},\quad \|P-P_m\|\leq \bar{C}(\eta r_m)^{k+\alpha},
\end{equation}
we have
\begin{equation}\label{e4.1}
\|u_m-P\|_{L^{\infty}(B_{\eta r_m})}> (\eta r_m)^{k+\alpha},
\end{equation}
where $\bar{C}$ and $0<\eta<1/2$ are to be specified later.

Let
\begin{equation*}
\tilde{x}=\frac{x}{r_m},\quad \tilde{u}_m(\tilde{x})=\frac{u_m(x)-P_m(x)}{r_m^{k+\alpha}}.
\end{equation*}
As before, $\tilde{u}_m$ are viscosity solutions of
\begin{equation}\label{e4.2}
    \tilde F_m(D^2\tilde u_m,D\tilde u_m,\tilde u_m,\tilde{x})=\tilde f_m~~\mbox{ in}~B_1,
\end{equation}
where
\begin{equation*}
  \begin{aligned}
\tilde{F}_m&(M,p,s,\tilde{x})\\
=&r_m^{-(k-2+\alpha)}
F_m(r_m^{k-2+\alpha}M+D^2P_m(x),r_m^{k-1+\alpha}p+DP_m(x),r_m^{k+\alpha}s+P_m(x),x)\\
&-r_m^{-(k-2+\alpha)}F_m(D^2P_m(x),DP_m(x),P_m(x),x),~\\
\tilde{f}_m&(\tilde{x})=r_m^{-(k-2+\alpha)}f_m(x)
-r_m^{-(k-2+\alpha)}F_m(D^2P_m(x),DP_m(x),P_m(x),x).\\
  \end{aligned}
\end{equation*}

Next, by \Cref{pr10.1}, for $m$ large enough,
\begin{equation*}
  \tilde u_m\in S^*_{\frac{1}{2}r_m^{-(k-2+\alpha)}\rho}(\lambda,\Lambda,r_mb_0,\bar{f}_m),
\end{equation*}
where
\begin{equation*}
\bar{f}_m(\tilde{x})=|r_m^{-(k-2+\alpha)}
\left(f_m(x)-F_m(D^2P_m,DP_m(x),P_m(x),x)\right)|+r_m^2c_0|u_m(x)|.
\end{equation*}
By the assumptions,
\begin{equation*}
  \begin{aligned}
\|\bar{f}_m\|_{L^{\infty}(B_1)}\leq& \frac{1}{m}+r_m^{1-\alpha}+c_0r_m^{k+2+\alpha}.
  \end{aligned}
\end{equation*}
Hence, for any $\varepsilon>0$, we can take $m_0$ large enough such that for any $m\geq m_0$,
\begin{equation*}
  \|\bar{f}_m\|_{L^{\infty}(B_1)}\leq \varepsilon_1,\quad
  \tilde C \left(\frac{4\rho_3}{r_m^{-(k-2+\alpha)}\rho}\right)^{\alpha_0/2}\leq \frac{\varepsilon}{2},
\end{equation*}
where $\alpha_0,\tilde C,\rho_3$ and $\varepsilon_1$ are as in \Cref{co10.2}. Set
\begin{equation*}
\delta=\left(\frac{4\rho_3}{r_m^{-(k-2+\alpha)}\rho}\right)^{1/2}.
\end{equation*}
By \Cref{co10.2}, for any $x_1,x_2\in B_{1/2}$ with $|x_1-x_2|\leq \delta$, by choosing $x_3\in B_{1/2}$
with $|x_1-x_3|=|x_2-x_3|=\delta$, we have for any $m\geq m_0$,
\begin{equation*}
|\tilde{u}_m(x_1)-\tilde{u}_m(x_2)|\leq |\tilde{u}_m(x_1)-\tilde{u}_m(x_3)|
+|\tilde{u}_m(x_2)-\tilde{u}_m(x_3)|\leq 2\tilde C\delta^{\alpha_0}\leq \varepsilon.
\end{equation*}
Thus, $\tilde{u}_m$ are equicontinuous. Then there exists $\tilde{u}\in C(\bar{B}_{1/2})$ such that
$\tilde{u}_m\to \tilde{u}$ in $L^{\infty}(B_{1/2})$.

As before, there exists a constant symmetric matrix $A$ such that
\begin{equation*}
D_MF_{0m}(0,0,0,0)\to A.
\end{equation*}
Now, we show that $\tilde{u}$ is a viscosity solution of
\begin{equation}\label{e4.4}
A^{ij}\tilde{u}_{ij}=0~~\mbox{ in}~B_{1/2}.
\end{equation}

Given $\tilde x_0\in B_{1/2}$ and $\varphi\in C^2$ touching $\tilde{u}$ strictly by above at $\tilde x_0$.
Then there exist a sequence of $\tilde x_m\to \tilde x_0$ such that $\varphi+c_m$ touch $\tilde{u}_m$ by
above at $\tilde x_m$ and $c_m\to 0$. By the definition of viscosity solution, for $m$ large enough,
\begin{equation*}
  \tilde{F}_m(D^2\varphi(\tilde x_m),D\varphi(\tilde x_m),\varphi(\tilde x_m)+c_m, \tilde x_m)\geq
  \tilde{f}_m(\tilde x_m).
\end{equation*}
Note that
\begin{equation*}
  \begin{aligned}
\tilde{f}_m&(\tilde x_m)\leq \tilde{F}_m(D^2\varphi(\tilde x_m),D\varphi(\tilde x_m),\varphi(\tilde
x_m)+c_m, \tilde x_m)\\
=&\frac{1}{r_m^{k-2+\alpha}}F_m(r_m^{k-2+\alpha}D^2\varphi+D^2P_m,
r_m^{k-1+\alpha}D\varphi+DP_m,r_m^{k+\alpha}(\varphi+c_m)+P_m,r_m\tilde x_m)\\
&-\frac{1}{r_m^{k-2+\alpha}}F_{0m}(r_m^{k-2+\alpha}D^2\varphi+D^2P_m,
r_m^{k-1+\alpha}D\varphi+DP_m,r_m^{k+\alpha}(\varphi+c_m)+P_m,r_m\tilde x)\\
&+\frac{1}{r_m^{k-2+\alpha}}F_{0m}(r_m^{k-2+\alpha}D^2\varphi+D^2P_m,
r_m^{k-1+\alpha}D\varphi+DP_m,r_m^{k+\alpha}(\varphi+c_m)+P_m,r_m\tilde x)\\
&-\frac{1}{r_m^{k-2+\alpha}}F_{0m}(D^2P_m,DP_m,P_m,r_m\tilde x_m)\\
\leq&\frac{1}{m}+F_{0m,M_{ij}}(\xi)\varphi_{ij}(\tilde x_m)+r_mF_{0m,p_i}(\xi)\varphi_{i}(\tilde x_m)
+r_m^2F_{0m,s}(\xi)\varphi(\tilde x_m),
  \end{aligned}
\end{equation*}
where
\begin{equation*}
\xi=(\theta r_m^{k-2+\alpha}D^2\varphi+D^2P_m,\theta r_m^{k-1+\alpha}D\varphi+DP_m,\theta
r_m^{k+\alpha}(\varphi+c_m)+P_m,r_m\tilde x_m)
\end{equation*}
for some $0<\theta<1$. Let $m\to \infty$, we have
\begin{equation*}
A^{ij}\varphi_{ij}(\tilde x_0)\geq 0.
\end{equation*}
Hence, $\tilde{u}$ is a subsolution of \cref{e4.4}. Similarly, we can prove that it is a viscosity supersolution
as well. That is, $\tilde{u}$ is a viscosity solution.

Next, we show $\tilde u_m\in C^{k-1,\alpha}(0)$ and their norms have a uniform bound for $m$ large
enough. Take $\delta$ small (to be specified later) and set $\hat u_m=\delta \tilde{u}_m$. Then $\hat u_m$
are viscosity solutions of
\begin{equation}\label{e4.12}
 \hat F_m(D^2\hat u_m,D\hat u_m,\hat u_m,\tilde{x})=\hat f_m~~\mbox{ in}~B_1,
\end{equation}
where
\begin{equation*}
  \begin{aligned}
\hat{F}_m(M,p,s,\tilde{x})=\delta\tilde{F}_m(\delta^{-1}M,\delta^{-1}p,\delta^{-1}s,\tilde{x}),\quad
\hat{f}_m(\tilde{x})=\delta\tilde{f}_m(\tilde{x})
  \end{aligned}
\end{equation*}
and we define $\hat{F}_{0m}$ similarly.

Since $r_m\to 0$ and $\|P_m\|\to 0$, for $m$ large enough (similarly in the following argument),
$\hat{F}_m$ are $\rho$-uniformly elliptic (although $\delta$ is very small) and
\begin{equation*}
\|\hat{u}_m\|_{L^{\infty}(B_1)}\leq \delta,\quad \|\hat{f}_m\|_{C^{k-2+\alpha}(0)}\leq \delta
\end{equation*}
In addition, $\hat{F}_m(0,0,0,\tilde x)\equiv 0$ obviously and it can be verified easily that $\hat{F}_m$
satisfies \cref{e.st.1} with $b_0$ and $c_0$. Next, for any $|M|,|p|,|s|\leq \rho$ and $\tilde x\in B_1$,
\begin{equation*}
  \begin{aligned}
|\hat F_m&(M,p,s,\tilde x)-\hat F_{0m}(M,p,s,\tilde x)|
&\leq \frac{2}{m}|\tilde x|^{k-2+\alpha}\leq \delta |\tilde x|^{k-2+\alpha}.
  \end{aligned}
\end{equation*}

Finally, we show
\begin{equation}\label{e4.13}
\|\hat F_{0m}\|_{C^{k-1}(\bar{\textbf{B}}_{\rho}\times \bar{B}_1)}\leq \hat K_F,
\end{equation}
where $\hat{K}_F$ depends only on $n,\rho$ and $K_F$.
In fact, note that
\begin{equation*}
  \begin{aligned}
\hat{F}_{0m}(M,p,s,\tilde{x})&=\hat{F}_{0m}(M,p,s,\tilde{x})-\hat{F}_{0m}(0,0,0,\tilde{x})\\
    &=\int_{0}^{1} F_{0m,M_{ij}}(\xi)M_{ij}+r_mF_{0m,p_i}(\xi)p_i+r_m^2F_{0m,s}(\xi)sdt,\\
  \end{aligned}
\end{equation*}
where
\begin{equation*}
\xi=t\left(\delta^{-1}r_m^{k-2+\alpha}M,\delta^{-1}r_m^{k-1+\alpha}p,
\delta^{-1}r_m^{k+\alpha}s,x\right)
+\left(D^2P_{m}(x),DP_{m}(x),P_{m}(x),x\right).
\end{equation*}
Hence,
\begin{equation*}
\|\hat F_{0m}\|_{C^{k-2}(\bar{\textbf{B}}_{\rho}\times \bar{B}_1)}\leq \hat K_F.
\end{equation*}
Next, from the definition of $\hat{F}_{0m}$, if any $(k-1)$-th derivative of $\hat{F}_{0m}$ involves one
derivative with respect to $M,p$ or $s$, it is bounded; if we take $(k-1)$-th derivative with respect to
$\tilde{x}$,
\begin{equation*}
  \begin{aligned}
D^{k-1}_{\tilde{x}}&\hat{F}_{0m}(M,p,s,\tilde{x})\\
=&r_m^{1-\alpha}D^{k-1}_{x}\Big(F_{0m}(\delta^{-1}r_m^{k-2+\alpha}M+D^2P_{m},
\delta^{-1} r_m^{k-1+\alpha}p+DP_{m},\delta^{-1}r_m^{k+\alpha}s+P_{m},x)\Big)\\
\leq& \hat K_F, ~\forall ~(M,p,s,x)\in \bar{\textbf{B}}_\rho\times \bar B_1.
  \end{aligned}
\end{equation*}
Therefore, \cref{e4.13} holds.

Choose $\delta$ small enough such that \Cref{th1.3} holds for $\delta$ and $\hat{K}_F$. By induction,
$\hat{u}_m\in C^{k-1,\alpha}(0)$ and hence $\tilde{u}_m\in C^{k-1,\alpha}(0)$ and their norms have a
uniform bound for $m$ large enough.

Since $u_m(0)=\cdots =D^{k-1}u_m(0)=0$ and $P_m\in \mathcal{HP}_{k-1}$,
\begin{equation*}
|\tilde u_m(\tilde x)|\leq C|\tilde x|^{k-1+\alpha},~\forall ~\tilde x\in B_1,
\end{equation*}
where $C$ is a constant independent of $m$. By taking $m\to \infty$,
\begin{equation}\label{e4.10}
|\tilde u(\tilde x)|\leq C|\tilde x|^{k-1+\alpha},~\forall ~\tilde x\in B_1.
\end{equation}

Since \cref{e4.4} is a linear equation with constant coefficients, $\tilde{u}\in C^{\infty}(B_{1/2})$. By
noting \cref{e4.10}, there exists $\tilde{P}\in \mathcal{HP}_{k}$ such that for any $0<\eta<1/4$,
\begin{equation*}
  \|\tilde{u}-\tilde P\|_{L^{\infty}(B_{\eta})}\leq C_1\eta^{k+1}\|\tilde{u}\|_{L^{\infty}(B_{1/2})}\leq
  C_1\eta^{k+1},
\end{equation*}
\begin{equation*}
A^{ij}\tilde P_{ij}\equiv 0
\end{equation*}
and
\begin{equation*}
  \|\tilde P\|\leq C_2,\quad \|\tilde{u}\|_{L^{\infty}(B_{1/2})}\leq C_2,
\end{equation*}
where $C_1$ and $C_2$ depend only on $n,\lambda,\Lambda$ and $k$. By taking $\eta$ small and
$\bar{C}$ large such that
\begin{equation*}
C\eta^{1-\alpha}\leq \frac{1}{2},\quad C_2\leq (\bar{C}-1)\eta^{k+\alpha}.
\end{equation*}
Then
\begin{equation}\label{e4.6}
  \|\tilde{u}-\tilde P\|_{L^{\infty}(B_{\eta})}\leq \frac{1}{2}\eta^{k+\alpha},\quad
   \|\tilde P\|\leq (\bar{C}-1)\eta^{k+\alpha}.
\end{equation}

Now, we try to construct a sequence of $\tilde{P}_m\in \mathcal{HP}_{k}$ such that
$\tilde{P}_m\rightarrow \tilde{P}$ as $m\rightarrow \infty$ and
\begin{equation}\label{e4.11}
D^iG_m(0)=0,~\forall ~1\leq i\leq k-2,
\end{equation}
where
\begin{equation*}
G_m(\tilde x)=\tilde F_{0m}(D^2\tilde P_m(\tilde x),D\tilde P_m(\tilde x),\tilde P_m(\tilde x),\tilde{x}).
\end{equation*}
The \cref{e4.11} implies
\begin{equation*}
D^i \left(F_{0m}(D^2Q(x),DQ(x),Q(x),x)\right)\Big|_{x=0}=0,~\forall ~1\leq i\leq k-2,
\end{equation*}
where
\begin{equation*}
Q_m(x)=P_m(x)+r_m^{k+\alpha}\tilde{P}_m(\tilde{x}).
\end{equation*}
Then \cref{e4.7} holds since $F_{0m}\in C^{k-1}$.

To prove \cref{e4.11}, let $\hat{P}_m\in \mathcal{HP}_{k}$ ($m\geq 1$) with $\|\hat{P}_m\|\leq 1$ to be
specified later and
\begin{equation*}
\tilde{P}_m(\tilde{x})=\hat{P}_m(\tilde{x})+\tilde{P}(\tilde{x}).
\end{equation*}
Since $\tilde{P}_m\in \mathcal{HP}_{k}$ and $\tilde F_{0m}(0,0,0,\tilde x)\equiv 0$,
\begin{equation*}
  \begin{aligned}
|G_m(\tilde x)|=& |\tilde F_{0m}(D^2\tilde{P}_m(\tilde x),D\tilde{P}_m(\tilde x),\tilde{P}_m(\tilde
x),\tilde x)-\tilde F_{m}(0,0,0,\tilde x)|\\
=& |\int_{0}^{1} \tilde F_{0m,M_{ij}}(\xi)\tilde{P}_{m,ij}+\tilde F_{0m,p_{i}}(\xi)\tilde{P}_{m,i}
+\tilde F_{0m,s}(\xi)\tilde{P}_{m}dt|\\
\leq& C|\tilde x|^{k-2},
  \end{aligned}
\end{equation*}
where
\begin{equation*}
\xi=(tD^2\tilde{P}_m(x),tD\tilde{P}_m(\tilde x),t\tilde{P}_m(\tilde x),\tilde x).
\end{equation*}
Hence, to verify \cref{e4.11} for $\tilde{P}_m$, we only need to prove $D^{k-2}G_m(0)=0$. Indeed, since
$\tilde{P}_m\in \mathcal{HP}_{k}$,
\begin{equation*}
  \begin{aligned}
D^{k-2}G_m(0)&=\tilde F_{0m,M_{ij}}(0,0,0,0)D^{k-2}\tilde{P}_{m,ij}\\
&=\tilde F_{0m,M_{ij}}(0,0,0,0)D^{k-2}\tilde{P}_{ij}+\tilde F_{0m,M_{ij}}(0,0,0,0)D^{k-2}\hat
P_{m,ij},
  \end{aligned}
\end{equation*}
Since
\begin{equation*}
  \begin{aligned}
\tilde F_{0m,M_{ij}}(0,0,0,0)D^{k-2}\tilde{P}_{ij}\rightarrow \tilde A^{ij}D^{k-2}\tilde{P}_{ij}=0
  \end{aligned}
\end{equation*}
and $\lambda\leq \tilde F_{0m,M_{ij}}(0)\leq \Lambda$, we can choose proper $\hat P_m$ such that
\begin{equation*}
D^{k-2}G_m(0)=0,~\forall ~m\geq 1~~\mbox{ and}~\|\hat P_m\|\rightarrow 0.
\end{equation*}
Hence, \cref{e4.7} holds for $Q_m$. Then,
\begin{equation*}
\|u_m-Q_m\|_{L^{\infty}(B_{\eta}r_m)}> (\eta r_m)^{k+\alpha}.
\end{equation*}
Equivalently,
\begin{equation*}
\|\tilde u_m-\tilde P-\hat{P}_m\|_{L^{\infty}(B_{\eta})}> \eta^{k+\alpha}.
\end{equation*}
Let $m\to \infty$, we have
\begin{equation*}
\|\tilde u-\tilde P\|_{L^{\infty}(B_{\eta})}\geq \eta^{k+\alpha},
\end{equation*}
which contradicts with \cref{e4.1}.~\qed~\\

Now, we give the~\\
\noindent\textbf{Proof of \Cref{le4.0}.} To prove \Cref{le4.0}, we only to prove the following. There exist a
sequence of $P_m\in \mathcal{HP}_{k}$ ($m\geq -1$) such that for all $m\geq 0$,
\begin{equation}\label{e4.8}
\|u-P_m\|_{L^{\infty }(B _{\eta^{m}r_0})}\leq (\eta ^{m}r_0)^{k+\alpha},
\end{equation}
\begin{equation}\label{e4.9-1}
|F_0(D^2P_m(x),DP_m(x),P_m(x),x)|\leq \hat{C}|x|^{k-1},~\forall ~x\in B_1,
\end{equation}
and
\begin{equation}\label{e4.9}
\|P_m-P_{m-1}\|_{\eta^{m}r_0}\leq \bar{C}(\eta ^{m}r_0)^{k+\alpha},
\end{equation}
where $\eta,r_0$ and $\bar{C}$ are as in \Cref{le4.1}.

We prove above by induction. Set $P_0=\bar{P}$ and $P_{-1}\equiv 0$. Choose $\delta_0$ such that
\begin{equation*}
\delta_0=\frac{1}{2} r_0^{k+\alpha}.
\end{equation*}
Then by \cref{e4.14} and noting
\begin{equation*}
\|u-P_0\|_{L^{\infty }(B_{r_0})}\leq \|u\|_{L^{\infty }(B_{r_0})}+\|P_0\|_{L^{\infty }(B_{r_0})}\leq
2\delta_0\leq r_0^{k+\alpha},
\end{equation*}
\crefrange{e4.8}{e4.9} hold for $m=0$. Suppose that the conclusion holds for $m\leq m_0$. By
\cref{e4.9},
\begin{equation*}
\|P_{m_0}\|\leq \sum_{i=1}^{m_0}\|P_i-P_{i-1}\|+\|P_0\|\leq
\bar{C}r_0^{k+\alpha}\frac{\eta^{k+\alpha}}{1-\eta^{k+\alpha}}+\frac{1}{2}r_0^{k+\alpha}\leq
\bar{C}r_0^{k+\alpha}.
\end{equation*}
By \Cref{le4.1}, the conclusion holds for $m=m_0+1$. By induction, the proof of \Cref{le4.0} is completed.
~\qed~\\

Now, we give the~\\
\noindent\textbf{Proof of \Cref{th1.3}.} As explained at the beginning of this section, we only need to prove
\textbf{Claim I}. That is, we assume $u\in C^{k-1,\alpha}(0)$ and need to prove $u\in C^{k,\alpha}(0)$.

Let $F_1(M,p,s,x)=F(M,p,s,x)-P_f(x)$ for $(M,p,s,x)\in \mathcal{S}^n\times \mathbb{R}^n\times
\mathbb{R}\times \bar{B}_1$. Then $u$ satisfies
\begin{equation*}
F_1(D^2u,Du,u,x)=f_1 ~~\mbox{in}~~B_1,
\end{equation*}
where $f_1(x)=f(x)-P_f(x)$. Thus,
\begin{equation*}
  |f_1(x)|\leq [f]_{C^{k-2,\alpha}(0)}|x|^{k-2+\alpha}\leq \delta|x|^{k-2+\alpha}, ~~\forall ~x\in B_1.
\end{equation*}
Define
\begin{equation*}
P_u(x)=\sum_{|\sigma|\leq k-1} \frac{1}{\sigma !} D^{\sigma} u(0) x^{\sigma}.
\end{equation*}
Set $u_1=u-P_u$ and $F_2(M,p,s,x)=F_2(M+D^2P_u(x),p+DP_u(x),s+P_u(x),x)$. Then $u_1$ satisfies
\begin{equation}\label{e4.15}
F_2(D^2u_1,Du_1,u_1,x)=f_1 ~~\mbox{in}~~B_1
\end{equation}
and
\begin{equation*}
  u_1(0)=|Du_1(0)|=\cdots=|D^{k-1}u_1(0)|=0.
\end{equation*}
We define $F_{01},F_{02}$ in a similar way.

Next, we try to construct $\bar{P}\in \mathcal{HP}_{k}$ such that \cref{e4.14} holds. Since we have
known the $C^{k-1,\alpha}(0)$ regularity,
\begin{equation*}
|F_{02}(0,0,0,x)|=|F_0(D^2P_u(x),DP_u(x),P_u(x),x)-P_f(x)|\leq \bar{C}|x|^{k-2},~\forall ~x\in B_1,
\end{equation*}
which implies
\begin{equation*}
D^{i}\left(F_{02}(0,0,0,x)\right)\Big|_{x=0}=0,~\forall ~0\leq i\leq k-3.
\end{equation*}
Since $\bar{P}\in \mathcal{HP}_{k}$,
\begin{equation*}
D^{i}\left(F_2(D^2\bar{P},D\bar{P},\bar{P},x)\right)\Big|_{x=0}
=D^{i}\left(F_2(0,0,0,x)\right)\Big|_{x=0}=0,~\forall ~0\leq i\leq k-3.
\end{equation*}

Next, we compute the $k-2$ order derivatives. Since $F_0(0,0,0,x)\equiv 0$,
\begin{equation*}
  \begin{aligned}
&\left|D^{k-2}_xF_0(D^2P_u(0),DP_u(0),P_u(0),0)\right|\\
&=\left|D^{k-2}_x F_0(D^2P_u(0),DP_u(0),P_u(0),0)-D^{k-2}_x F_0(0,0,0,0)\right|\\
 &\leq K_F\|P_u\|\leq C\delta,\\
  \end{aligned}
\end{equation*}
which implies
\begin{equation*}
\left|D^{k-2}\left(F_2(0,0,0,x)\right)\Big|_{x=0}\right|\leq C\delta.
\end{equation*}
Note
\begin{equation*}
  \begin{aligned}
D^{k-2}&\left(F_2(D^2\bar{P},D\bar{P},\bar{P},x)\right)\Big|_{x=0}
=F_{2,M_{ij}}(0,0,0,0)D^{k-2}\bar{P}_{ij}+D^{k-2}\left(F_2(0,0,0,x)\right)\Big|_{x=0}\\
  \end{aligned}
\end{equation*}
and $\lambda I\leq D_MF_2\leq \Lambda I$. Then we can choose proper $\bar{P}$ such that
\begin{equation*}
D^{k-2}\left(F_2(D^2\bar{P},D\bar{P},\bar{P},x)\right)\Big|_{x=0},\quad \|\bar{P}\|\leq C\delta,
\end{equation*}
where $C$ depends only on $n,\lambda, \Lambda,\rho,k$ and $K_F$. Therefore, \cref{e4.14} holds by
taking $\delta$ small enough. By \Cref{le4.0}, $u_1$ and hence $u$ is $C^{k,\alpha}$ at $0$. \qed~\\

Finally, we give the~\\
\noindent\textbf{Proof of \Cref{co1.0}.} By \Cref{th1.2}, $u\in C^{2,\alpha}(0)$ and there exists $P\in\mathcal{P}_2$ such that
\begin{equation*}
|u(x)-P(x)|\leq C|x|^{2+\alpha},~\forall ~x\in B_1
\end{equation*}
and
\begin{equation}\label{e10.6}
\|P\|\leq \bar C\delta,\quad  F(D^2P,DP(0),P(0),0)=f(0).
\end{equation}

For $r>0$ to be specified later, let
\begin{equation*}
y=\frac{x}{r}, \quad  \tilde{u}(y)=\frac{u(x)-P(x)}{r^2}, \quad \tilde{f}(y)=f(x)-F(D^2P,DP(x),P(x),x).
\end{equation*}
Then $\tilde{u}$ is a solution of
\begin{equation}\label{e10.7}
\tilde{F}(D^2\tilde{u},D\tilde{u},\tilde{u},x)=\tilde{f}\quad\mbox{in}~B_1,
\end{equation}
where
\begin{equation*}
\tilde{F}(M,p,s,y)=F(M+D^2P,rp+DP(x),r^2s+P(x),x)-F(D^2P,DP(x),P(x),x).
\end{equation*}
In addition, define
\begin{equation*}
\tilde{F}_0(M,p,s,y)=F_0(M+D^2P,rp+DP(x),r^2s+P(x),x)-F_0(D^2P,DP(x),P(x),x).
\end{equation*}

Next, we show that \cref{e10.7} satisfies the assumptions of \Cref{th1.3}. First, $\tilde{F}_0(0,0,0,y)\equiv 0$ clearly. In addition, by \cref{e10.6} and taking $\delta$ small enough, $\tilde{F}$ and $\tilde{F}_0$ are $\rho/2$-uniformly elliptic and
\begin{equation*}
K_{\tilde{F}}=\|\tilde F_{0}\|_{C^{k-1}(\bar{\textbf{B}}_{\rho/2}\times \bar{B}_1)}\leq CK_F.
\end{equation*}
Moreover, it can be verified easily that for any $|M|,|p|,|s|\leq \rho/2,~y\in B_1$,
\begin{equation*}
|\tilde F(M,p,s,y)-\tilde F_0(M,p,s,y)|\leq C\delta |x|^{k-2+\alpha}
\leq C\delta r^{k-2+\alpha}|y|^{k-2+\alpha}.
\end{equation*}
Finally, by noting $\tilde{f}(0)=0$,
\begin{equation*}
  \|\tilde u\|_{L^{\infty}(B_1)}\leq Cr^{\alpha},\quad \|\tilde f\|_{C^{k-2+\alpha}(0)}\leq Cr.
\end{equation*}

Therefore, we can choose $r$ small enough such that the assumptions of \Cref{th1.3} are satisfies. Then $\tilde{u}\in C^{k,\alpha}(0)$ and hence $u\in C^{k,\alpha}(0)$.~\qed~\\

\section*{Statements and Declarations}
\textbf{Competing Interests:} All authors declare that they have no conflicts of interest.

\printbibliography
\end{document}